\documentclass[a4paper,11pt]{amsart}

\usepackage{amssymb,amsmath, pinlabel, mathabx, amscd, tikz, amsthm, xcolor}

\usepackage{lineno}
\usepackage[normalem]{ulem} 
 \usepackage{hyperref}
\hypersetup{colorlinks=true, linkcolor=cyan,citecolor=cyan, urlcolor= cyan}
\usepackage{mathtools}

\input xy
\xyoption{all}

\usepackage[latin1]{inputenc}
\usepackage[T1]{fontenc}
\usepackage{amsfonts,amssymb}
\usepackage{amsmath}
\usepackage{pstricks}
\usepackage{amsthm}
\usepackage[all]{xy}
\usepackage{graphicx}  
\usepackage{comment} 
\usepackage{enumitem}
\usepackage{xcolor}

\usepackage{marginnote}

\definecolor{green}{rgb}{0.0, 0.5, 0.0}
\definecolor{yellow}{rgb}{0.8, 0.33, 0.0}

\setlist[itemize,enumerate]{leftmargin=*}

\newcommand{\g}{\bar{\beta}}
\newcommand{\Z}{\mathbb{Z}}
\newcommand{\R}{\mathbb{R}}
\newcommand{\Q}{\mathbb{Q}}

\renewcommand{\k}[0]{{k}}

\newcommand{\cV}{\mathcal{V}}
\newcommand{\cT}{\mathcal{T}}

\newcommand{\Trop}{\textnormal{Trop}}
\newcommand{\ord}{\mathrm{ord}}
\newcommand{\set}[1]{\left\{#1\right\}}
\renewcommand{\S}[0]{{S}}

\newcommand{\trop}{\mathrm{trop}}
\newcommand{\GF}{\mathrm{GF}}

\newcommand{\de}{\mathbf{i}}
\newcommand{\ex}{\mathbf{e}}
\newcommand{\ic}{\mathbf{c}}

\newcommand{\si}{\mathbf{s}}

\numberwithin{equation}{section}

\newtheorem{theorem}{Theorem}[section]
\newtheorem{lemma}[theorem]{Lemma}
\newtheorem{proposition}[theorem]{Proposition}
\newtheorem{corollary}[theorem]{Corollary}

\theoremstyle{definition}
\newtheorem{definition}[theorem]{Definition}
\newtheorem{notation}[theorem]{Notation}
\newtheorem{remark}[theorem]{Remark}
\newtheorem{example}[theorem]{Example}

\title{Resolving singularities of curves with one toric morphism}

\author{Ana Bel\'en de Felipe}
\address{
Departamento de Matem\'aticas, Estad\'{\i}stica e I.O.
Secci\'on de Matem\'aticas, Universidad de La Laguna. Apartado de Correos 456.
38200 La Laguna, Tenerife, Espa\~na.}
\email{afelipe@ull.edu.es}

\author{Pedro D. Gonz\'alez P\'erez} 
\address{Instituto de Matem\'atica Interdisciplinar, Departamento de \'Algebra, 
Geometr\' \i a y Topolog\'\i a,  Facultad de Ciencias Matem\'aticas,
Universidad Complutense de Madrid, Plaza de las Ciencias 3, Madrid 28040, Espa\~na.}
   \email{pgonzalez@mat.ucm.es}

\author{Hussein Mourtada} 
\address{Universit\'e de Paris, Sorbonne Universit\'e, CNRS, 
Institut Math\'ematiques de Jussieu-Paris Rive Gauche,
 F-75013, Paris, France.}
   \email{hussein.mourtada@imj-prg.fr}

 \date{\today}

\begin {document}
\maketitle 

{\em \hfill This paper is dedicated to Bernard Teissier. \smallskip} 

\begin{abstract}
We give an explicit positive answer,  in the case of reduced curve singularities, to a question of B. Teissier about the existence of a toric embedded resolution after reembedding.  In the case of a curve singularity $(C,O)$ contained in a non singular surface $S$ such a reembedding may be defined in terms of a sequence of maximal contact curves of the minimal embedded resolution of $C$. We prove that there exists a toric modification, after reembedding, which provides an embedded resolution 
of $C$. We use properties of the semivaluation space of $S$ at $O$ to describe 
how the dual graph 
of the minimal embedded resolution of $C$ may be seen on the local tropicalization of $S$ associated to this reembedding. 
\end{abstract}

%--------------------------------------------------------------------------------------------

\tableofcontents

%--------------------------------------------------------------------------------------------
\section*{Introduction}

In \cite{GT}, Goldin and Teissier proved that one can resolve the singularities of a plane branch (i.e.,\ an analytically irreducible plane curve singularity) with one toric morphism, after reembedding it in a possibly higher dimensional affine space.  
Teissier asked then the following question in \cite[Section 5]{T4}, see also \cite{T1}.

\subsubsection*{Question}
Given a reduced and equidimensional algebraic or formal
space $X$ over an algebraically closed field $\k$, is it true that for every point $x \in X$
there is a local formal embedding of $(X,x)$  into an affine space $(\mathbb{A}^m,0)$ and a toric
structure on $ \mathbb{A}^m$ such that $(X,x)\subset (\mathbb{A}^m,0)$ can be resolved by one toric morphism ?

\medskip 

This means that there exist local coordinates $u_1, \dots, u_m$ centered at $x$ and
an open (\'etale or formal) neighborhood $U$ of $x \in \mathbb{A}^m(k)$, such that 
there is a proper birational toric map $\pi: Z \to U$ with respect to the coordinates $u_1, \dots, u_m$
with $Z$ non singular and such that the strict transform $X'$ of $X \cap U$ is 
non singular and transversal to the non dense orbits at every point of $\pi^{-1} (x) \cap X'$.  An embedding satisfying the property of the question will be called {\it torific}, a terminology which combines toric and terrific.

\medskip 
We consider also the problem of existence of \textit{torific embeddings
of a triple} $(x, X,\mathbb{A}^n(k)) $ where  
$ x \in X \subset \mathbb{A}^n(k)$, and $n$ is the embedding dimension of $x \in X$, that is,  if there exists a formal embedding of  $\mathbb{A}^n(k)$ into an affine space $\mathbb{A}^m(\k)$ endowed with a toric structure, 
such that a toric morphism of $\mathbb{A}^m(\k)$ induces an embedded resolution of 
the triple $(x, X,  \mathbb{A}^n(k))$.

\medskip

Besides the case of plane branches there are some partial answers to Teissier's question:

\medskip

- Aroca, G\'omez-Morales, and Shabbir considered a notion of \textit{Newton non-degenerate} ideals 
of the ring of polynomials $\k [X_0, \dots, X_n]$ (see 
\cite{AGS}).
This is a condition on the initial ideals 
with respect to a weight vector lying in the \emph{tropicalization} of $I$. In this case, the original embedding of the singularity is already torific. 
The notion of Newton non-degenerate ideal is
related to that of sch\"on variety in \cite{Te1}. It is a generalization of the notions of  hypersurface and complete intersection singularities  which are non degenerate with respect to their Newton polyhedra, which was introduced by Khovanskii and Kouchnirenko \cite{Kou,Kho}. 
Cueto, Popescu-Pampu and Stepanov have proven that the ideals defining the natural embeddings of splice type surface singularities are Newton non-degenerate  (see \cite{CPS}). In the second version of the preprint \cite{CPS}, 
they have deduced the existence of a torific embedding 
of a reduced complex analytic plane curve singularity, as 
an application of their results on surface singularities.

- Lejeune-Jalabert and Reguera  proved in 
 \cite{LR99} that 
sandwiched surface singularities
admit natural torific embeddings, which were called \textit{toric environments}. 

 - Gonz\'alez P\'erez generalized Goldin and Teissier's result to the case of an irreducible germ of quasi-ordinary 
hypersurface singularity  (see  \cite{GP-Fourier}).

- More generally, Tevelev  proved that if 
$\k$ is an algebraically closed field of characteristic zero and $X \subset \mathbb{P}^n(k)$ 
is a projective algebraic variety, then there exists an integer $m \gg 0$ such that if 
 $X \subset \mathbb{P}^m (k)$ is the Veronese reembedding of order $m$, there exist 
projective coordinates $(z_0: \dots : z_m)$ on $\mathbb{P}^m (k)$ such that 
the intersection of $X$ with the torus $T^m (k) = (k^*)^m$ defined by this choice of coordinates, is dense in $X$
and there exists  an equivariant map of toric varieties $\pi: Z \to \mathbb{P}^m (k)$, with $Z$ non singular, such that  the strict transform of $X$ is non singular and transversal
to the non dense toric orbits in $Z$ (see \cite{Te2}). 
The proof of this result uses embedded resolution of singularities which is not known to be available in general
when the field $k$ has positive characteristic.

\medskip 
The main result of \cite{GT} motivates Teissier's strategy towards the proof of  local uniformization, which is a very
local version of resolution, 
by a method based on the comparison of a given singular germ by deformation with a space whose resolution is easy and blind to the characteristic, namely using toric methods  (see \cite{Tei09,T1, T2}).
Note here that the valuative machinery that makes things work in \cite{GT}, still works to some extent in higher dimension \cite{CMT21} but does not lead to similar conclusions on resolution of singularities.  

\medskip 

In addition, a resolution obtained by a torific embedding is often easier to use than classsical Hironaka-type resolutions for computing subtle invariants of singularities like motivic or topological zeta functions, monodromy zeta functions, log canonical thresholds and jumping numbers of  multiplier ideals (see \cite{MVV, GGGR21, GG14}). The resolutions of singularities obtained by the classical approach are often complex to handle and it is a difficult problem to link invariants of resolution of singularities like Hironaka's order of ideals to subtle  invariants of singularities such as those that we have just mentioned.

\medskip 
Let us explain the main contributions of this paper. 
Along this article, we assume that we are working over an algebraically closed field $\k$ of 
arbitrary characteristic.

\medskip 

 In Section  \ref{trc} we prove the existence of functions 
defining a torific embedding of a reduced curve singularity of arbitrary 
embedding dimension
in two different ways (see Corollary \ref{enum:diffvect} and Theorem  \ref{thm:toremrescurves}).

\medskip 

In Section \ref{resolving_reduced} we 
consider a reduced plane curve singularity $C$ at 
a point $O$.
We take  a minimal generating sequence $(x_0, \dots, x_m)$ of the divisorial valuations 
defined by those prime exceptional divisors of the minimal embedded resolution $\psi$ of $C$ 
which intersect
the strict transform of $C$.
We assume 
that the strict transform of the  branch $L_i$ defined by $x_i =0$ 
does not intersect the strict transform of $C$  by $\psi$, for $i=0, \dots, m$. 
Then, we prove that the tuple of functions $(x_0, x_1, \dots, x_m)$
defines a torific embedding of $C$ at $O$ 
(see Theorem \ref{th:resgenplane}).

\medskip 

In Section \ref{sec:TEpairs} we prove that the tuple of  functions $(x_0, \dots, x_m)$ considered above, defines
 torific embedding  of the triple $(O, C, \mathbb{A}^2(k))$,
 that is, 
 there exists a toric 
 modification of $\mathbb{A}^{m+1}(\k)$
 which induces an embedded resolution of $(O, C, \mathbb{A}^2(k))$ (see Theorem \ref{th-res}). 
 In addition, we show how 
 the dual graph of the total transform of $\cup_{j=0}^m L_j$ by the minimal embedded resolution $\psi$ of $C$ appears on 
 the support of the fan defining this modification.

 \medskip 
 
Our main result 
realises the program of the geometric approach to resolution of singularities explained in \cite{Mou20,Mo4}. This program suggests that one can find a torific embedding of a singularity $(X,x)$ from first detecting a finite number of special divisorial valuations which may be called essential and then finding an embedding of $X$ which is torific along each of these divisorial valuations. Given an essential valuation $\nu$, finding an embedding of $X$ which is torific along $\nu$ is a valuation theoretical issue; then a torific embedding may be found by a simple procedure (concatenation of the packages of functions defining the torific embeddings along each essential valuation); such a torific embedding is far from being minimal in general.  

\medskip

In this article, for a plane curve singularity, 
we find these essential valuations on the minimal embedded resolution of the curve and we characterize them by using the Eggers-Wall tree. 
We use the geometric characterization of the 
generating sequences of tuples of divisorial valuations by Spivakovsky  \cite[Section 8]{S}  in the case of 
one valuation and Delgado, Galindo, and Nu\~nez \cite{DGN} in the case of several valuations (see Proposition \ref{gen-cur}).
The embedding that we find has no superfluous information as in the concatenation procedure, it is actually minimal for the triple $(O,C,\mathbb{A}^2)$. One can also detect the essential valuations on a graph which is associated with the jet schemes of the curve singularity \cite{Mou11,LJMR}. This latter graph makes sense also for higher dimensional singularities \cite{CM21} and the \textit{loc.cit.} program suggests that one can detect essential valuations on it \cite{Mou20}. 

\medskip 

In order to prove the main results in  Section \ref{sec:TEpairs}
we consider the \emph{local tropicalization}  of $\mathbb{A}^2(k)$
associated with  the tuple of functions $(x_0, \dots, x_m)$. 
The definition and properties of local tropicalizations in a general setting 
was developped in \cite{PS13} by Popescu-Pampu and Stepanov.
Further applications of 
local tropicalization are obtained recently 
by  Cueto,  Popescu-Pampu, and Stepanov in the case of surface singularities of splice type (see \cite{CPS}).

\medskip 

The \textit{finite local tropicalization} of $\mathbb{A}^2(k)$
associated with  the tuple of functions $(x_0, \dots, x_m)$
consists of 
the tuples $(\nu (x_0), \dots,   \nu(x_m)) \in \R^m_{\geq 0}$, where $\nu$
runs though the semivaluations  of $\mathbb{A}^2(k)$
at the closed point $O$ of $C$, which take finite values on $x_0, \dots, x_m$.
We show that this set is the support of a fan $\cT$ in $\R^{m+1}$ consisting of at most two dimensional cones 
(see Theorem \ref{loc-trop}). 
For instance, if $E$ is an exceptional prime divisor in the minimal embedded resolution of
$\psi$ of $C$ and if $\nu_{E}$ denotes its associated divisorial valuation
then $( \nu_{E} (x_0), \dots, \nu_{E} (x_m)) \in \Z^{m+1}_{>0}$ belongs to the support of $\cT$.
Then, we study the initial ideals associated to the ideal ${I}$ defining 
the embedding $\mathbb{A}^2(k) \hookrightarrow \mathbb{A}^{m+1}(\k)$  with respect to weight vectors $w $ in the support of $\cT$. 
The key step is to prove that 
the ideal ${I}$ is Newton non degenerate 
(see Proposition \ref{prop:ISisNnd}). This is shown by  using the properties of expansions in terms of generating sequences. 
The main theorem  of this section
is obtained by applying the result  of \cite{AGS} mentioned above. 
We show that 
if $Y$ denotes the image of $\mathbb{A}^2(k)$ after reembedding, then 
there exists a toric modification $\pi$ of $\mathbb{A}^{m+1}(\k)$, defined by a regular fan supported on $\R_{\geq 0}^{m+1}$
which induces a regularization of the fan $\cT$,  
and such that  the restriction of $\pi$ to the strict transform $Y'$ of $Y$ is an embedded resolution of $C$
(see Theorem \ref{th-res}). If in addition the regularization of the fan $\cT$ is the minimal one, then the restriction of $\pi$ to $Y'$ is 
 the minimal  embedded resolution of $C$ (see Theorem \ref{prop:div}). 
In particular, we prove that 
the projectivization of the minimal regularization of the fan 
$\cT$ is isomorphic to the dual graph of the total transform of $\cup_{j=0}^m L_j$ by the minimal embedded resolution $\psi$ of $C$ (see Proposition \ref{prop:exc-dis}).
 
\medskip 

The description of the fan $\cT$ is very explicit in terms 
of the embedding of the \textit{Eggers-Wall tree} of $C$ with respect to the smooth branch  $L_0$ defined by $x_0 = 0$ 
in the \textit{space of normalized semivaluations} $\cV_{L_0}$ with respect to $L_0$. 
This embedding was described by Garc\'\i a Barroso, Popescu-Pampu, and Gonz\'alez P\'erez (see \cite[Section 8]{GGP19} and   \cite[Remark 5.37]{GGP19b}).
The space $\cV_{L_0}$ is one of the valuative trees considered in  Favre and Jonsson's book \cite{FJ}, 
see also the presentation of this theory in Jonsson's survey  \cite[section 7]{Jon15}.

\medskip

Besides providing an answer to  Teissier's question for plane curve singularities, our results exhibit the beautiful interplay between  semivaluation spaces, resolution of singularities, toric geometry and tropical geometry, which may shed some light on higher dimensional generalizations.

\subsubsection*{Notation}
A curve $C$ is an affine scheme $\text{Spec }\mathcal O$, where $\mathcal O$ is a complete noetherian 
local ring of dimension one over $\k$ with residue field $\k$. We call the ring $\mathcal O$ 
the local ring of $C$ and denote by $\mathcal M$ the maximal ideal of $\mathcal O$. 
The curve $C$ is reduced if the local ring $\mathcal O$ is reduced. 
If $\mathcal O$ is an integral domain, the curve $C$ is said to be a branch. 
The local ring $\mathcal O$ of a reduced curve $C$ has a finite number 
of minimal prime ideals, say $\mathcal P_j$ for $j=1,\ldots,r$, 
each one defining a branch $C_j$ with local ring $\mathcal O/\mathcal P_j$. Then we write $C=\cup_{j=1}^r{C_j}$ 
and say that $C$ is a curve with $r$ branches.
The embedding dimension of a curve $C$ is the dimension of 
$\mathcal M/\mathcal M^2$ as $\k$-vector space. A plane curve is a curve of embedding dimension at most two.

\subsubsection*{Acknowlegment} 
We are grateful to Javier Fern\'andez de Bobadilla for
sharing some ideas which were helpful in the formulation of Theorem \ref{thm:toremrescurves}.  We thank  Bernard Teissier and Alicia Dickenstein 
for useful discussions and suggestions,
 Patrick Popescu-Pampu for his comments on a preliminary version of the paper, and the referee for the careful reading and remarks.

\medskip 

The first author was supported by ERCEA Consolidator Grant 615655-NMST, the Basque Government through the BERC 2018-2021 program, BCAM Severo Ochoa excellence accreditation SEV-2017-0718, grant MTM2016-80659-P, and grant PID2019-105896GB-I00 funded by MCIN/AEI/10.13039/501100011033. The second author was supported by the Spanish grants PID2020-114750GB-C32 and MTM2016-76868-C2-1-P. The third author was supported by the French grant
Projet ANR LISA, ANR-17-CE40-0023.
We also thank the hospitality of the Institute of Mathematics of the University of Barcelona (IMUB) and the Institut de Math\'ematiques de Jussieu-Paris Rive Gauche (IMJ-PRG).

%------------------------------------------------------------------------------------------------

\section{Toric resolutions of reduced curve singularities}  \label{trc}

\subsection{Toric modifications, their critical and discriminant locus}

In this section we describe the modification induced by an ambient toric modification on an orbit closure in 
a toric variety. See the textbooks \cite{F,CLS,E,O} for standard notions about toric geometry. 

\medskip

Let $\Sigma$ be a fan with respect to a rank $n$ lattice $N$, with dual lattice $M$.
If $\sigma \subset N_\R := N \otimes {\R}$ is a cone we denote by 
$\check{\sigma} \subset M_{\R} := M \otimes {\R}$ the dual cone, and by $\sigma^\perp$ the 
orthogonal cone. 

\medskip
 
If $\sigma \in\Sigma$ then the semigroup $\check{\sigma} \cap M$ is finitely generated and the 
semigroup algebra
\[
\k [ \check{\sigma} \cap M] =
\Big\{ \sum_{finite} a_v \chi^v \mid a_v \in \k, \, v \in  \check{\sigma} \cap M \Big\}
\]
is a $k$-algebra of finite type. It is the coordinate ring of the affine toric variety $X_\sigma$. 

\medskip

In addition, if  the cone $\sigma$ is of dimension $n$, then 
there is a unique minimal system of generators $v_1, \dots, v_{m_\sigma}$  of 
the semigroup $\check{\sigma} \cap M$. Setting, $x_i := \chi^{v_i}$ for $i=1, \dots, m_{\sigma}$, 
defines an embedding of the affine toric variety 
\begin{equation} \label{emb-tor}
X_{\sigma} \hookrightarrow \k^{m_\sigma},
\end{equation} 
which sends the $0$-dimensional orbit to the origin of $\k^{m_\sigma}$.  
In this case,  the defining ideal of this $0$-dimensional orbit is 
generated by the monomials $\chi^v$, for $v \in  \check{\sigma} \cap M \setminus \{ 0\}$. 
The completion of $\k [ \check{\sigma} \cap M]$ with 
respect to this ideal is the formal power series ring 
\[
\k [[ \check{\sigma} \cap M]]=
\Big\{ \sum a_v \chi^v \mid a_v \in \k, \, v \in  \check{\sigma} \cap M \Big\}.
\]

\medskip

We denote by $X_{\Sigma, N}$ the normal toric variety associated 
with the fan $\Sigma$ and the lattice $N$, or simply by $X_{\Sigma}$, 
if 
the lattice $N$ is clear from the context. If $\sigma \in \Sigma$ we denote by $X_{\sigma}$ 
the affine open toric subvariety of $X_{\Sigma}$, and by $O(\sigma)_\Sigma$, or simply by 
$O(\sigma)$, the corresponding orbit by the action on $X_\Sigma$ of the torus $T_N=\mathrm{Hom}(M,\k^*)$. 
Abstractly, the orbit $O(\sigma)$ is the torus $\mathrm{Hom}(M(\sigma),\k^*)$, where 
$M(\sigma)$ denotes the lattice spanned by $\sigma^\perp\cap M$. In particular, we have $T_N=O(\set{0})$. 
The map which sends a cone $\sigma\in\Sigma$ to its corresponding orbit $O(\sigma)$ defines a 
bijection of $\Sigma$ with the set of orbits of the torus action on $X_\Sigma$. 

\medskip 
 
The closure of the orbit $O(\sigma)$ in $X_\Sigma$ is a normal toric variety. 
In order to describe its associated fan, we consider the lattices $N_\sigma$ spanned by $\sigma \cap N$
and the quotient $N (\sigma) := N/ N_\sigma$. The dual lattice of $N(\sigma)$ is precisely the lattice 
$M(\sigma)$ introduced before. Denote 
by $\mathrm{Star}(\sigma)_\Sigma$ the set of cones in $N(\sigma)_\R$ which are 
images by the canonical projection 
$N_\R \to N(\sigma)_\R = N_\R / (N_\sigma)_\R$
of the cones of the fan $\Sigma$ which contain $\sigma$ as a face. The set $\mathrm{Star}(\sigma)_\Sigma$ 
is a fan with respect to the lattice $N(\sigma)$.

\begin{lemma}(see \cite[Section 3.1]{F})
The closure of the orbit  $O(\sigma)$ in $X_{\Sigma}$ is the normal toric variety 
$X_{\mathrm{Star}(\sigma)_\Sigma, N(\sigma)}$
associated with the fan  $\mathrm{Star}(\sigma)_\Sigma$  and the lattice $N(\sigma)$. 
\end{lemma}

\begin{remark} \label{rem:f}
The torus $T_{N (\sigma)} = \mathrm{Hom} (M(\sigma), \k^*)$ of the toric variety 
$X_{\mathrm{Star}(\sigma)_\Sigma, N(\sigma)}$ is equal to $O(\sigma)$. 
More generally, the orbits of the action of the torus $T_{N (\sigma)}$ on 
$X_{\mathrm{Star}(\sigma)_\Sigma, N(\sigma)}$ are orbits of the action of $T_N$ on $X_\Sigma$. 
\end{remark}

\begin{notation} \label{quadrant}
We denote by $\Sigma_{0,n}$ the fan consisting of the faces of the cone 
$\sigma_{0,n}:= \R^n_{\geq 0}$, with respect to the lattice $N:=\Z^n$. 
We often write simply $\sigma_{0}:= \sigma_{0,n}$ and $\Sigma_0:=\Sigma_{0,n}$, 
if $n$ is clear from the context. We denote by $e_1, \dots, e_n$ the canonical 
basis of the lattice $N$, which spans the cone $\sigma_{0,n}$. We denote by 
$\check{e}_1, \dots,\check{e}_n$ the dual basis of $M$, which spans 
the dual cone $\check{\sigma}_{0,n}$. The toric variety $X_{\Sigma_{0,n}}$ 
is isomorphic to $\k^n$ by (\ref{emb-tor}). 
It is equipped with coordinates $(X_1, \dots, X_n) := (\chi^{\check{e}_1}, \cdots, \chi^{\check{e}_n})$. 
\end{notation}

\begin{example} \label{ex:theta}
Let us fix an integer $s$ with $1 \leq s \leq n-1$. The cone 
$\sigma := \R_{\geq 0} e_{s+1} + \cdots  + \R_{\geq 0} e_{n}$  belongs to the fan $\Sigma_{0,n}$ (see Notation 
\ref{quadrant}). The lattice $N_\sigma$ has basis $e_{s+1}, \dots,  e_{n}$. Let us denote by $\bar{u} \in N(\sigma)$ 
the image of a vector $u \in N$ by the canonical projection $N \to N(\sigma)$. Then, the lattice $N(\sigma)$ has 
basis $\bar{e}_{1}, \dots, \bar{e}_{s}$. Thus, the closure of the orbit $O(\sigma)$ is the coordinate subspace 
$\k^{s} \hookrightarrow \k^n$ with coordinates $(X_1, \dots, X_s)$, and defined by $X_{s+1}=\cdots=X_{n}=0$.
\end{example}

\medskip 

Let $\Sigma'$ be a fan subdividing $\Sigma$ with respect to the lattice $N$. We denote by $\pi_{\Sigma}^{\Sigma'}\colon X_{\Sigma'} \to X_\Sigma$ the associated toric \textit{modification}. 
The \textit{exceptional locus} of the modification $\pi_{\Sigma}^{\Sigma'}$ consists of those subvarieties of 
$X_{\Sigma'}$ which are mapped to subvarieties of $X_{\Sigma}$ of smaller dimension. The \textit{discriminant locus} 
of the modification $\pi_{\Sigma}^{\Sigma'}$ is the image of the exceptional locus. These loci can be described in terms of orbits as follows (see \cite[Section 1.1]{GP-Fourier} and also 
\cite[Page 100]{GSLJ}).

\begin{proposition} \label{prop:exc}
The exceptional locus $\mathrm{Exc} \, (\pi_{\Sigma}^{\Sigma'})$ of $\pi_{\Sigma}^{\Sigma'}$ 
is the union of orbits 
$O(\sigma)_{\Sigma'}$,  for $\sigma  \in \Sigma' \setminus \Sigma$.  
The discriminant locus of $\pi_{\Sigma}^{\Sigma'}$ is the union of orbits 
$O(\tau)_\Sigma$, for $\tau \in \Sigma \setminus \Sigma'$.
\end{proposition}

\medskip 

Let us fix a cone $\sigma \in \Sigma$. If $\sigma$ belongs also to $\Sigma'$ then the orbit 
$O(\sigma)_\Sigma \subset X_\Sigma$ is not contained in the discriminant locus of $\pi_\Sigma^{\Sigma'}$. 
We have also an orbit $O(\sigma)_{\Sigma'} \subset X_{\Sigma'}$ and the restriction of $\pi_\Sigma^{\Sigma'}$ 
defines an isomorphism of orbits $ O(\sigma)_{\Sigma'} \to O(\sigma)_\Sigma$.
The fan $\mathrm{Star}(\sigma)_{\Sigma'}$ is a subdivision of $\mathrm{Star}(\sigma)_{\Sigma}$. 
The following lemma  is a particular case of Lemma 3.3.21 of \cite{CLS}. 

\begin{lemma} \label{disc}
Let $\Sigma'$ be a subdivision of the fan $\Sigma$. If $\sigma \in \Sigma' \cap \Sigma$ 
then the restriction of $\pi_\Sigma^{\Sigma'}$ to the closure of $O(\sigma)_{\Sigma'}$ in the source and to 
the closure of $O(\sigma)_\Sigma$ in the target is the toric modification
\begin{equation} \label{f:modi}
\pi_{\mathrm{Star}(\sigma)_{\Sigma}}^{\mathrm{Star}(\sigma)_{\Sigma'}} \colon X_{\mathrm{Star}(\sigma)_{\Sigma'}} \to X_{\mathrm{Star}(\sigma)_{\Sigma}}.
\end{equation}
\end{lemma}

\begin{example} \label{uno}
Let us specialize to the case $n=3$ in Notation \ref{quadrant}. 
Consider the subdivision $\Sigma$ of $\Sigma_0$ at the vector $w= (4, 6, 13)$. It has three-dimensional cones 
$\sigma_{i, j} = \R_{\geq 0} e_i + \R_{\geq 0} e_j + \R_{\geq 0} w$, 
for $1 \leq i < j \leq 3$. The exceptional locus of $\pi_{\Sigma_0}^{\Sigma}$ is the union of 
orbits $O(\sigma)_{\Sigma}$ such that the relative interior 
of the cone $\sigma$ is contained in $\R^3_{> 0}$, while the discriminant locus of $\pi_{\Sigma_0}^{\Sigma}$ is reduced to the origin, that is, the $0$-dimensional orbit. 
Let us consider the cone $\sigma = \R_{\geq 0} e_3$. 
The closure of the orbit  $O(\sigma)_{\Sigma_0}$ is the coordinate plane of ${\k}^3$ defined 
by $X_3 = 0$. Abstractly, it is the toric variety defined by the cone 
$\R^2_{\geq 0} = \R_{\geq 0} \bar{e}_1 + \R_{\geq 0} \bar{e}_2$, 
with respect to the lattice $N(\sigma) = \Z \bar{e}_1 + \Z \bar{e}_2$, 
where we use the notation of Example \ref{ex:theta}. 
The only three-dimensional cones of $\Sigma$ which contain $\sigma$ are $\sigma_{1,3}$ and 
$\sigma_{2,3}$. These two cones intersect along the cone $\R_{\geq 0} w+ \R_{\geq 0} e_3$. 
The fan $\mathrm{Star}(\sigma)_{\Sigma}$ is the subdivision of the cone $\R^2_{\geq 0}$ along the ray spanned by $2 \bar{e}_1 + 3 \bar{e}_2$. This ray is the image 
of $\R_{\geq 0} w+ \R_{\geq 0} e_3$ under the canonical projection 
$N_\R=\R^3 \to N(\sigma)_\R=\R^2$, which sends $(a_1, a_2, a_3) \mapsto (a_1, a_2)$. 
\end{example} 

%---------------------------------------------------------------------------------------------------------------

\subsection{Sufficient conditions for the existence of a toric embedded resolution}

Let $C$ be a reduced curve, not necessarily a plane one. We denote by $\mathcal{O}$ its local ring, with maximal 
ideal $\mathcal{M}$. If $x_1, \dots, x_n$ is a set of generators of the maximal ideal $\mathcal{M}$, 
then we have a surjection
\begin{equation} \label{one}
\k [[ X_1, \dots, X_n]] \to \mathcal{O}, \mbox{ such that } X_i \mapsto x_i, \mbox{ for } i=1, \dots, n. 
\end{equation}
This surjection defines an embedding of germs $C\hookrightarrow \k^n$, which maps the special point of $C$ to 
the origin of $\k^n$. 
In this section, we give sufficient conditions on the functions $x_1, \dots, x_n$ in 
order to guarantee that an embedded resolution of  a branch $C \subset \k^n$ 
can be obtained by a toric modification of $\k^n$. 
\medskip

In what follows we consider $\k^n$ as the affine toric variety $X_{\Sigma_{0}}$ and use Notation 
\ref{quadrant}. 
The following definition is slightly more general than the one considered in 
 \cite{GT}.

\begin{definition} \label{ter}
Let $C$ be a reduced curve embedded in $\k^n$ by the map (\ref{one}). A toric modification 
$\pi_{\Sigma_0}^\Sigma \colon X_\Sigma \to X_{\Sigma_0}= \k^n$ defined by a regular subdivision 
$\Sigma$ of $\Sigma_0$ is a {\it toric embedded resolution} of $C$ if no branch of $C$ 
is contained in the discriminant locus of $\pi_{\Sigma_0}^\Sigma$, and 
the strict transform $\tilde{C}$ of $C$ is smooth and transversal
to the orbit stratification of the exceptional locus of $\pi_{\Sigma_0}^\Sigma$. 
\end{definition}
Recall that the strict transform $\tilde{C}$ of $C$ by $\pi_{\Sigma_0}^\Sigma$ is the closure in 
$X_{\Sigma}$ of 
$(\pi_{\Sigma_0}^\Sigma)^{-1} (C) \setminus  \mathrm{Exc} \,  (\pi_{\Sigma_0}^{\Sigma})$.
\begin{remark} \label{Rem: hyper}
In the paper \cite{GT} the case when the curve $C$ is contained in  a coordinate hyperplane
was not considered. In order to deal with this case we introduced the hypothesis 
about the discriminant locus of the modification in Definition \ref{ter}.
The transversality condition in Definition \ref{ter} should be considered with respect to the modifications of the minimal 
coordinate subspaces containing the branches of $C$, which are toric morphisms by 
Lemma \ref{disc}. 
Notice that the minimal coordinate subspace of $X_{\Sigma_0}=\k^n$ which 
contains a branch $C_j$ of $C$ 
is the closure $X_{\mathrm{Star}(\sigma)_{\Sigma_0}}$
of an orbit $O(\sigma)_{\Sigma_0}$, for some $\sigma \in \Sigma_0$ (see Example \ref{ex:theta}). By Proposition \ref{prop:exc} this coordinate subspace is not contained in 
the discriminant locus of $\pi_{\Sigma_0}^{\Sigma}$
if and only if $\sigma \in \Sigma$. Then, the toric modification
$\pi_{\Sigma_0}^{\Sigma}$ is a toric embedded resolution of $C_j$ according to
Definition \ref{ter}. This means that 
the map
$\pi_{\mathrm{Star}(\sigma)_{\Sigma_0}}^{\mathrm{Star}(\sigma)_{\Sigma}}$, defined in \eqref{f:modi},
is a toric embedded resolution of 
$C_j \subset X_{\mathrm{Star}(\sigma)_{\Sigma_0}}$
according to the definition considered in \cite{GT}.
\end{remark}

Before stating our results we need to introduce preliminary materials on arcs.  
See \cite{I} for more on arcs on toric varieties.

\medskip

An \textit{arc} on $X_{\Sigma_0}$ is a morphism 
$
\eta\colon \mathrm{Spec} \, \k [[t]] \to X_{\Sigma_0}$, 
which corresponds to a homomorphism of $\k$-algebras 
$\eta^{*}\colon \k[\check{\sigma}_0 \cap M] \to \k[[t]]$.
The arc $\eta$ has its \textit{generic point in the torus} of $k^n$ if 
the series $\chi^v \circ \eta := \eta^* (\chi^v )$ is a nonzero element of $\k[[t]]$ 
for any $v \in \check{\sigma}_0 \cap M$. The arc is \textit{centered at the origin} of $\k^n$ if 
$\chi^v \circ \eta$ belongs to $(t)\k[[t]]$, for any nonzero $v\in\check{\sigma}_0 \cap M$. 
In this case, the arc $\eta$ induces a map of local $k$-algebras 
$\k[[\check{\sigma}_0 \cap M]] \to \k[[t]]$ 
which we denote also by $\eta^*$. 

\medskip 

If $\eta$  has its generic point in the torus $T_N$, then it has an associated 
\textit{order vector} $\ord (\eta) \in N$. It is defined as the lattice homomorphism 
$M \to \Z$ which sends $v \in M$ to the order of the Laurent series 
$\chi^v \circ \eta \in \k ((t))$. Notice that $\chi^v \circ \eta \in \k [[t]]$, 
for any $v \in \check{\sigma}_0 \cap M$. This implies that $\ord (\eta) \in \sigma_0$. 
The coordinates of the order vector $\ord (\eta)$ with respect to the basis 
$e_1, \dots, e_n$ are equal to 
\[
(\ord_t (X_1 \circ \eta), \cdots, \ord_t (X_n \circ \eta) ). 
\]

\begin{definition} \label{etaC}
Let $C$ be a branch embedded in $\k^n$ by the map  (\ref{one}). We associate to $C$ an arc 
$\eta_C\colon\bar{C} \to \k^n$ which is the composition of a normalization 
$\bar C\to C$ of $C$ with the inclusion $C \hookrightarrow \k^n$.
\end{definition}

In terms of the coordinates $(X_1,\dots, X_n)$ of $\k^n$ and a regular parameter $t$ of $\bar C$, 
the arc $\eta_C$ in Definition \ref{etaC} corresponds to a local homomorphism of $k$-algebras
\[
\eta_C^*\colon\k [[ X_1, \dots, X_n ]] = \k[[\check{\sigma}_0 \cap M]]  \to \k[[t]], 
\]
whose kernel is the \textit{defining ideal} $I_C$ of the branch $C$. 
By definition, the arc $\eta_C$ is centered at the origin of $\k^n$. If in addition, $\eta_C$ has its 
generic point in the torus of $\k^n$, then $\ord(\eta_C)$ belongs to the interior of $\sigma_0$.

\subsubsection{The case when $\eta_C$ has its generic point in the torus}

We start by discussing what happens to the order vector of an arbitrary arc $\eta$ on $X_{\Sigma_0}$ 
with generic point in the torus, when a toric modification is made. 

\begin{remark} \label{order-vector}
Let $\Sigma$ be a regular subdivision of $\Sigma_0$. By the valuative criterion of 
properness there is a unique arc 
$\tilde{\eta}\colon \mathrm{Spec }\, \k [[t]] \to X_{\Sigma}$ 
such that $\pi_{\Sigma_0}^\Sigma \circ \tilde{\eta} = \eta$. 
Since the modification $\pi_{\Sigma_0}^\Sigma$ is the identity on the torus, 
the lifted arc $\tilde{\eta}$ has also its generic point in the torus $T_N$ and by 
definition we get $
 \ord (\eta) = 
 \ord( \tilde{\eta})$. 
Let $\sigma \in \Sigma$ be a cone of dimension $n$ which contains the order vector 
$\ord ( \eta )$. Denote by $v_1, \dots, v_n$ the basis of the lattice $N$ which spans 
the regular cone $\sigma$. The affine toric variety $X_\sigma$ is isomorphic to $\k^n$ 
with coordinates $(U_1, \dots, U_n) = (\chi^{\check{v}_1}, \dots, \chi^{\check{v}_n})$, 
where $\check{v}_1, \dots, \check{v}_n$ is the dual basis of $v_1, \dots, v_n$.
Then, it follows that the coordinates of the vector 
$\ord (\tilde{\eta}) $ with respect to the basis $v_1, \dots, v_n$ of $N$ are equal to 
$(\ord_t (U_1 \circ \tilde{\eta} ) , \dots, \ord_t (U_n \circ \tilde{\eta} ) )$. 
\end{remark}

The following proposition provides a sufficient condition to guarantee the existence of 
a toric embedded resolution of the branch $C$, when $\eta_C$ has generic point 
in the torus of $\k^n$.  

\begin{proposition} \label{prop:prim}
Let $C$ be a branch embedded in $\k^n$ by the map  (\ref{one}). 
If the arc $\eta_C$ has its generic point in the torus of $\k^n$ and $\ord (\eta_C)$ 
is primitive with respect to the lattice $N$, then 
any regular subdivision $\Sigma$ of the fan $\Sigma_0$ such that 
the ray $\rho_C$ spanned by $\ord (\eta_C)$ belongs to $\Sigma$ 
defines a toric embedded resolution of $C$. 
\end{proposition}
\begin{proof}
Let us denote the toric modification $\pi_{\Sigma_0}^\Sigma$ simply by $\pi$ and 
the arc $\eta_C$ by $\eta$. 
Let $\tilde{\eta}$ be the unique lifting of $\eta$ to $X_{\Sigma}$. Then, the strict transform 
$\tilde{C}$ of $C$ is equal to the image of $\tilde{\eta}$, and $\tilde{\eta}$ 
is the composition of the normalization 
of $\tilde{C}$ with the inclusion $\tilde{C} \hookrightarrow X_\Sigma$.  

Since the order vector $\ord (\eta )$ is contained in the interior of the cone $\sigma_0$, Proposition 
\ref{prop:exc} implies that the orbit $O(\rho_C)$  is contained in the exceptional locus of $\pi$. 

Let $\sigma$ be a regular cone of $\Sigma$ spanned by a basis $v_1, \dots, v_n$ of $N$ such 
that $v_1 = \ord (\eta) $. Notice that we use here the assumption that $\ord (\eta)$ is 
primitive in the lattice $N$. 
The affine toric variety $X_\sigma$ is isomorphic with $\k^n$ with coordinates $(U_1, \dots, U_n)$. 
Then, it follows from  Remark \ref{order-vector}, that $\ord (\eta) = \ord (\tilde{\eta})$ and 
\begin{equation} \label{cond}
(\ord_t (U_1 \circ \tilde{\eta} ) , \dots, \ord_t (U_n \circ \tilde{\eta} ) ) = (1, 0, \dots, 0), 
\end{equation}
since $(1, 0, \dots, 0)$ is the vector of coordinates of $\ord (\eta)$ with respect to the basis 
$v_1, \dots, v_n$. The parametrization $\tilde{\eta}$ of $\tilde{C}$, in terms of the coordinates 
$(U_1, \dots, U_n)$, is given by $( U_1 \circ \tilde{\eta}, \dots, U_n \circ \tilde{\eta})$. 
Therefore, the condition (\ref{cond}) implies that the strict transform $\tilde{C}$ of $C$ is smooth 
and transversal to the orbit $O({\rho_C})$.
\end{proof}

\subsubsection{The case when $\eta_C$ does not have generic point in the torus}

We now deal with the generalization of Proposition \ref{prop:prim} to the case when $\eta_C$ does not 
have generic point in the torus of $\k^n$. Then, some of the  series $X_j \circ \eta $ are identically 
zero, say  for $j=s+1 , \dots, n$. 
The smaller coordinate subspace of $k^n$ which 
contains the branch $C$ is 
$ X_{s+1} = \cdots = X_n = 0$.
We are in the situation of Remark \ref{Rem: hyper}. If 
$
\sigma =  \R_{\geq 0} e_{s+1} + \cdots  + \R_{\geq 0} e_n
$, then the arc $\eta_C$ has its generic point in the torus 
$O(\sigma)_{\Sigma_0}$, $\eta_C$ factors through the orbit closure $X_{\mathrm{Star}(\sigma)_{\Sigma_0}} = \k^s$,
and it is centered at the origin of $\k^s$. 
That is, we have
$X_j \circ \eta_C \in (t) \, \k [[t]] \setminus \{ 0 \}$, for $j = 1 \dots, s$.
The order vector $\ord(\eta_C)$ is an element of the lattice $N(\sigma)$ which belongs to the interior 
of the image of the cone $\sigma_0$ under the canonical projection $N_\R\to N(\sigma)_\R$ (see Example 
\ref{ex:theta} and Remark \ref{order-vector}).

\begin{proposition} \label{subspace}
With the previous hypothesis and notation, let $\Sigma$ be a regular 
fan subdividing $\Sigma_0$ such that $\sigma \in \Sigma$. 
Assume that $\ord (\eta_C)$ is a primitive vector in the lattice $N(\sigma)$ and 
it spans a ray $\rho_C$ of the fan $\mathrm{Star}(\sigma)_{\Sigma}$. 
Then, $\pi_{\Sigma_0}^\Sigma$ is a toric embedded resolution of $C$. 
\end{proposition}

\begin{proof}
Since $\sigma \in \Sigma$ the image of $\eta_C$, which is equal to the branch $C$, 
is not contained in the discriminant locus of $\pi$ (see  Proposition \ref{prop:exc}).
This implies that the strict transform $\tilde C$ of $C$ exists. 
Recall that $\eta_C$ factors through $X_{\mathrm{Star}(\sigma)_{\Sigma_0}} = \k^s$, has generic 
point in the torus of this variety, and it is centered at the origin of $\k^s$. Proposition 
\ref{prop:prim} implies that the modification \eqref{f:modi} is a toric embedded resolution of 
$C \subset X_{\mathrm{Star}(\sigma)_{\Sigma_0}} $. By Remark  \ref{Rem: hyper} this implies the result.
\end{proof}

\begin{example} 
We keep the notations of Example \ref{uno}.
Let $C_1$, $C_2$ be the branches in $\k^3$ given by $\eta_{C_1}^* (X_1) =t^{4}$, 
$\eta_{C_1}^* (X_2) =t^{6}$, $\eta_{C_1}^* (X_3) =t^{13}$, and $\eta_{C_2}^* (X_1) =t^{2}$, 
$\eta_{C_2}^* (X_2) =t^{3}$, $\eta_{C_2}^* (X_3) =0$. 
If $\Sigma'$ is any regular subdivision of $\Sigma$, then
by Proposition \ref{prop:prim}
the toric modification defined by $\Sigma'$
is an embedded resolution of $C_1$. 
The branch $C_2$ is contained in the closure of 
the orbit of $\k^3$ associated with $\sigma= \R_{\geq 0} e_3$.
The cone spanned by $2 \bar{e}_1 + 3 \bar{e}_2 \in N(\sigma)$
belongs to the fan $\mathrm{Star}(\sigma)_{\Sigma'}$.
The cone $\sigma \in \Sigma_0$ belongs also to $\Sigma'$.
Since the order vector of $C_2$ is the primitive vector 
$2 \bar{e}_1 + 3 \bar{e}_2$,
the modification $\Sigma'$ induces also 
an embedded resolution of $C_2$ by Proposition 
\ref{subspace}. Therefore, in this case 
$\pi_{\Sigma_0}^{\Sigma'}$ is a toric embedded resolution of $C_1 \cup C_2$.
\end{example}

%-------------------------------------------------------------------------------------------------------

\subsection{Toric resolutions after reembedding}\label{sec:torenvironm}

In this section we prove that certain choices of functions define suitable 
\textit{torific embeddings} for a reduced curve singularity. We mean by this that the curve 
singularity can be resolved by one toric modification, after reembedding it in an 
affine space of possibly higher dimension.
This result is a generalization of a 
theorem of Goldin and Teissier, see \cite[Th. 6.1, Corollary]{GT}.

\begin{definition}\label{def:semigroup}
Let $C$ be a branch embedded in $\k^n$ by the map  (\ref{one}). The semigroup $\Gamma_C$ 
of the branch $C$ is 
$\Gamma_C = \{ \ord_t ( h \circ \eta_C ) \mid h \in \k [[ X_1, \dots, X_n ]], \, h \notin I_C \}.$
\end{definition}

The semigroup $\Gamma_C$ is a finitely generated semigroup of $(\Z_{\geq 0}, +)$ (see \cite{Z}).

\begin{lemma} \label{lem:gcd}
Let $C$ be a branch embedded in $\k^n$ by the map  (\ref{one}). Let $y_1, \dots, y_m$ be 
elements of the maximal ideal of $\k [[ X_1, \dots, X_n ]] $ such that: 
\begin{enumerate}
\item $y_1, \dots, y_m \notin I_C$,
\item the numbers $\ord_t (y_j \circ \eta_C) \in \Z_{>0}$, $j=1, \dots, m$, generate 
the semigroup $\Gamma_C$. 
\end{enumerate}
Then, the vector $(\ord_t (y_1 \circ \eta_C), \dots, \ord_t (y_m \circ \eta_C))$
is primitive in the lattice $\Z^m$. 
\end{lemma}
\begin{proof}
Since the local homomorphism $\mathcal O\to k[[t]]$ corresponding to the normalization $\bar C\to C$ 
induces an isomorphism between the fields of fractions, we must have that the group generated by 
$\Gamma_C$ is equal to $\Z$. This is equivalent to the condition 
$\gcd(\ord_t (y_1 \circ \eta_C), \dots, \ord_t (y_m \circ \eta_C))=1$. 
\end{proof}

With the hypothesis of Lemma \ref{lem:gcd}, if we
take $y_{m+1}, \dots, y_s$ in the maximal ideal of $\k [[ X_1, \dots, X_n ]]$
then the values 
$\ord_{t} (y_{\ell} \circ \eta_{C})$, for $\ell \in \{1, \dots, s \}$, belong to $\Gamma_{C} \cup \{ \infty \}$, 
and 
the finite values generate 
the semigroup $\Gamma_C$. We use this fact 
constantly afterwards.

\medskip 

Let $y_1, \dots, y_m$ be generators of the maximal ideal of $\k[[ X_1, \dots, X_n]]$. 
Then, the surjection
\[
\k [[ Y_1, \dots, Y_m ]] \to \k [[ X_1, \dots, X_n ]] \mbox{ such that } Y_j \mapsto y_j, 
\mbox{ for } j =1, \dots, m.
\]
defines an embedding 
\begin{equation} \label{re-emb}
\k^n \hookrightarrow \k^m.
\end{equation}
If $C$ is a reduced curve embedded in $\k^n$ by the map (\ref{one}) then, by composing the embedding $C \hookrightarrow \k^n$ 
with (\ref{re-emb}), we get an embedding 
\begin{equation} \label{re-emb-curv}
C \hookrightarrow \k^m. 
\end{equation}
If $C_j$ is a branch of $C= \cup_{l=1}^r  C_l$, then $C\hookrightarrow k^n$ induces an embedding 
$C_j\hookrightarrow k^n$. 
We denote  by $\eta_j \colon \mathrm{Spec} \, k[[t_j ]] \to \k^m $ the composition of $\eta_{C_j}$ with the map 
(\ref{re-emb}). It corresponds to a local homomorphism of $k$-algebras 
 \[
\eta_j^*\colon\k [[Y_1, \dots, Y_m ]] \to \k [[t_j]], \mbox{such that } Y_s \mapsto y_s \circ \eta_{C_j}, 
\mbox{ for } s=1, \dots, m. 
\]

\begin{corollary} \label{cor:resgen}
Let $C= \cup_{j=1}^r  C_j$ be a reduced curve with $r\geq1$ branches embedded in $\k^n$ by the map
(\ref{one}). Let $y_1, \dots, y_{m} \in \k [[ X_1, \dots, X_n ]] $ be such that: 
\begin{enumerate}
\item \label{enumresgen:me} $y_1, \dots, y_m$ generate the maximal ideal of  $\k [[X_1, \dots, X_n]]$,
\item\label{enumresgen:gptorus} $y_1, \dots, y_m \notin I_{C_j}$, for $1 \leq j \leq r$,
\item\label{enumresgen:primitive} the numbers $\ord_{t_j} (y_s \circ \eta_{C_j})\in\Z_{>0}$, 
$s= 1, \dots, m$, generate the group $\Z$, for $1 \leq j \leq r$.
\end{enumerate}
Let $C \hookrightarrow \k^{m}$ be the embedding defined as (\ref{re-emb-curv}) with 
respect to  $y_1, \dots, y_{m}$. 
Let $\Sigma$ be a regular subdivision of the fan ${\Sigma}_{0,m}$ 
such that the ray $\rho_j$ spanned by $\ord(\eta_j)$ belongs to $\Sigma$, for 
$1\leq j \leq r$. Denote by $\tilde{C}_j$ the strict transform of $C_j$ under 
$\pi^{\Sigma}_{{\Sigma}_{0,m}}$. 
If the special points of $\tilde{C}_j$, for $j=1, \dots,r$, are pairwise different,
then the toric modification $\pi^{\Sigma}_{{\Sigma}_{0,m}}$ is a toric 
embedded resolution of $C\subset \k^m$.
\end{corollary}

\begin{proof}
For $1\leq j \leq r$, by condition \eqref{enumresgen:gptorus}, the arc $\eta_j$ has its 
generic point in the torus of $\k^m$. Note that the coordinates of the vector 
$ \ord (\eta_j) $ with respect to the canonical basis 
$e_1, \dots, e_m$ of $N =\Z^m$ (see Notation \ref{quadrant}) are equal to 
\[
(\ord_{t_j} (Y_1 \circ \eta_j), \dots, \ord_{t_j} (Y_m \circ \eta_j)) = 
( \ord_{t_j} (y_1 \circ \eta_{C_j} ), \dots, \ord_{t_j} (y_m \circ \eta_{C_j} )).
\]
The vector $\ord (\eta_j)$ is primitive in the lattice $N$ by condition 
\eqref{enumresgen:primitive}. Then, we apply Proposition \ref{prop:prim} 
to the branch $C_j$ embedded in $\k^m$ by (\ref{re-emb-curv}) with respect to $y_1, \dots, y_{m}$, for $1\leq j \leq r$. 
It follows that the strict transform $\tilde C_j$ of $C_j$ is smooth and transversal 
to the orbit $O(\rho_j)_{\Sigma}$, 
which is contained in the exceptional locus of $\pi_{{\Sigma}_{0, m}}^\Sigma$. 
Since the special points of $\tilde{C}_j$ are pairwise different, the branches 
$\tilde C_j$ and $\tilde C_{j'}$ do not intersect if $j\neq j'$.
\end{proof}

\begin{remark} \label{enum:takegensem} 
Let $y_1, \dots, y_m$ be generators of the maximal ideal of 
$\k [[ X_1, \dots, X_n ]]$ such that 
for any $1 \leq j \leq r$ the semigroup $\Gamma_{C_j}$
is generated by the numbers
$\ord_{t_j} (y_{\ell} \circ \eta_{C_j})$, for $1 \leq \ell \leq m$
such that $\ord_{t_j} (y_{\ell} \circ \eta_{C_j}) \ne \infty$.
We can choose such functions in such a way that the three conditions of
Corollary \ref{cor:resgen} are satisfied. 
The condition \eqref{enumresgen:me} holds by assumption.
Let us choose a function $y$ in the maximal ideal  of 
$\k [[ X_1, \dots, X_n ]]$ such that $y \notin  I_{C_j}$
for  $1 \leq j \leq r$ and take an integer $p>0$.
 If $p$ is large enough 
then  for any $1 \leq j \leq r$ we have that:
\[
\ord_{t_j} ((y_i + y^p) \circ \eta_{C_j})
= 
\left\{ 
\begin{array}{lcl}
\ord_{t_j} (y_i\circ \eta_{C_j}) & \mbox{if} &  
\ord_{t_j} (y_i \circ \eta_{C_j}) \in \mathbb{N}
\\
\ord_{t_j} (y^p \circ \eta_{C_j}) & \mbox{if} &  
\ord_{t_j} (y_i \circ \eta_{C_j}) = \infty.
\end{array}
\right.
\]
Thus, if we replace $y_\ell$ by $y_\ell + y^p$, for $1 \leq \ell \leq m$, conditions \eqref{enumresgen:gptorus} and \eqref{enumresgen:primitive} of 
Corollary \ref{cor:resgen} are also satisfied (see Lemma \ref{lem:gcd}).
\end{remark}

\begin{corollary} \label{enum:diffvect}
Let $C= \cup_{j=1}^r  C_j$ be a reduced curve with $r\geq1$ branches embedded in $\k^n$ by  the map
(\ref{one}). There exist $m \geq 1$ and functions  $y_1, \dots, y_{m} $ such that
the conditions  \eqref{enumresgen:me}, 
\eqref{enumresgen:gptorus} and \eqref{enumresgen:primitive} of Corollary \ref{cor:resgen}
are satisfied and in addition, the vectors
\begin{equation} \label{eq:diffv}
( \ord_{t_j} (y_1 \circ \eta_{C_j}), \dots, \ord_{t_j} (y_{m}  \circ \eta_{C_j}) ) \mbox{, for } j=1, \dots, r, 
\end{equation} 
are pairwise different.  
Then, the modification $\pi^{\Sigma}_{{\Sigma}_{0,m}}$ of Corollary \ref{cor:resgen}
is a toric 
embedded resolution of $C \subset \k^m$.
\end{corollary}

\begin{proof} Assume first that we are given functions 
$y_1, \dots, y_{m} $ verifying the conditions of Corollary \ref{cor:resgen}
and such that the vectors \eqref{eq:diffv} are pairwise different.
This implies that the cones 
$\rho_1,\ldots,\rho_r$, and their orbits, are also pairwise different. Since the special 
point of $\tilde C_j$ belongs to the orbit $O(\rho_j)_{\Sigma}$, for $j=1, \dots, r$, 
these special points are also pairwise different. By Corollary \ref{cor:resgen} 
the modification $\pi^{\Sigma}_{{\Sigma}_{0,m}}$ is a toric embedded resolution of $C\subset \k^m$.

Let us take functions $y_1, \dots, y_{m'} $ such that the conditions  \eqref{enumresgen:me}, 
\eqref{enumresgen:gptorus} and \eqref{enumresgen:primitive} of Corollary \ref{cor:resgen} are satisfied (see Remark \ref{enum:takegensem}).
Up to relabelling the branches of $C$ we can assume that 
\begin{equation}\label{eq:wlog}
\ord_{t_1} (y_1 \circ \eta_{C_1}) \geq \dots \geq \ord_{t_{r}} (y_1 \circ \eta_{C_{r}}).
\end{equation} 
For $i=1, \dots, r-1$ we take $f_i \in I_{C_i}$ such that $f_i \notin I_{C_j}$ 
for any $1 \leq j \leq r$ and $j \ne i$.
Then we choose an integer $\ell_i$ large enough such that 
\begin{equation}\label{eq:wlog2}
\ord_{t_j} (f_i \circ \eta_{C_j}) <  \ord_{t_j} (y_1^{\ell_i} \circ \eta_{C_j}) , \mbox{ for } j=i+1, \dots, r.
\end{equation}
Set $m= m' + r-1$ and consider the functions: 
\[
y_1, \dots, y_{m'}, y_{m'+1} := f_1 + y_1^{\ell_1}, \dots,y_{m'+r-1} := f_{r-1} + y_1^{\ell_{r-1}}.
\]
Then, the vectors \eqref{eq:diffv} associated with the branches of $C$, with 
respect to 
the functions $y_1,\dots, y_{m}$,  are pairwise different. This follows from the inequalities 
\[
\ord_{t_i} (y_{m'+i} \circ \eta_{C_i}) =  \ord_{t_i} (y_1^{\ell_i} \circ \eta_{C_i})  \stackrel{\eqref{eq:wlog}}{\geq}  \ord_{t_j} (y_1^{\ell_i} \circ \eta_{C_j})  \stackrel{\eqref{eq:wlog2}}{>} \ord_{t_j} (f_i \circ \eta_{C_j}) =  \ord_{t_j} (y_{m'+i} \circ \eta_{C_j}),
\]
for $i=1, \dots, r-1$ and $j=i+1, \dots, r$. \end{proof}

\begin{example}\label{ex:stplanesing}
 The integers $\g_0:=8$, 
$\g_1:= 12$, $\g_2 := 26$ and $\g_3 := 53$ define the 
minimal system of generators of the semigroup of a plane branch
(see \cite[Sec 3.2]{T3}).
By the method explained in \textit{loc.cit.}
we get that a defining equation 
of such a plane branch $C$ is 
\[ ((Y^2-X^3)^2-X^5Y)^2-X^{10}(Y^2-X^3) =0.\] 
Then, $y_0=X$, $y_1=Y$, 
$y_2=Y^2-X^3$, and $y_3=(Y^2-X^3)^2-X^5Y$ 
verify the conditions of Corollary \ref{cor:resgen} (see Remark 
\ref{enum:takegensem}). But observe that $y_0$, $y_1$ and $y_3$
verify also these conditions because $\gcd(8,12,53)=1$. 
By Corollary \ref{cor:resgen}, any regular fan subdividing ${\Sigma}_{0,3}$ 
which contains the ray spanned by $(8,12,53)$ defines a toric embedded 
resolution of $C$. Next we study the behaviour of the branch $C\subset\k^3$, 
and also of the plane $\k^2\subset\k^3$, under the monomial map $\pi_{{\sigma}_{0,3}} ^\sigma$ 
corresponding to 
the regular cone $\sigma\subset{\sigma}_{0,3}$ generated by the vectors 
$v_1=(1,1,4)$, $v_2=(2,3,13)$, and $v_3=(8,12,53)$.
Consider $X_{{\Sigma}_{0,3}}=\k^3$ with coordinates $(X_1,X_2,X_3)$ 
as in Notation \ref{quadrant}, 
and $X_\sigma$ with coordinates $(U_1,U_2,U_3)$, where $U_i=\chi^{\check{v}_i}$. The 
equation that defines $\k^2\subset\k^3$ is 
$ X_3-(X_2^2-X_1^3)^2+X_1^5X_2=0$.
In addition, the points of  $C \subset \k^3$ must satisfy the following equation $
X_3^2-X_1^{10}(X_2^2-X_1^3)=0$. 
The points in the strict transform $\tilde C$ of the curve $C$ 
by the monomial map $\pi_{{\sigma}_{0,3}} ^\sigma$
satisfy the following 
two equations: 
\begin{align}
    U_2U_3^5-(1-U_1)^2+U_1^2U_2U_3^4=0,\label{eq:stplane}\\
    U_3^2-U_1^4(1-U_1)=0,\label{eq:stBTcurve}
\end{align}
where \eqref{eq:stplane} is the equation of the strict transform of $\k^2$. 
Starting with \eqref{eq:stBTcurve} and then substituting in \eqref{eq:stplane}, 
we obtain a parametrization of $\tilde C$ of the form 
$$U_1=1-t^2+\ldots,\quad U_2=1+\ldots,\quad U_3=t+\ldots,$$
where we have omitted higher order terms (compare with \eqref{cond}). As expected, 
$\tilde C$ is smooth and transverse to $U_3=0$. Notice that the intersection of the 
surfaces defined by \eqref{eq:stBTcurve} and \eqref{eq:stplane}, contains also
an exceptional component given by $U_1=1,U_2=t,U_3=0$, and  the 
strict transform 
of the plane \eqref{eq:stplane} is singular at the point $(1,0,0)$.

One may consider this toric  resolution 
as 
a sort of \textit{toric embedded normalization} of $C$, since we do not care about 
the effect of this modification on the original plane $\k^2 \supset C$. 
\end{example}

In Corollary \ref{enum:diffvect} we  showed the existence of functions defining 
a suitable embedding. We consider a different choice of functions in Theorem \ref{thm:toremrescurves} below. We are grateful to Javier Fern\'andez de Bobadilla for sharing some ideas which were helpful in the formulation of the result. We will need the following elementary lemma: 

\begin{lemma} \label{lem:basis}
Denote by $N'$ the sublattice of $N$ spanned by $v_1, \dots, v_{s-1}$. If
$v_1, \dots, v_{s-1}$ are part of a basis of the lattice $N$ 
and if the canonical image of  a vector  $w \in N$ is a primitive vector 
in the lattice $N/ N'$, then $v_1, \dots, v_{s-1}, w$ are part of a basis 
of the lattice $N$. 
\end{lemma}
\begin{proof}
Let $v_1, \dots, v_{s-1}, v_s, \dots, v_n$ be a basis of the lattice $N$. 
Denote by $\bar{u}$ the class of a vector $u \in N$ in the quotient lattice $N/N'$. 
We have that $\bar{v}_s, \dots, \bar{v}_n$ is a basis of $N/N'$. 
By hypothesis $\bar{w}$ is a primitive vector in the lattice $N/N'$, 
hence there are elements $ w_s:= w, w_{s+1}, \dots, w_n\in N$ such that 
$\bar{w}_s = \bar{w},  \bar{w}_{s+1}, \dots, \bar{w}_n$ is a basis of $N/N'$. 
If we expand $\bar{w}_j = \sum_{l=s}^n a_{l,j} \bar{v}_l$, for $j=s,\dots,n$, 
it follows that the matrix $ A= (a_{l,j})$ has determinant $\pm1$. 
It is enough to check that $v_1, \dots, v_{s-1}, w_{s}, w_{s+1}, \dots, w_n$ 
is a basis of $N$. The matrix whose columns are the coordinates of these vectors 
in terms of $v_1, \dots,v_n$ has a block structure of the form 
\[
\begin{pmatrix}
I_{s-1}  & * 
\\
0 & A
\end{pmatrix}
\]
where $I_{s-1}$ denotes the identity matrix of size $s-1$. This matrix is unimodular, 
since $A$ is. Therefore $v_1, \dots, v_{s-1}, w_{s}, w_{s+1}, \dots, w_n$ is a basis of 
$N$ as claimed. 
\end{proof}

\begin{theorem} \label{thm:toremrescurves}
Let $C= \cup_{j=1}^r  C_j$ be a reduced curve with $r$ branches embedded in $\k^n$ by the map (\ref{one}). 
Let $y_1,\dots,y_m,y_{m+1},\ldots,y_{m+r} \in \k [[X_1, \dots, X_n]]$ be such that: 
\begin{enumerate}
\item $y_1, \dots, y_m$ generate the maximal ideal of  $\k [[X_1, \dots, X_n]]$,
\item\label{enum:gptorus} $y_1, \dots, y_m \notin I_{C_j}$, for $1 \leq j \leq r$,
\item\label{enum:primitive} the numbers $\ord_{t_j} (y_s \circ \eta_{C_j})\in\Z_{>0}$, 
$s= 1, \dots, m$, generate the group $\Z$, for $1 \leq j \leq r$,
\item\label{enum:coordsubspace} $y_{m+j} \in \bigcap_{1 \leq i \leq r}^{ i \ne j} I_{C_i}$ 
and $y_{m+j} \notin I_{C_j}$, for $1 \leq j \leq r$.  
\end{enumerate}
Let $C \hookrightarrow \k^{m+r}$ be the embedding defined as (\ref{re-emb-curv}) with 
respect to  $y_1, \dots, y_{m+r}$. 
Then, there exists a toric modification of $X_{{\Sigma}_{0,m+r}}=\k^{m+r}$ which is 
a toric embedded resolution of $C$.
\end{theorem}
\begin{proof}  
We denote by $\Sigma_0$ the fan ${\Sigma}_{0,m+r}$ and use Notation \ref{quadrant}. 
Condition \eqref{enum:coordsubspace} implies that for $1 \leq j \leq r$, the branch 
$C_j \subset \k^{m+r} $ is contained in the coordinate subspace:
\[
Y_{m+1} =  \cdots = Y_{m+j-1 } =   Y_{m+j+1 } = \cdots  =Y_{m+r} = 0.
\]
This subspace is the closure of the orbit $O ({{\sigma}_j})_{\Sigma_0}$, where 
the cone ${\sigma}_j \in {\Sigma_0}$ is spanned by  the vectors
$e_{m+1},  \dots , e_{m+j-1 },   e_{m+j+1 }, \dots , e_{m+r},$
of the canonical basis $N=\Z^{m+r}$. 
Denote by $\bar{e}_{j, i}$ the image of $e_i$ under the canonical projection $N \to N({\sigma}_j)$, 
for $i=1, \dots, m+r$. 
By conditions \eqref{enum:gptorus}-\eqref{enum:primitive}, 
the arc $\eta_j$ (i.e.,\ the composition of $\eta_{C_j}$ with $k^n\hookrightarrow k^{m+r}$) 
has its generic point in the torus $O ({{\sigma}_j})_{\Sigma_0}$, and the order vector $\ord (\eta_j )$ 
is primitive in the lattice $N({\sigma}_j)$. We have that
\begin{equation} \label{f:etaj}
\ord (\eta_j ) =w_{j, 1} \bar{e}_{j, 1}+ \cdots + w_{j, m}  \bar{e}_{j,m} + w_{j, m+j}  \bar{e}_{j,m+j}, 
\end{equation}
in terms of the basis of $\bar{e}_{j,1}, \dots, \bar{e}_{j,m}, \bar{e}_{j,m+j}$ of $N({\sigma}_j)$, 
with all the $w_{j,i}$ in $\Z_{>0}$.  
Let us consider the vectors: 
\[
w_0: =  w_{1, m+1} e_{m+1} + \cdots + w_{r, m+r} e_{m+r}, 
\quad
w_j :=  w_{j, 1} {e}_{1}+ \cdots + w_{j, m}  {e}_{m} , \mbox{ for } 1 \leq j \leq r. 
\]
Notice that the canonical image of $w_0 + w_j$ in $N({\sigma}_j)$ 
is precisely the primitive vector $\ord (\eta_j) \in N({\sigma}_j)$. 
The cone 
\[\theta_j:=\R_{\geq 0} (w_0 + w_j) + {\sigma}_j \subset N_\R\]
is regular for the lattice $N$ by Lemma \ref{lem:basis}, 
and its canonical projection in $N({\sigma}_j)_\R$ is equal to the ray 
$\rho_j:=\R_{\geq0}\ord(\eta_j)\subset N({\sigma}_j)_\R$.

Let us prove that the cones $\theta_j$ for $j =1, \dots, r$ are pairwise different.
We assume that there are integers $ 1 \leq s < l \leq r$ such that 
$\theta := \theta_s = \theta_l$, and proceed by contradiction. 
By construction, the cone $\theta$ is of dimension $r$. 
Since ${\sigma}_s$ and ${\sigma}_l$ are faces of $\theta$, it follows that 
\[ {\sigma}_s + {\sigma}_l = \R_{\geq 0} e_{m+1} + \cdots + \R_{\geq 0} e_{m+ r} 
\]
is an $r$ dimensional cone contained in $\theta$. Assume that there is a vector 
$w \in \theta \setminus ({\sigma}_s + {\sigma}_l)$. 
Since $\theta$ is contained in the cone $\R_{\geq 0} e_{1} + \cdots + \R_{\geq 0} e_{m+ r}$,
this  implies that $w,  e_{m+1} , \dots, e_{m+ r}$ are linearly independent,  
contradicting the assumption that $\theta$ has dimension $r$. 
Hence, we get that $\theta = \R_{\geq 0} e_{m+1} + \cdots + \R_{\geq 0} e_{m+ r}$.
But, then the canonical projection of 
$\theta$ onto $N({\sigma}_l)_\R$ is equal to the ray 
spanned by $\bar{e}_{l, m+l}$. 
This ray is different from the ray $\rho_l$ because 
the coefficients $w_{l, i}$  appearing in the 
expansion (\ref{f:etaj}) of $\ord (\eta_l)$ are non zero. 
This contradiction shows 
that the cones $\theta_j$, for $1\leq j \leq r$ 
are pairwise different, as claimed.

Assume that $\Sigma$ is a regular subdivision of ${\Sigma_0}$ which contains the cones
$\theta_j$, for $j = 1, \dots, r$. 
Then, we have ${\sigma}_j \in \Sigma$ since ${\sigma}_j$ is a face of $\theta_j$, and 
$\rho_j$ belongs to the fan $\mathrm{Star} ({\sigma}_j)_{\Sigma}$. 
This implies that the toric modification $\pi_{{\Sigma_0}}^{\Sigma}$ 
is a toric embedded resolution of the branches $C_j$ for $j=1, \dots, r$, by Proposition \ref{subspace}.
The strict transform of $C_j$ by this modification is smooth and 
transversal to the orbit associated with $\theta_j$, for $1 \leq j \leq r$.
Since these cones are pairwise different, it follows 
that the strict transform of $C$ by $\pi_{{\Sigma_0}}^\Sigma$ is smooth.

In order to finish the proof it remains to show that such a regular fan $\Sigma$ exists. 
It is enough to show the existence of a fan $\Sigma'$ subdividing ${\Sigma_0}$ such that 
the regular cone $\theta_j$ belongs to $\Sigma'$, for $1 \leq j \leq r$, and then 
obtain $\Sigma$ as a regular subdivision of $\Sigma'$.

Next, we show that the set consisting of the cones $\theta_j$, for $j=1, \dots, r$, and their faces 
is a fan, which contains cones of dimension at most $r$. Then, there exists a fan $\Sigma'$ subdividing 
${\Sigma_0}$ and containing it, and this ends the proof. Set $1 \leq j, j' \leq r$ and let us check 
that $\theta_j \cap \theta_{j'}$ is a common face of $\theta_j$ and $\theta_{j'}$.
If $w_j= w_{j'}$ then $\theta_j \cap \theta_{j'} = {\sigma}_j \cap {\sigma}_{j'}  + \R_{\geq 0} (w_j +w_0)$.
If $w_j \ne w_j'$ then $\theta_j \cap \theta_{j'} = {\sigma}_j \cap {\sigma}_{j'}$.
In both cases the assertion holds by a direct computation when we express a vector 
$v \in \theta_j \cap \theta_{j'}$ in terms of the basis $e_1, \dots, e_{m+r}$. 
\end{proof}

In the proof of Theorem \ref{thm:toremrescurves}, we show the existence of a fan $\Sigma'$ subdividing 
${\Sigma_0}$, which can be obtained explicitly as follows. 

\begin{remark}
Let $N_1$ (resp. $N_2$) be the lattice spanned by $e_1, \dots, e_m$ 
(resp. $e_{m+1}, \dots, e_{m+r}$), and let ${\Sigma}_1$ (resp. ${\Sigma}_2$) be 
the fan of faces of the cone spanned by $e_1, \dots, e_m$ (resp. $e_{m+1}, \dots, e_{m+r}$). 
The lattice $N$ is equal to the direct sum of its sublattices $N_1$ and $N_2$. 
The fan ${\Sigma_0} $ is the \textit{direct} Minkowski sum ${\Sigma}_1 \oplus {\Sigma}_2$
of the fans ${\Sigma}_1$ and ${\Sigma}_2$, in the sense that any cone $\tau \in {\Sigma_0}$ is 
the Minkowski sum of unique cones $\tau_1 \in {\Sigma}_1$ and $\tau_2 \in {\Sigma}_2$. 
Let $\Sigma_1'$ be a regular subdivision of ${\Sigma}_1$ containing 
the rays spanned by 
$w_j$, for $j= 1, \dots, r$.
Let us denote by $\Sigma_2'$ the star subdivision of ${\Sigma}_2$ at the vector $w_0$. 
The  $r$-dimensional cones  of $\Sigma_2'$ are 
${\sigma}_j + \R_{\geq 0} w_0 $, for $j=1, \dots, r$.
The $(r+1)$-dimensional cone 
$\tau_j := \R_{\geq 0} w_j +   \R_{\geq 0} w_0 +{\sigma}_j $
belongs to the direct Minkoski sum of fans $ \Sigma_1 '\oplus \Sigma_2'$
and contains the cones ${\sigma}_j$ and $\R_{\geq 0} w_j + \R_{\geq 0} w_0$ as faces. 
Apply,  for $1 \leq j \leq r$,  a star subdivision at the vector
$\R_{\geq 0} (w_j + w_0)$. 
This combinatorial operation only affects the cones
which contain $\R_{\geq 0} w_j +  \R_{\geq 0} w_0$ as a face
(see \cite[Lemma 11.1.3]{CLS}).
In particular, it replaces
the cone $\tau_j $  by the cones 
$\R_{\geq 0} w_j + \R_{\geq 0} (w_0 + w_j) +  {\sigma}_j $,
$\R_{\geq 0} w_0 + \R_{\geq 0} (w_0 + w_j) +  {\sigma}_j $,
and $ \theta_j =\R_{\geq 0} (w_0 + w_j) +  {\sigma}_j $, and their faces.
Therefore, the fan $\Sigma'$ constructed by this process, starting from 
$\Sigma_1' \oplus \Sigma_2'$, contains the cones $\theta_j$,  for $1 \leq j \leq r$.
\end{remark}

%--------------------------------------------------------------------------------------------------------------

\section{Resolving a plane curve singularity with one toric morphism} \label{resolving_reduced}

In this section, we prove that a singular curve 
$C = \cup_{j=1}^r C_j$ with $r\geq1$ branches 
admits a non degenerate reembedding,
in such a way that an embedded 
resolution is obtained by one toric morphism. 
This result is a generalization of 
theorems of Goldin and Teissier in \cite{GT}. 

\subsection{Background on plane curve singularities and divisorial valuations}\label{sec:dival}

In this section we explain the basic notation and conventions used in the rest of the paper 
about plane curve singularities and divisorial valuations. 
We denote by $S$ the affine space $\mathrm{Spec}\, \k[X, Y]$ and 
by $O$ its closed point corresponding to the maximal ideal $(X,Y)$. 
We denote by $R$ the ring $k [[ X, Y]]$ and by $m_R$ its maximal ideal as before.

\medskip

Let $C = \cup_{j=1}^r C_j$ be a  plane curve with $r\geq1$ branches. 
Denote by $\mathcal{O}$ its local ring. If $x, y$ generate the maximal ideal of 
$\mathcal{O}$, then the embedding of $C \subset  \mathrm{Spec}\,R$ is defined by a local homomorphism of 
$k$-algebras
\begin{equation} \label{one-2}
R \longrightarrow \mathcal{O}, \mbox{ such that } X \mapsto x, Y \mapsto y.
\end{equation}
If $D$ is a closed subscheme of $\mathrm{Spec}\,R$
defined by a principal ideal of $R$, then we call any generator $f\in R$ of this ideal a 
\textit{defining function} of $D$ and we write $D=Z(f)$. If $f\in m_R$ generates the 
kernel of the map \eqref{one-2}, then $f$ is a defining function of $C$ and $\set{C_j}_{j=1}^r$ 
is in bijection with the set of irreducible factors of $f$ in $R$.

\medskip

If $D_1=Z(f_1)$ and $D_2=Z(f_2)$, we denote by $(D_1 \cdot D_2) := \dim_{\k} {R}/(f_1, f_2)$ 
the {\it intersection multiplicity} of $D_1$ and $D_2$. 
If $D_1$ is a branch, then $(D_1 \cdot D_2)$ is equal to the order of the series 
$f_2\circ \eta_{D_1}\in\k[[t]]$, where $\eta_{D_1}:\mathrm{Spec}\,\k[[t]]\to  \mathrm{Spec}\,R$
is the arc associated to $D_1$ (see Definition \ref{etaC}). 
 
A pair $(x,y) \in {R}^2$ is a {\it local coordinate system} on $\mathrm{Spec}\,R$ if $x,y$ generate 
the maximal ideal $m_R$. Then, if $L = Z(x)$ and $L' =Z(y)$ we say
that $(L,  L')$ is a {\it cross} at $O$.

 \begin{remark} \label{WPT}
  Once the local coordinate system $(x, y)$ is fixed, the Weierstrass Preparation Theorem 
allows us to express any branch $A \ne Z(x)$ in the form $A = Z (x_A)$ for a 
monic element $x_A \in \k[[x]][y]$ such that $\deg_y x_A = (Z(x) \cdot A)$. 
\end{remark}

If $A$ is a plane branch on $S$ then the semigroup $\Gamma_A$ is finitely generated. 
We mention below some classical properties (see \cite[Lemma 1.2 and Theorem 3.2]{GBP}). 

\begin{notation} \label{not:GA}
 Let $A$ be a plane branch on $S$. Denote by $\g_0,\g_1,\dots,\g_{g(A)}$ 
 the minimal generating system of the semigroup $\Gamma_A$, where $\g_0 < \g_1< \cdots < \g_{g(A)}$.
We set 
$e_0 := \g_0$, $e_j = \mathrm{gcd} (e_{j-1}, \g_j)$, and 
$n_j = e_{j-1} / e_j $  for $j = 1, \dots, g(A)$.  
We set $n_0:=0$. 
\end{notation}

\begin{proposition} \label{GA} $\,$ 
\begin{enumerate}
\item For every $j\in\set{1,\ldots,g(A)}$  one has $n_j > 1$ and 
\begin{equation}
\label{rel-sem}
n_j \g_j = b_0^j \g_0 + b_1^j  \g_1 + \cdots + b_{j-1}^j \g_{j-1} 
\mbox{ with } 0 \leq b_i^j < n_i , \mbox{ for } i= 1, \dots, j-1.
\end{equation}

\item We have the inequalities: 
\begin{equation}
\label{gg-ineq}
n_j \g_j < \g_{j+1}, \mbox{ for }  j = 1, \dots, g(A)-1.
\end{equation}

\item There exist a local coordinate system $(x_0, x_1)$ on $S$ and irreducible elements 
$x_j \in \k [[x_0]][x_1]$
with $\deg_{x_1} x_j = n_1 \dots n_{j -1}$, for $j \in \set{2,\ldots,g(A)}$, such that $(Z(x_j) \cdot A) = \g_j$ for $j \in \set{0,\ldots,g(A)}$.

\end{enumerate}
\end{proposition}

%--------------------------------------------------------------------------------------------------------------

A \textit{model} of $(S,O)$ is a proper birational morphism 
$\psi\colon (\S(\psi), E(\psi)) \to (\S, O)$, obtained as a composition of blow ups of a  finite set of
infinitely near points of $O$. If this set is empty then $\psi$ is the identity map of $S$ and $E(\psi) =\emptyset$. Otherwise, the preimage $E(\psi)= \psi^{-1} (O)$, seen as a reduced 
divisor on $\S(\psi)$, is called the \textit{exceptional divisor} of $\psi$. 
It has simple normal crossings and its irreducible components are projective lines 
which are called the {\it exceptional prime divisors} of $\psi$. 
We denote by $E_O$ the prime exceptional divisor created by blowing up $O$ in $S$.

\medskip

The \textit{dual graph} $G(\psi)$ of the model $\psi$ is a tree
whose set of vertices in bijection with the set of exceptional prime divisors of $\psi$, 
and where two vertices are joined by an edge if and only if the corresponding 
exceptional prime divisors intersect. 
Recall that the \textit{valency} of a vertex  in a tree is the number of edges 
incident to it. A vertex of a tree is an \textit{end} (resp. a  \textit{ramification vertex}) 
if it has valency equal to one (resp. $\geq 3$). 
We denote by $\mathcal E_\psi$ the set of exceptional prime divisors of $\psi$ which  
correspond to ends of $G(\psi)$.

\begin{notation}
If $C$ is a curve on $\S$ we denote by $C^\psi$ its strict transform on $\S(\psi)$. 
If the model $\psi$ is clear from the context 
we often denote by the same letter a curve $C$ and its strict transform $C^\psi$ on this model. 
This applies also in the case of an exceptional prime divisor $E$ which was created in 
another model dominated by $\psi$ and its strict transform $E^\psi$. 
For instance, we will often denote $E_O^\psi$ simply by $E_O$. 
\end{notation}

\begin{definition}
Let $E$ be an irreducible component of $E(\psi)$. A plane branch $C$ on $S$ is 
a \textit{curvetta} of $E$ at the model $\psi$ if $C^\psi$ is smooth and 
transversal to $E$ at a smooth point of $E(\psi)$. We denote by $\mathcal C_E (\psi)$ 
the set of curvettas of $E$ at the model $\psi$. We call $\nu_{E}$ the 
\textit{divisorial valuation} on $(S,O)$ defined by $E$, that is, $\nu_E$ maps 
$h \in {R} \setminus \{ 0 \}$ to the order of vanishing along $E$ of $h \circ \psi$. 
\end{definition}

\begin{comment}
Since $\nu_{E}$ is centered at $O$ and ${R}$ is noetherian, the quotient 
$\mathcal{P}_{\phi}^{\nu_E} /\mathcal{P}_{\phi}^{\nu_E}$ is a $\k$-module of finite type. 
Following \cite{T3}, we define the graded algebra 
$$gr_{\nu_E}{R}=\bigoplus_{\phi\in\Z_{\geq0}} 
\frac{\mathcal{P}^{\nu_E}_{\phi}}{\mathcal{P}^{\nu_E}_{\phi+1}}.$$ 
Observe that the sum on the right is in fact indexed by the semigroup 
$\Gamma_{\nu_E} := \nu_E({R}\setminus\set{0})$, which is a 
subsemigroup of $\Z_{\geq0}$ and hence a finitely generated semigroup. 
As a consequence, $gr_{\nu_{E} }{R}$ is finitely generated as $\k$-algebra. 
We call $in_{\nu}$ the natural map $in_{\nu_E}:{R}\to gr_{\nu_E}{R}$, which sends 
$h\in {R}$ to its initial form with respect to the valuation $\nu$, that is 
$in_{\nu}(h)=h~\mbox{mod}~\mathcal{P}^E_{\nu(h)+1}.$
\end{comment}

\begin{remark} \label{rem:val-int}
The value of $\nu_{E}$ at $h \in R \setminus \{ 0 \}$ is characterized in terms 
of intersection multiplicity of curvettas as follows: 
\begin{equation} \label{val-int}
\nu_E (h)=\min\{ (Z(h)  \cdot D) \mid D\in\mathcal C_{E}(\psi) \}. 
\end{equation}
Furthermore, the minimum in \eqref{val-int} is attained whenever the strict 
transforms $Z(h)^\psi $ and $D^\psi$ intersect $E$ at different points (see \cite[Section 7]{S}). 
\end{remark}

\begin{definition}\label{def:sigma}
A model $\psi$ is an \textit{embedded resolution} of a curve $C$ on $S$ 
if the divisor ${C}^\psi + E(\psi)$ has simple normal crossings. 
If $C=\cup_{j=1}^r{C_j}$ is a curve with $r$ branches, we denote by 
$E_{C_j}$ the unique exceptional prime divisor of $\psi$ that intersects 
$C_j^\psi$, for $j=1,\ldots,r$. 
The \emph{dual graph of the total transform} of $C$ under $\psi$ 
is the tree $G(\psi, C)$  obtained by adding to $G(\psi)$, 
a new vertex corresponding to $C_j^\psi$, and an edge joining it to 
the vertex of $G(\psi)$ corresponding to $E_{C_j}$, for $j=1,\ldots,r$.
\end{definition}

\begin{definition}\label{def:minresdivval}
A model $\psi$ is the \textit{minimal resolution} of the divisorial valuation 
$\nu_E$ if $\psi$ is the composition of 
a finite sequence of point blow-ups above $O$,  
\[
\psi\colon \S(\psi) = \S_n\xrightarrow{\psi_{n}} \S_{n-1}\xrightarrow{\psi_{n-1}}\cdots\xrightarrow{\psi_2} 
\S_1\xrightarrow{\psi_1} \S_0=\S
\]
such that
$E$ is the exceptional divisor of $\psi_n$ and the center of $\psi_{i+1}$ 
in $\S_{i}$ belongs to the exceptional divisor of $\psi_i$ for $i=1,\ldots,n-1$. 
\end{definition}

\begin{figure}[h!]
    \begin{center}
\begin{tikzpicture}[scale=0.7]
\tikzstyle{every node}=[font=\small]

 \begin{scope}[shift={(0,0)},scale=0.7]
 
\draw [-, color=black, thick](0,0) -- (1,0);
\draw [-, color=black, thick, dashed](1,0) -- (2,0);
\draw [-, color=black, thick](2,0) -- (4,0);
\draw [-, color=black, thick, dashed](3,0) -- (6,0);
%\draw [-, color=black, thick](6,0) -- (10,0);

\draw [-, color=black, thick](5.8,0) -- (8.2,0);
\draw [-, color=black, thick, dashed](7,0) -- (10,0);

\draw [-, color=black, thick](10,0) -- (11,0);

\draw [-, color=black, thick, dashed] (3,0)--(3,3);
\draw [-, color=black, thick] (3,0)--(3,1);
\draw [-, color=black, thick] (3,2)--(3,3);
\node[draw,circle, inner sep=1.pt,color=black, fill=black] at (3,0.8){};
\node[draw,circle, inner sep=1.pt,color=black, fill=black] at (3,2.2){};
\node[draw,circle, inner sep=1.pt,color=black, fill=black] at (3,3){};

\draw [-, color=black, thick, dashed] (7,0)--(7,3);
\draw [-, color=black, thick] (7,0)--(7,1);
\draw [-, color=black, thick] (7,2)--(7,3);
\node[draw,circle, inner sep=1.pt,color=black, fill=black] at (7,0.8){};
\node[draw,circle, inner sep=1.pt,color=black, fill=black] at (7,2.2){};
\node[draw,circle, inner sep=1.pt,color=black, fill=black] at (7,3){};

\draw [-, color=black, thick, dashed] (11,0)--(11,3);
\draw [-, color=black, thick] (11,0)--(11,1);
\draw [-, color=black, thick] (11,2)--(11,3);
\node[draw,circle, inner sep=1.pt,color=black, fill=black] at (11,0.8){};
\node[draw,circle, inner sep=1.pt,color=black, fill=black] at (11,2.2){};
\node[draw,circle, inner sep=1.pt,color=black, fill=black] at (11,3){};

\node[draw,circle, inner sep=1.pt,color=black, fill=black] at (0,0){};
\node[draw,circle, inner sep=1.pt,color=black, fill=black] at (0.8,0){};
\node[draw,circle, inner sep=1.pt,color=black, fill=black] at (2.2,0){};
\node[draw,circle, inner sep=1.pt,color=black, fill=black] at (3,0){};
 \node[draw,circle, inner sep=1.pt,color=black, fill=black] at (3.8,0){};
 \node[draw,circle, inner sep=1.pt,color=black, fill=black] at (6,0){};
\node[draw,circle, inner sep=1.pt,color=black, fill=black] at (7,0){};
\node[draw,circle, inner sep=1.pt,color=black, fill=black] at (7.8,0){};
\node[draw,circle, inner sep=1.pt,color=black, fill=black] at (10.2,0){};
\node[draw,circle, inner sep=1.7pt,color=black, fill=black] at (11,0){};
\node [right] at (11,0) {$E$};

\end{scope}
\begin{scope}[shift={(10,0)},scale=0.7]
 
\draw [-, color=black, thick](0,0) -- (1,0);
\draw [-, color=black, thick, dashed](1,0) -- (2,0);
\draw [-, color=black, thick](2,0) -- (4,0);
\draw [-, color=black, thick, dashed](3,0) -- (6,0);

\draw [-, color=black, thick](5.8,0) -- (8.2,0);
\draw [-, color=black, thick, dashed](7,0) -- (10,0);

\draw [-, color=black, thick](10,0) -- (11,0);

\draw [-, color=black, thick](11,0) -- (12,0);
\draw [-, color=black, thick, dashed](12,0) -- (13,0);
\draw [-, color=black, thick](13,0) -- (14,0);

\node[draw,circle, inner sep=1.pt,color=black, fill=black] at (11.8,0){};
\node[draw,circle, inner sep=1.pt,color=black, fill=black] at (13.2,0){};
\node[draw,circle, inner sep=1.7pt,color=black, fill=black] at (14,0){};
\node [right] at (14,0) {$E$};

\draw [-, color=black, thick, dashed] (3,0)--(3,3);
\draw [-, color=black, thick] (3,0)--(3,1);
\draw [-, color=black, thick] (3,2)--(3,3);
\node[draw,circle, inner sep=1.pt,color=black, fill=black] at (3,0.8){};
\node[draw,circle, inner sep=1.pt,color=black, fill=black] at (3,2.2){};
\node[draw,circle, inner sep=1.pt,color=black, fill=black] at (3,3){};

\draw [-, color=black, thick, dashed] (7,0)--(7,3);
\draw [-, color=black, thick] (7,0)--(7,1);
\draw [-, color=black, thick] (7,2)--(7,3);
\node[draw,circle, inner sep=1.pt,color=black, fill=black] at (7,0.8){};
\node[draw,circle, inner sep=1.pt,color=black, fill=black] at (7,2.2){};
\node[draw,circle, inner sep=1.pt,color=black, fill=black] at (7,3){};

\draw [-, color=black, thick, dashed] (11,0)--(11,3);
\draw [-, color=black, thick] (11,0)--(11,1);
\draw [-, color=black, thick] (11,2)--(11,3);
\node[draw,circle, inner sep=1.pt,color=black, fill=black] at (11,0.8){};
\node[draw,circle, inner sep=1.pt,color=black, fill=black] at (11,2.2){};
\node[draw,circle, inner sep=1.pt,color=black, fill=black] at (11,3){};

\node[draw,circle, inner sep=1.pt,color=black, fill=black] at (0,0){};
\node[draw,circle, inner sep=1.pt,color=black, fill=black] at (0.8,0){};
\node[draw,circle, inner sep=1.pt,color=black, fill=black] at (2.2,0){};
\node[draw,circle, inner sep=1.pt,color=black, fill=black] at (3,0){};
 \node[draw,circle, inner sep=1.pt,color=black, fill=black] at (3.8,0){};
 \node[draw,circle, inner sep=1.pt,color=black, fill=black] at (6,0){};
\node[draw,circle, inner sep=1.pt,color=black, fill=black] at (7,0){};
\node[draw,circle, inner sep=1.pt,color=black, fill=black] at (7.8,0){};
\node[draw,circle, inner sep=1.pt,color=black, fill=black] at (10.2,0){};
\node[draw,circle, inner sep=1.pt,color=black, fill=black] at (11,0){};

\end{scope}
\end{tikzpicture}
\end{center}
 \caption{The shape of the dual graph in the case $\ell_{E} = 0$ on the left (resp.  $\ell_{E} > 0$ on the right).
 The vertex corresponding to $E$ is marked by a bigger bullet.}
\label{fig:dualgraphs}
   \end{figure}
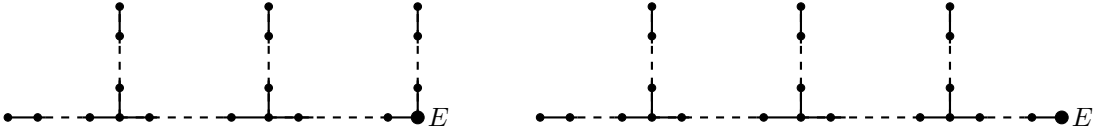

If $\psi$ is the minimal resolution of $\nu_E$ and $A$ is any curvetta in 
$\mathcal{C}_E (\psi)$, then by definition $\psi$ is an embedded resolution 
of $A$. However, it may not be the minimal one. 
Denote by $\pi\colon (S(\pi), E(\pi)) \to (S,O)$ the minimal embedded 
resolution of $A$ and by $E_A$ the irreducible component of $E(\pi)$ 
which intersects $A^\pi$. Since by definition $\psi$ dominates $\pi$, 
the strict transform $E_A^\psi$ is a component of $E(\psi)$. We will 
distinguish the following two cases: 
\begin{enumerate}
\item $E_A^\psi = E$, that is,  the minimal resolution of $\nu_{E}$ equals 
the minimal embedded resolution of any curvetta of $E$. In this case we 
set $\ell_E :=0$. 
\item $E_A^\psi \ne E$, in this case the model $S(\psi)$ is obtained 
from $S(\pi)$ after blowing up $\ell_{E} > 0$ additional infinitely 
near points of $A$ over $O$.  We indicate in Figure \ref{fig:dualgraphs}
the shapes of the dual graph of $G(\psi)$ in both cases.
\end{enumerate}

We introduce the notion of generating sequence of a finite set of divisorial valuations 
according to Delgado, Campillo, Galindo and N\'u\~nez (see \cite{CG, DGN}), based upon the work of 
Spivakovsky who described the case of one valuation (see \cite{S}).

\medskip

Let $E_1, \dots E_s$ be exceptional prime divisors on a model $\psi$ of $(S,O)$.
We denote by $V = ( \nu_{E_1}, \dots, \nu_{E_s})$ be the associated tuple of divisorial valuations. 
Given a nonzero $h\in R$, we denote $V(h) =(\nu_{E_1} (h),\ldots,\nu_{E_s} (h))\in\Z^s_{\geq 0}$. 
For any element $\gamma\in \Z^s_{\geq 0}$, we consider the \textit{valuation ideal} 
$\mathcal{P}^V_{\gamma}$ of $R$ defined by 
$
\mathcal{P}^V_{\gamma} = \{ h \in {R} \mid V (h) \geq \gamma \}
\cup \set{0}$,
where $\leq$ is the partial ordering over $\Z^s$ given by $\delta \leq \gamma$ 
if and only if $\gamma-\delta\in\Z^s_{\geq 0}$.

\begin{definition}\label{def:suite_2} (See \cite{S,CG, DGN}).
The sequence $x_0, \dots, x_m \in m_R$ is a  \textit{generating sequence} of $V$  if for each 
$\gamma \in \Z^s_{\geq 0}$ the ideal $\mathcal{P}^V_{\gamma}$  is  
generated by those monomials in $x_0, \dots, x_m$  which belong to 
$\mathcal{P}^V_{\gamma}$. 
We say that it is \textit{minimal} if no proper subsequence of it is a generating sequence of
$V$.  
The {\it minimal resolution} of $V$ is a model $\psi$ of $(S,O)$ such that $\psi$ 
dominates the minimal resolution of $\nu_{E_i}$ for all $i\in\set{1,\ldots,s}$, and 
$\psi$ is minimal with this property. 
\end{definition}

Observe that $\mathcal{P}^V_{\underline 1}=m_R$ where 
$\underline{1} = (1, \dots, 1)$, and thus a generating sequence 
 $x_0, \dots, x_m$ generates the maximal ideal $m_R$. 
 If $s=1$, that is, when $V= ( \nu_{E}) $ we speak simply about generating sequences of the
divisorial valuation $\nu_{E}$.

\begin{remark}\label{rem:globalgivesgs}
A generating sequence of $V$ is also a generating sequence 
of $\nu_{E_i}$, for $1\leq i\leq s$. Indeed, if $p_i\colon\Z^s_{\geq0}\to\Z_{\geq0}$ 
denotes the canonical projection onto the $i^\text{th}$-coordinate and $\phi\in\Z_{\geq0}$, 
then $p_i^{-1}(\phi)$ has a unique minimal element $\gamma=\min p_i^{-1}(\phi)$ and then
$\mathcal{P}^V_{\gamma}=\mathcal{P}^{\nu_{E_i}}_{\phi}$.
\end{remark}

\begin{definition}  
Let $\psi$ be a model of $(S,O)$ different from the identity map and the blow up of $O$. Write $\mathcal{E}_\psi = \{ E_0, \dots, E_{s} \}$, where 
the indices are compatible with the order in which the divisors were created. 
We say that 
$L_0, \dots, L_s$ is 
a \emph{sequence of maximal contact curves} of $\psi$ if $L_i$ is a curvetta of $E_i$
at the model $\psi$, for all $i \in \{0, \dots,s \}$. If $\psi$ is either the identity map of $S$ or the blow up of $O$,  a \emph{sequence of maximal contact curves} of $\psi$ is $L_0, L_1$
where $(L_0, L_1)$ is a cross at $O$.
\end{definition}

The following result gives a geometric characterization of a 
minimal generating sequence of a finite set of divisorial valuations 
(see \cite[Section 8]{S} for the case of one valuation 
and \cite{DGN} for the general case). 

\begin{proposition}\label{gen-cur}  
Let $V$ be a tuple of divisorial valuations.
Let $L_0, \dots, L_{g(V)}$ be 
a sequence of maximal contact curves of 
the minimal resolution of $V$.
Take a defining function $x_i \in m_R$ of the branch $L_i$ for $i \in \{ 0, \dots, g(V)\}$. 
Then, $x_0, \dots, x_{g(V)}$ is a minimal 
generating sequence of $V$, and any minimal generating sequence of $V$ is of this form. 
\end{proposition}

It follows that if  $E $ is an exceptional prime divisor
 then any minimal generating sequence of $\nu_{E}$ has the same lenght, that is, it is of 
the form $x_0, x_1, \dots, x_{g(E)}$, for some $g(E) \geq 1$. 
\begin{example} \label{Spi2}
If $E=E_O$ then, for all $h\neq0$, $\nu_E(h)$ is the largest $n\in\Z_{\geq0}$ 
such that $h\in m_R^n$. Any minimal generating sequence of $\nu_E$ is of the 
form $x_0,x_1\in m_R$ where $(x_0,x_1)$ defines a local coordinate system on $S$. 
In this case $\ell_{E} =1$, and we have $g(E)=1$.
\end{example} 

\begin{proposition}[see \cite{S}] \label{GA2} 
Let us consider a minimal generating sequence $x_0,\ldots,x_{g(E)}$ of $\nu_E$. 
Denote $L_j = Z(x_j)$, for $j \in \{0, \ldots, g(E)\}$.
Let $A$ be any curvetta in $\mathcal{C}_E (\psi)$ and take a defining function $x_A$
of the branch $A$. 
Denote by $\g_0, \dots, \g_{g(A)}$ the minimal system of generators of the 
semigroup $\Gamma_A$  (see Notation \ref{not:GA}).
Then:
\[
 \nu_E (x_j ) =  \g_j, \mbox{ for  }  j=0, \dots,g(A), \mbox{ and }
\nu_E (x_{A} ) =  n_{g(A)} \g_{g(A)} + \ell_E .
\]
In addition, if $\ell_E =0$ then one has $g(E) = g(A)$, while 
if $\ell_E >0$ then $g(E) = g(A)+1$ and in this case
$L_{g(E)}$ is a curvetta in $\mathcal{C}_E (\psi)$.
\end{proposition}

\begin{notation} \label{GA3} 
If $\ell_E > 0$, we set  $n_{g(E)}:= 1$
and $\g_{g(E)} := n_{g(A)} \g_{g(A)} + \ell_E $.
We get the formulas 
$\nu_E (x_A) = n_{g(E)}  \g_{g(E)} $ and 
$\deg_{x_1} x_A = n_1 \dots n_{g(E)}$, which hold also when $\ell_E = 0$. 
\end{notation}

\begin{definition} \label{def:mccred}
Let $C= \cup_{j=1}^r  C_j$ be a reduced singular curve at $(S,O)$ and denote by
$\psi_C$ its minimal embedded resolution.
A sequence $L_0, \dots, L_{m}$ 
of maximal contact curves
of $\psi_C$ is 
\textit{generic} for $C$  if 
\begin{equation} \label{f:generic}
(C^{\psi_C} \cdot L_i^{\psi_C}) =0, \mbox{ for } 0 \leq i \leq m.
\end{equation}
\end{definition}

\medskip

Notice that condition \eqref{f:generic} holds if and only if $\psi_C$ is the minimal embedded resolution of $C \cup D$ where $D := \cup_{j=0}^m L_j$.

\medskip

The following lemma is a consequence of Proposition \ref{gen-cur} and Remark \ref{rem:globalgivesgs}. 

\begin{lemma} \label{lem:mgs}
Let $L_0, \dots, L_{m}$ be a sequence
of \textit{maximal contact curves} of $\psi_C$. 
If $j \in \{0, \dots, m\}$
then we can extract from 
$L_0, \dots, L_{m}$ 
a sequence of maximal contact curves of $\psi_{L_j}$.
\end{lemma}

\begin{proof}
If $L_j$ is smooth we may take $L_0, L_1$ defining a cross at $O$. 
If $L_j$ is singular, we denote by $E$ the exceptional prime divisor $E_{L_j}$ of 
the minimal resolution of $V$.
By Proposition \ref{gen-cur} and Remark \ref{rem:globalgivesgs}
one can extract  a minimal generating sequence $x_0, \dots, x_{g(E)}$ 
of  $\nu_{E}$ from 
a minimal generating sequence 
of $V$.
Let us denote by  $L_0, \dots, L_{g(E)}$ the branches defined by this sequence. 
Then, we apply Proposition \ref{gen-cur} again.
If $\ell_E = 0$ (resp. $\ell_E > 0$) then let 
$L_0, \dots, L_{g(E)}$ (resp. $L_0, \dots, L_{g(E)-1}$) be
a sequence of maximal contact curves of $\psi_{L_j}$, according to the cases considered in Proposition \ref{GA2}.
\end{proof}

\begin{example}\label{ex1}
Let $C$ be the plane curve defined by $f=f_1f_2$, where $f_1=Y^2-X^3-X^4$ and 
$f_2=Y^2-X^3-X^5$. Denote by $C_i$ the branch defined by $f_i$, for $i=1,2$. 
Let $\psi$ be the minimal embedded resolution of $C$. We obtain $\psi$, after 
blowing up the origin, by successively blowing up the point where the strict 
transforms of $C_1$ and $C_2$ intersect. The exceptional divisor $E(\psi)$ has 
5 irreducible components. We call them $D_0\ldots,D_4$, where the indices 
are compatible with the order in which the divisors where created. 
The divisors corresponding to the ends of $G(\psi)$ are $E_0:=D_0$, $E_1:=D_1$, 
and $E_2:=D_4$. 
The strict transforms $C_1^\psi$ and $C_2^\psi$ intersect $E_2$. 
Set $E:=E_2$. By construction, the minimal resolution of $\nu_E$ is $\psi$. 
We have $g(E)=2$ and $\ell_{E}=2$. 
For any $a\in\k$, the plane branch defined by $h_a=Y^2-X^3+aX^4$ is a 
curvetta of $E$. We may take as a minimal generating sequence of $\nu_E$ 
the functions $x_0:=X$, $x_1:=Y$ and $x_2:=h_a$. 

On the one hand, if $a=-1$ then $h_a=f_1$ 
and $(C_1^\psi\cdot Z(h_a)^\psi)=\infty$. On the other hand, if $a=0$ then 
$(C_2^\psi\cdot Z(h_a)^\psi)=1$. Hence, 
the functions 
$x_0$, $x_1$, and $x_2$ define a sequence of maximal contact curves of $\psi_C$ for any $a \in \k$, while 
the conditions 
\eqref{f:generic} 
are satisfied just for  $a\in\k\setminus\set{0,-1}$.
\end{example}

%-------------------------------------------------

\begin{example}\label{ex2}
Let $C$ be the plane curve defined by $f=f_1f_2$, where $f_1=(Y^2-X^3)^2-X^5Y$ 
and $f_2=(Y^2-X^3)^2-X^6Y$. Denote by $C_i$ the branch defined by $f_i$, for $i=1,2$. 
In order to obtain the minimal embedded resolution $\psi$ of $C$, we start with the 
minimal embedded resolution of $C_1$ and then we successively blow up the points 
at which the total transform of $C_2$ fails to have simple normal crossing support. 
The exceptional divisor $E(\psi)$ has 7 irreducible components. We call them 
$E_0\ldots,E_6$, where the indices are compatible with the order in which the 
divisors were created. With the notation of Definition \ref{def:sigma}, we have 
$E_{C_1}= E_4$ and $E_{C_2}= E_6$.

Let us set $V= (\nu_1,\nu_2)$, where $\nu_i$ is the divisorial valuation 
defined by $E_{C_i}$, for $i=1,2$. By construction, the minimal resolution 
of $V$ is $\psi$. The divisors $E_0$, $E_1$, and $E_5$ are those 
which correspond to ends of $G(\psi)$. 
 For any $a\in\k$, the branches defined by 
$x_0=X$, $x_1=Y$, and $x_2=Y^2-X^3+aX^4$ form a 
minimal generating sequence of $V$.

Consider now $C'$ be the curve defined by $f=f_1f_2f_3$  with $f_3=Y^2-X^3$. Then, for any $a\in\k^*$, $x_0=X$, 
$x_1=Y$, and $x_2=Y^2-X^3+aX^4$ define a sequence of maximal contact curves of 
$\psi_{C'}$.
\end{example}

%--------------------------------------------------------------------------------------------------------------

\subsubsection{Expansions in terms of generating sequences} \label{sec:expgensec}

We recall how to expand a function in terms of a minimal 
generating sequence  $x_0, \dots, x_{g(E)}$ of a divisorial valuation $\nu_{E}$, 
and how to determine from it the value of $\nu_{E}$ on this function. 
The elements $x_2, \dots, x_{g(E)}$ can be chosen as 
monic polynomials in the ring $\k [[x_0]] [x_1]$ (see  Remark \ref{WPT}).

\begin{proposition} (see \cite[Lemma 7.2]{PP})  \label{expansions-gen}
Let us consider positive integers $G\geq 1$, and $N_j>1$ for $j=1,\dots,G-1$, 
and let $f_j \in \k [[ x]] [y] $ be monic polynomials in $y$ with 
\[
\deg_y f_1 =1 \mbox{ and } 
\deg_y f_j = N_1 \dots N_{j-1},  \mbox{  for } j = 2, \dots, G.
\] 
Then, any polynomial $f \in \k [[ x]] [y]$ has a unique finite expansion in the form: 
\begin{equation} \label{exp}
f=\sum_{I =(i_1,\dots, i_{G})} c_{I }(x) \cdot f_1^{i_1} \cdots f_G^{i_G}, 
\end{equation} 
where $c_I (x)\in\k [[ x ]]$, $0\leq i_G \leq  \left\lfloor \frac{\deg_{y}(f)}{\deg_{y}(f_G)} \right\rfloor$, 
and $0\leq i_j<N_j$ for $j=1,\ldots,G-1$.
\end{proposition}

We call the expansion \eqref{exp} the \textit{$(x, f_1, \dots, f_G)$-adic expansion of $f$}. 

\medskip
Next we apply  Proposition \ref{expansions-gen} with respect to the minimal generating sequence 
$x_0, \dots, x_{g(E)}$ of $\nu_{E}$ taking into account the properties of the values 
$\nu_{E} (x_0) = \g_0, \dots, \nu_{E} (x_{g(E)} ) = \g_{g(E)}$
(see Proposition \ref{GA2}, Notation \ref{GA3} and Proposition \ref{GA}).
See also \cite[Rem. 8.16]{S} for details.

\begin{proposition} \label{expansions}
Let $A$ be a curvetta in $\mathcal{C}_E (\psi)$. The $(x_0, \dots, x_{g(E)})$-adic expansion 
of $x_{A}$ is of the form:
\begin{equation}\label{binomial-G}
x_{A}=x_{g(E)}^{n_{g(E)}} -\theta^A_{I_0} \cdot  x_{0}^{b_{0}}
x_{1}^{b_{1}}\cdots x_{g(E)-1}^{b_{g(E) -1}}
+p_A (x_0,\ldots,x_{g(E)}),
\end{equation}
where $\theta^A_{I_0} \in\k$, $I_0 = (b_{0},b_1, \dots, b_{g(E) -1}, 0)$, the equality 
\begin{equation} \label{binomial-G2}
{n_{g(E)}} \g_{g(E)} = {b_{0}} \g_0 + \cdots + {b_{g(E)-1}} \g_{g(E) -1},
\text{ with } 0 \leq b_{j} < n_j \text{ for } 1 \leq j <g(E)-1, 
\end{equation}
holds, where
$
p_A(x_0,\ldots,x_{g(E)})  = 
\sum_{I = ({i_0}, \dots, {i_{g(E)}})} \theta_I^{A} \cdot x_{0}^{i_0} \,  x_{1}^{i_1} \cdots x_{g(E)}^{i_{g(E)}}, 
$
with $\theta_I^A  \in \k$, and whenever $\theta_I^A \ne 0$ we have 
\begin{equation*}
{n_{g(E)}} \g_{g(E)} < i_0 \g_0 + \cdots + i_{g(E)} \g_{g(E)}, \text{ with } 0\leq i_j<n_j 
\text{ for } j=1,\ldots,g(E).
\end{equation*}
In addition, if $\ell_E =0$ then $\theta^A_{I_0} \in\k^*$.
\end{proposition}

We have also a similar result for the expansion of an element $x_{j+1}$ of 
the generating sequence in terms of the previous ones.

\begin{proposition}\label{prop:adicexpxj}
The $(x_0, \dots, x_j)$-adic expansion of $x_{j+1}$, for $1 \leq j \leq g(E) -1$, is of the form: 
\begin{equation}\label{binomial2}
x_{j+1}=x_{j}^{n_j}-\theta_j \cdot  x_{0}^{b_{0}^j}
x_{1}^{b_{1}^j}\cdots x_{j-1}^{b_{j-1}^j}
+p_j(x_0,\ldots,x_j),
\end{equation}
where $\theta_j\in\k^*$, the relation \eqref{rel-sem} holds, where 
$
p_j(x_0,\ldots,x_j) = 
\sum_{I = ({i_0^j}, \dots, {i_j^j})} \theta_I^{j} \cdot x_{0}^{i_0^j} \,  x_{1}^{i_1^j} \cdots x_{j}^{i_j^j} 
$
with $\theta_I^{j} \in \k$, and whenever $\theta_I^{j} \ne 0$ we have 
$n_j \g_j < \sum_{t=0}^j  i_t^j \g_t$ and $0\leq i_t^j<n_t$ for $t=1,\ldots,j$.
\end{proposition} 

\begin{remark}  \label{rem:expxj}
One can use the expansion of $h \in \k [[x_0]][x_1]$ in terms of the generating sequence 
of $\nu_{E}$ to determine the value of $\nu_{E} (h)$. Indeed, if $h \in \k[[ x_0]] [x_1]$ and if 
\[
h = \sum_{I = (i_1, \dots, i_{g(E)})}  c_{I }(x_0) \cdot x_1^{i_1} \cdots x_{g(E)}^{i_{g(E)}}, 
\]
is the $(x_0, \dots, x_{g(E)})$-adic expansion of $h$ then 
$
\nu _{E} (h) = \min_I  \{   \nu_{E} ( c_{I }(x_0) \cdot x_1^{i_1} \cdots x_{g(E)}^{i_{g(E)}} ) \}, 
$
where $ \nu_{E} ( c_{I }(x_0) \cdot x_1^{i_1} \cdots x_{g(E)}^{i_{g(E)}} )  ) = \ord_{x_0} (c_I) 
 \cdot  \g_0 + i_1 \g_1 + \cdots +  i_{g(E)} \g_{g(E)}$.
 This may be seen as a consequence of \cite[Th. 4.125]{Bu19}.
\end{remark}

%--------------------------------------------------------------------------------------------------------------
\subsection{The embedding defined by a sequence of maximal contact curves}\label{sec:planebranch}

We keep the notations of section \ref{sec:dival}. The curve $C$ is embedded in $S=\k^2$ by the map \eqref{one-2}. Recall that for any branch $C_j$ of $C$ we consider an arc $\eta_{C_j}$ on $k^2$ 
defined in terms of the normalization of $C_j$ and 
the given embedding of $C_j \subset \k^2$ (see Definition \ref{etaC}).

\begin{definition} \label{def:D}
Let $L_0, \dots, L_{m} $ be a sequence of maximal contact curves of $\psi_C$, 
which is generic for $C$ (see Definition \ref{def:mccred}). 
We set $D := \cup_{j=0}^m L_j$. 
Choose a defining function $x_i \in m_R$ of the branch $L_i$, for $i =0, \dots, m$. 
The sequence $x_0, \dots, x_m$ generates 
the maximal ideal of ${R}$. 
Hence it defines an embedding $\k^2 \hookrightarrow \k^{m+1}$ 
(see \eqref{re-emb}) corresponding to the surjective homomorphism
\begin{equation} \label{emb-plane}
\k [[ X_0, \dots, X_m ]] \longrightarrow R, \mbox{ such that } X_i \mapsto x_i, 
\mbox{ for } i=0, \dots, m. 
\end{equation}
 We call $I_S$ the kernel of the surjective homomorphism \eqref{emb-plane}, that is, the defining ideal of $S\subset\k^{m+1}$.  We get from this an embedding $C \hookrightarrow \k^{m+1}$ as in \eqref{re-emb-curv}.
\end{definition}
With the notation of Section \ref{sec:torenvironm}, for any branch $C_j$ of $C$ 
we have an arc $\eta_j$ on $k^{m+1}$  
that corresponds to the local homomorphism of $k$-algebras 
\[
\k [[X_0, \dots, X_m ]] \to \k [[t_j]], \mbox{ such that } X_i \mapsto x_i \circ  \eta_{C_j}, 
\mbox{ for } i=0, \dots, m. 
\]
By \eqref{f:generic}, we have that $x_0,\ldots,x_m\notin I_{C_j}$, 
and as a consequence the arc $\eta_j$ has its generic point in the torus of $k^{m+1}$.

\begin{proposition} \label{prop-semi}
Keep the previous notation. We have:
\begin{enumerate}
\item\label{enum:generators} If $C_j$ is a branch of $C$, then the positive integers 
$\ord_{t_j} (x_0 \circ \eta_{C_j} ), \ldots, \ord_{t_j} (x_m \circ  \eta_{C_j} )$, 
generate the semigroup $\Gamma_{C_j}$.
\item\label{enum:samediv} 
Let $C_j$ and $C_{j'}$ be two branches of $C$. If $\ord(\eta_j)=\ord(\eta_{j'})$, then $C^{\psi_C}_j$ and $C^{\psi_C}_{j'}$ intersect the same irreducible component of $E(\psi_C)$.
\end{enumerate}
\end{proposition}
\begin{proof}
Denote the minimal embedded resolution $\psi_C$ of $C$ simply by $\psi$. 
The branch $C_j$ belongs to $\mathcal{C}_{E_{C_j}} (\psi)$, where we recall that 
$E_{C_{j}}$ is the irreducible component of $E(\psi)$ which intersects $C^\psi_j$. 
By Remark \ref{rem:val-int} and \eqref{f:generic}, we deduce that 
\begin{equation}\label{eq:ordervect}
\nu_{E_{C_j}} (x_i) = (L_i \cdot C_j) = \ord_{t_j} (x_i \circ \eta_{C_j} ), 
\mbox{ for } 0 \leq i \leq m.    
\end{equation}
By Proposition \ref{gen-cur}, $x_0, \dots, x_{m}$ is a generating sequence of the 
divisorial valuation $\nu_{E_{C_j}}$, thus this sequence contains a minimal 
generating sequence of $\nu_{E_{C_j}}$. By Proposition \ref{GA2}, the set 
$\set{ \ord_{t_j} (x_i \circ \eta_{C_j} ) \mid 0\leq i\leq m}$  
contains the minimal generating system of $\Gamma_{C_j}$. 
This proves \eqref{enum:generators}.

If $C_j$ and $C_{j'}$ are two branches of $C$ and $\ord(\eta_j)=\ord(\eta _{j'})$, then 
$\nu_{E_{C_j}}(x_i)=\nu_{E_{C_{j'}}}(x_i)$ for 
$i=0,\ldots,m$, by \eqref{eq:ordervect}. Since $x_0,\ldots,x_m$ is a generating sequence 
of both valuations $\nu_{E_{C_j}}$ and $\nu_{E_{C_{j'}}}$, and their values on this sequence coincide, 
the valuation ideals must be equal. 
 This implies that $\nu_{E_{C_j}}=\nu_{E_{C_{j'}}}$ and thus 
 $E_{C_j} = E_{C_{j'}}$.
\end{proof}

The following proposition provides a geometric interpretation of  a coefficient appearing on a
$(x_0, \dots, x_{g(E)})$-adic expansion of the defining function $x_A$, 
of a curvetta $A \in \mathcal{C}_E (\psi)$.

\begin{proposition} \label{thetas}
Let $\psi\colon (S(\psi), E(\psi)) \to (S,O)$ be the minimal resolution of $\nu_{E}$ 
and let $A$ and $B$ be two curvettas in $\mathcal{C}_E (\psi)$. 
Let $\theta^A_{I_0}$ and $\theta^B_{I_0}$ be the coefficients of the term 
$x_{0}^{b_{0}}x_{1}^{b_{1}}\cdots x_{g(E)-1}^{b_{g(E) -1}}$ in the expansions 
\eqref{binomial-G} of $x_A$ and $x_B$, respectively. 
The following are equivalent:
\begin{enumerate}
\item\label{enum:thetas1} 
The strict transforms of $A$ and $B$ on $S(\psi)$ intersect 
$E$ at different points.
\item\label{enum:thetas2} The coefficients $\theta^A_{I_0}$ and 
$\theta^B_{I_0}$ are different.
\end{enumerate}
\end{proposition}

\begin{proof}
$\bullet$ \textbf{The case $\ell_{E} > 0$.} 
In this case $L_{g(E)} = Z(x_{g(E)} )$ is a curvetta of $E$. Take local coordinates 
$(w, s)$ at the point of $S(\psi)$ where $L_{g(E)}$ and $E$ meet, such that 
$Z(w) = E$, and $x_{g(E)}\circ\psi=w^{\nu_{E}(x_{g(E)})}\cdot s$, 
that is, $Z(s) = L_{g(E)}^\psi$. Recall that $\nu_E(x_{g(E)})=\g_{g(E)}$ (see Notation \ref{GA3}). 
Let us compute the total transform of $A$. 
First take $j\in\set{0,\ldots,g(E)-1}$. 
By construction, the exceptional prime divisor $E$ and 
the strict transform $L_j^\psi$ do not intersect. 

Recall that $\nu_E(x_j)=\g_j$.
Therefore, there exists $\epsilon_{0,j}\in\k^*$ such that
$
x_j \circ \psi = w^{\g_j}(\epsilon_{0,j}+w(\cdots))$.
Set $\epsilon_{I_0}:=\prod_{j=0}^{g(E)-1}\epsilon_{0,j}^{b_j}\in\k^*$.
From the expansion \eqref{binomial-G} of $x_A$ we obtain the following:
\[
x_{A} \circ \psi =w^{\g_{g(E)}} \left( s -\theta^A_{I_0} \epsilon_{I_0} + w (\cdots) \right).
\]
Then, the strict transform of $A$ is $A^\psi=Z(s-\theta^A_{I_0}\epsilon_{I_0}+w(\cdots))$ 
and its intersection with $E$ is the point of coordinates $(0,s(A))$ where 
$s(A)=\theta^A_{I_0}\epsilon_{I_0}$. If we apply this to $B$ we get that 
$s(B)=\theta_{I_0}^B\epsilon_{I_0}$. 
Since $\epsilon_{I_0} \ne 0$, it follows that \eqref{enum:thetas1} 
and \eqref{enum:thetas2} are equivalent.

$\bullet$ \textbf{The case $\ell_{E} =0$.}  
 In what follows, we use $g$ instead of $g(E)$. In this case $L_{g} = Z(x_{g} )$ is not a curvetta of $E$. The model $\psi$ dominates 
the minimal embedded resolution 
$\varphi\colon (S(\varphi), E(\varphi)) \to (S,O)$ of $L_{g}$.

We address first the case where $g>1$. We denote by $R_{g}$ the exceptional 
prime divisor of $\varphi$ which intersects $L_{g} ^\varphi$. Take local coordinates 
$(w, s)$ at the point of $S(\varphi)$ where $L_{g}^\varphi$ and $R_{g}$ meet, 
such that $Z(w) = R_{g}$, and $x_{g} \circ \varphi= w^{\nu_{R_{g}} (x_{g} ) }\cdot s,$ 
and then, $Z(s) = L_{g}^\varphi$. According to Proposition \ref{gen-cur}, 
the minimal generating system of the semigroup of $L_g$ is 
$\nu_{R_g} (x_0),\nu_{R_g} (x_1),\ldots,\nu_{R_g} (x_{g-1})$. 
By \cite[Proposition 5.4]{GBP},
we have that 
$\nu_{R_g} (x_j) = {\g_j}/{n_g}$
for $j\in \{0, \dots, g-1 \}$. 
Using this and Proposition \ref{GA2}, 
we get that $\nu_{R_g} (x_g) =  \frac{n_{g-1}}{n_g} \g_{g-1}$ and thus 
$x_g\circ\varphi=w^{n_{g-1}\g_{g-1}/n_g}\cdot s$. 
As in the case $\ell_{E}>0$, for $j=0,\ldots,g-1$, we have $
x_j \circ \varphi = w^{\g_j/n_g}(\delta_{0,j}+w(\cdots))$  with 
$ \delta_{0,j}\in\k^*$.
Set $\delta_{I_0}=\prod_{j=0}^{g-1}\delta_{0,j}^{b_j}\in\k^*$. 
Taking into account that $n_{g-1}\g_{g-1}<\g_g$ 
(see \eqref{rel-sem}) and denoting by $m_g$ the difference 
$\g_g-n_{g-1}\g_{g-1}$, from the expansion \eqref{binomial-G} of $x_A$ we obtain 
that
\[
x_{A}\circ \varphi =
w^{n_{g-1}\g_{g-1}}\left(s^{n_g}-w^{m_g}\left(\theta^A_{I_0}\delta_{I_0}+w(\cdots)\right)\right).
\]
Hence $A^\varphi=Z(h_A)$, where $h_A=s^{n_g}-w^{m_g}
\left(\theta^A_{I_0}\delta_{I_0}+w(\cdots)\right)\in\k[[w,s]]$. Notice that the integers $n_g$ and $m_g$ 
must be coprime (otherwise $A^\varphi$ would not be a branch) and the compact edge of the Newton polygon of 
$h_A$ is the line segment that joins $(0,n_g)$ and $(m_g,0)$. Now, let us consider the fan 
${\Sigma_0}:={\Sigma}_{0,2}$ (see Notation \ref{quadrant}) and the vector $u=(n_g,m_g)$. 
Let $\pi:S(\pi) \to S(\varphi)$ be the toric modification, defined with respect to the local 
coordinates $(w, s)$, associated to the minimal regular subdivision $\Sigma$ of ${\Sigma_0}$ 
containing the ray $\R_{\geq0}u$. 
This morphism is the minimal embedded resolution of $A^\varphi$ and $\psi = \varphi \circ \pi$. 
We refer to \cite{GP} for details. 
In the chart corresponding to the cone $\R_{\geq0}u+\R_{\geq0}v\in\Sigma$, where $v=(c,d)$ 
and $cm_g-dn_g=1$, this morphism is given by 
$w =  w_1^c s_1^{n_g}$, and  $s  =  w_1^d s_1^{m_g}$.

We get that $h_A(w_1^c s_1^{n_g},w_1^d s_1^{m_g})=w_1^{dn_g} s_1^{n_gm_g}
\left(1-\theta^A_{I_0}\delta_{I_0}w_1+ s_1\left(\cdots\right)\right)$. 
The exceptional divisor $E$ on this chart is $Z(s_1)$.
Thus, $A^\psi$ intersects $E$ at the point with coordinates $(w_1, s_1) = (w_1(A),0)$, 
where $w_1(A)=1/\theta^A_{I_0}\delta_{I_0}$. If we apply the previous construction 
to $B\in\mathcal{C}_{E}(\psi)$, 
we get that $w_1(B)$ equals $1/\theta^B_{I_0}\delta_{I_0}$. 
It follows that \eqref{enum:thetas1} and \eqref{enum:thetas2} are equivalent.

Assume now that $g=1$. 
The expansion \eqref{binomial-G} of $x_A$ is of the form \[x_A=x_1^{\g_0}-\theta_{I_0}^A x_0^{\g_1}+
\sum_{\g_0 i_0+\g_1 i_1>\g_0\g_1}{\theta^A_{I}x_0^{i_0}x_1^{i_1}}.\]
Let ${\Sigma_0}$ be as above and $u=(\g_0,\g_1)$. Then, $\psi$ is the toric modification, 
defined with respect to the local coordinates $(x_0,x_1)$, associated to 
the minimal regular subdivision of ${\Sigma_0}$ containing $\R_{\geq0}u$. 
Choosing the chart  as above, we see that \eqref{enum:thetas1} 
and \eqref{enum:thetas2} are also equivalent in this case.
\end{proof}

We prove now the announced generalization a theorem of Goldin and Teissier:

\begin{theorem} \label{th:resgenplane}
Let $C= \cup_{j=1}^r  C_j$ be a  reduced singular plane curve with $r\geq1$ branches embedded in $\k^2$ by the map \eqref{one-2}. Denote by $\psi_C$ the minimal embedded resolution of $C$. Let $x_0, \dots, x_m$ be defining functions of a sequence  $L_0, \dots, L_m$ of 
 maximal contact curves  of $\psi_C$ which is generic for $C$. 
Let $C \hookrightarrow \k^{m+1}$ be the embedding defined by 
the functions $x_0, \dots, x_{m}$. 
Then, any regular subdivision $\Sigma$ of the fan ${\Sigma}_{0,m+1}$ 
such that $\Sigma$ contains the ray $\rho_j$ spanned by $\ord(\eta_j)$, 
for $1\leq j \leq r$, defines a toric embedded resolution of $C$. 
\end{theorem}

\begin{proof}
The functions $x_0,\ldots,x_m$ satisfy conditions \eqref{enumresgen:me}, 
\eqref{enumresgen:gptorus}, and \eqref{enumresgen:primitive} of Corollary 
\ref{cor:resgen} (see the discussion at the beginning of section \ref{sec:planebranch} 
and Proposition \ref{prop-semi}.\eqref{enum:generators}). 
Therefore, the result follows directly from Corollary \ref{cor:resgen} 
if $r=1$, or $r>1$ and the order vectors $\ord(\eta_j)$, for $1\leq j\leq n$, 
are pairwise different  (see Corollary \ref{enum:diffvect}). 
So let us assume that $\ord(\eta_j)=\ord(\eta_{j'})$ 
for some $1\leq j<j'\leq r$. It was shown in the proof of Corollary \ref{cor:resgen} 
that the strict transforms $\tilde C_j$ and $\tilde C_{j'}$ are smooth and 
transversal to the orbit $O(\rho)_\Sigma$, where $\rho:=\rho_j=\rho_{j'}$. 
Next we show that the special points of $\tilde C_j$ and $\tilde C_{j'}$ are 
different. This is enough to end the proof.

Denote by $\tilde\eta_j$ (resp. $\tilde\eta_{j'}$) the lifting of the arc $\eta_j$ 
(resp. $\eta_{j'}$) to $X_\Sigma$. Take $\sigma\in\Sigma$ spanned by a basis 
$v_0 \in \rho,v_1,\ldots,v_m$ of the lattice $\Z^{m+1}$. The affine toric variety $X_\sigma$ 
is isomorphic to $\k^{m+1}$ with coordinates $(U_0,\ldots,U_m)$, 
where $U_i=\chi^{\check v_i}$ for $i=0,\ldots,m$. 
The parametrization $\tilde\eta_j$ of $\tilde C_j$, in terms of these coordinates, 
is $(U_0\circ\tilde\eta_j,\ldots,U_m\circ\tilde\eta_j)$, and we have an analogous 
expression for the parametrization $\tilde\eta_{j'}$ of $\tilde C_{j'}$. 
Since $\ord(\eta_j)=\ord(\eta_{j'})$, we have the identities for $0\leq i\leq m$:
\begin{align*}
X_i \circ \eta_j    = x_i \circ \eta_{C_j}     = d_i t_j^{\alpha_i}+\cdots,\quad & \quad 
X_i \circ \eta_{j'} = x_i \circ \eta_{C_{j'}}  = d'_i t_{j'}^{\alpha_i}+\cdots,
\end{align*}
where $d_i,d'_i\in\k^*$, we have omitted higher order terms, 
and $(\alpha_0,\ldots,\alpha_{m+1})$ are the coordinates of $\ord(\eta_j)$ 
with respect to the canonical basis of $\Z^{m+1}$.

Let us assume that the special points of $\tilde C_j$ and $\tilde C_{j'}$ coincide 
and proceed by contradiction. By assumption, there exist 
$ c_0, c_0', c_1\ldots, c_m\in\k^*$ such that:
\begin{align*}
U_0 \circ \tilde\eta_j &= c_0 t_j+\cdots,  
&  U_i \circ \tilde\eta_j &= c_i+\cdots,\text{ for }1\leq i\leq m;\\
U_0 \circ \tilde\eta_{j'} &= c_0' t_{j'}+\cdots,
&  U_i \circ \tilde\eta_{j'} &= c_i+\cdots,\text{ for }1\leq i\leq m,
\end{align*}
where we have omitted higher order terms (see \eqref{cond}).

For $1\leq i\leq m$, denote by $(a_{0,i}, \dots, a_{m,i}) $ the coordinates of $v_i$ 
with respect to the canonical basis of $\Z^{m+1} =N$. Taking into account 
the monomial expression of the toric modification in the chart of $X_\Sigma$ 
defined by the cone $\sigma$, we obtain: 
\begin{align*}
d_i =  c_0^{\alpha_i} c_1^{a_{i,1}} \cdots  c_m^{a_{i,m}} , 
\quad & 
d_i' = (  c_0')^{\alpha_i}  c_1^{a_{i,1}} \cdots  c_m^{a_{i,m}}, 
\mbox{ for } 0\leq i \leq m.    
\end{align*}
Set $\kappa_i := c_1^{a_{i,1}} \ldots  c_m^{a_{i,m}}\in\k^*$, for $i=0, \dots, m$. 
Without loss of generality, we may assume that $x_0,\ldots,x_g$ is a minimal generating 
sequence of $\nu:=\nu_{E_{C_j}}$ as in Proposition \ref{gen-cur}. 
In comparison with Notation \ref{GA3}, we have that $\alpha_i=\bar\beta_i$, for $i=0,\ldots,g$. 
We consider a $(x_0,\dots,x_g)$-adic expansion \eqref{binomial-G} of a defining 
function $x_{C_j}$ of $C_j$. Substituting $x_i$ by the series $x_i \circ \eta_{C_j}$, 
for $i=0,\dots,g$, in the expansion \eqref{binomial-G}, provides an expression
\[
0 =  x_{C_j} \circ \eta_{C_j}
=(d_{g}^{n_{g}} -\theta^{C_j}_{I_0} \cdot  d_{0}^{b_{0}}
d_{1}^{b_{1}}\cdots d_{g-1}^{b_{g -1}}) t_j^{n_{g} \bar \beta_g} + \cdots,
\]
where we have omitted the terms of higher order on the right side. It follows that:
\begin{equation} \label{ell2}
d_{g}^{n_{g}} -\theta^{C_j}_{I_0} \cdot  d_{0}^{b_{0}}
d_{1}^{b_{1}}\cdots d_{g-1}^{b_{g -1}} = 0. 
\end{equation}
Now substitute $d_i= \kappa_i c_0^{\alpha_i}$ in \eqref{ell2}.
Taking into account the condition \eqref{binomial-G2} we get the equality:
\[
( \kappa_{g}^{n_{g}} -\theta^{C_j}_{I_0} \cdot  \kappa_{0}^{b_{0}}
\kappa_{1}^{b_{1}}\cdots \kappa_{g-1}^{b_{g -1}} ) c_0^{n_g \bar{\beta}_g} = 0. 
\]
Since $c_0 \ne 0$, we conclude that $\theta^{C_j}_{I_0}= \kappa_{g}^{n_{g}} / \kappa_{0}^{b_{0}}
\kappa_{1}^{b_{1}}\cdots \kappa_{g-1}^{b_{g -1}}$. 
According to Proposition \ref{prop-semi}.\eqref{enum:samediv}, we have that $E_{C_j}=E_{C_{j'}}$. 
Hence, by the same argument applied to $C_{j'}$, it follows that $\theta^{C_j}_{I_0}= \theta^{C_{j'}}_{I_0}$. 
This gives a contradiction with Proposition \ref{thetas}.
\end{proof}

\begin{example}
Let $C$ be the curve of Example \ref{ex1}. The branches $L_0,L_1,L_2$ defined by $x_0=X$, 
$x_1=Y$, and $x_2=Y^2-X^3+X^4$, respectively, form a sequence of maximal contact curves 
of $\psi_C$ which satisfies \eqref{f:generic}. Let $C\hookrightarrow k^3$ be the embedding defined 
by this sequence. With the notation of Example \ref{ex1}, for $1\leq i\leq 2$, we have 
$(L_0,C_i)=2$, $(L_1,C_i)=3$, and $(L_2,C_i)=8$. Hence $\ord(\eta_1)=\ord(\eta_2)=(2,3,8)$. 
By Theorem \ref{th:resgenplane}, any regular subdivision of ${\Sigma}_{0,3}$ containing 
the ray spanned by $v_3:=(2,3,8)$ defines a toric embedded resolution of $C$. Let us 
study the strict transforms of $C\subset\k^3$ and $S=\k^2\subset\k^3$ under the monomial 
map $\pi:=\pi^\sigma_{\sigma_{0,3}}$ associated to the regular cone 
$\sigma\subset\sigma_{0,3}$ generated by the vectors 
$v_1=(1,1,2)$, $v_2=(2,3,7)$, and $v_3$. We take coordinates in $X_{{\Sigma}_{0,3}}=k^3$, 
and also in $X_\sigma$, as in Example \ref{ex:stplanesing}. 

On the one hand, the surface $S\subset\k^3$ has equation 
$X_3-(X_2^2-X_1^3+X_1^4)=0$. On the other hand, the points of $C\subset\k^3$ 
satisfy the previous equation and $(X_3-2X_1^4)(X_3-X_1^4-X_1^5)=0$. 
Thus the points in the strict transform $\tilde C$ of $C$ under 
$\pi$ satisfy the following two equations:
\begin{align}
    1-U_1-U_2U_3^2+U_1^2U_2^2U_3^2=0,\label{eq:stplane2}\\
    (1-2U_1^2U_2)(1-U_1^2U_2-U_1^3U_2^3U_3^2)=0.\nonumber
\end{align}
One checks easily that the strict transform $\tilde S$ of $S$, which is given by equation 
\eqref{eq:stplane2}, is smooth. The points $(1,t,0)\in\tilde S$, for $t\in\k$, 
are those in the intersection with the exceptional locus of $\pi$. 
For $t=1/2$ (resp. $t=1$) we find the intersection with $\tilde C_1$ (resp. $\tilde C_2$).
\end{example}

%-------------------------------------------------------------------------------------------------------------

\section{Embedded resolution and local tropicalization}\label{sec:TEpairs}

In this section we consider a singular curve $C \subset S = \k^2$ and the embedding $S= \k^2  \subset \k^{m+1}$ 
defined by the map \eqref{emb-plane}. 
The goal of this section is to prove that there exists a regular subdivision $\Sigma$ of the positive quadrant 
${\Sigma}_{0, m+1}$
such that the restriction of $\pi_{\Sigma}$ to the strict transform of $S$ 
is the minimal embedded resolution of the given plane curve singularity $C\subset S$. 
In order to prove this we describe the \textit{local tropicalization} of $S$ 
associated to this embedding.  

\medskip 

Throughout this section we keep the notation introduced at the beginning 
of section \ref{sec:planebranch}. In particular, 
$L_0=Z(x_0), \dots, L_{m}=Z(x_m) $ is a sequence of maximal contact curves of $\psi_C$ which is generic for $C$. 
We will consider below the auxiliary curve $D = \cup_{j=0}^m L_j$.

\subsection{The local tropicalization of $S$ defined by the maximal contact curves}\label{sec:localtrop}

We start by recalling the definition of a semivaluation of a local domain with values in the ordered set 
$\overline{\R}_{\geq 0}:=\R_{\geq0} \cup \{ \infty \}$. 

\begin{definition} \label{valdef}
Let $R$ be a local domain containing the field $\k$. 
A \emph{semivaluation} of $R$ is a function 
$\nu\colon R\rightarrow \overline\R_{\geq 0}$ such that:
\begin{enumerate}
    \item \label{add} $\nu(fg)=\nu(f) + \nu(g)$ for all $f,g \in R$; 

    \item \label{ineq} $\nu(f+g) \geq \min(\nu(f), \nu(g))$  for all $f,g \in R$;

    \item \label{const} $\nu(0)=\infty$ and $\nu(\lambda)=0$ for all $\lambda\in \k^*$.
\end{enumerate} 
The semivaluation $\nu$ is a \emph{valuation} if it 
takes the value $\infty$ only at $0$, and it is \emph{trivial} if 
$\nu(R)=\set{0,\infty}$.
\end{definition}

Let us come back to the case $R=\k [[ X, Y]]$. 
If $A$ is a plane branch defined by a function $x_A \in R$, then we denote
$\nu(A) := \nu(x_A)$ for every semivaluation $\nu$ of $R$. Associated to a plane branch $A$, 
we have the \emph{intersection semivaluation} $I^A$,  defined by
$I^A (h) = (A \cdot Z(h))$ for any $h \in R$,
the \emph{vanishing order semivaluation} $\ord^A$, whose value 
on a function $h \in R$ is the order of vanishing of $h$ along the branch $A$, and the \emph{trivial valuation} $\mathrm{triv}^A$ which takes the value $\infty$ on any function vanishing on $A$ and zero otherwise. 
The trivial valuations $\mathrm{triv}^O$ and  $\mathrm{triv}^S$ associated to the point $O$ and the surface $S$ respectively, are defined analogously to $\mathrm{triv}^A$. The \emph{center} of a semivaluation $\nu$ is the algebroid subvariety defined by 
the prime ideal $\nu^{-1} ( (0, \infty])$. The center of $\nu$ is $S$ if $\nu = \mathrm{triv}^S$, otherwise the center is $O$ or a branch.
If $\nu$ is non-trivial and its center is a branch 
$A$ then $\nu$ is proportional to $\ord^A$.

\medskip
The \emph{semivaluation space} $\cV$ of $R$ (or of the smooth surface 
germ $S$) is the set of semivaluations of $R$, endowed 
with the topology of pointwise convergence, that is, with the restriction 
of the product topology of $[0,\infty]^{R}$. 
We denote by $\cV_{L_0} = \{ \nu \in \cV \mid \nu (x_0)  =1 \}$ 
the subspace of \emph{normalized semivaluations with respect to $L_0$} and by $\cV^*$ the subspace of non-trivial semivaluations of $\cV$. 
The semivaluation space $\cV$ has been intensively studied after Favre 
and Jonsson's book  \cite{FJ}, see also \cite{Jon15, GGP19}. 
The normalized semivaluation space $\cV_{L_0}$, which is
also called the \emph{(relative) valuative tree},
is an \emph{$\R$-tree rooted} at $\ord^{L_0}$ (see \cite[section 3.9]{FJ} and \cite{GGP19}). 
One has the following property:

\begin{lemma} (see \cite[section 6]{GGP19})  \label{rem:intsemivalLj}  
Let  $\nu\in\cV_{L_0}$ and let $A$ be a plane branch different from $L_0$.  
If  $\nu(A)=\infty$ then we have that  $\nu = I_{L_0}^{A} := \frac{1}{(L_0\cdot A)} \cdot I^{A}$. If $\nu(A)=0$ then  $\nu $ is equal to $\ord^{L_0}$.
\end{lemma}
 
We consider the following notion of \emph{local tropicalization} introduced in  \cite{PS13} (see also \cite{St}).

\begin{definition}
The \emph{local tropicalization} $\Trop_{\geq 0} (I)$ of an ideal 
$I\subset\k [[ X_0, \dots, X_m ]]$ is the set of tuples 
$(\mu (X_0), \dots, \mu (X_m)) \in \overline\R_{\geq 0}^{m+1},$
for $\mu$ running through the semivaluations $\mu$ 
of $\k [[X_0, \dots, X_m]]$ such that 
$\mu_{| I} = \infty $. 
\end{definition}

We have also a notion of tropicalization map which is defined in terms of  $x_0, \dots, x_m \in m_R$:

\begin{definition}
The
\emph{tropicalization map} defined by  the tuple $(x_0, \dots, x_m)$ is 
\[
\begin{array}{crlc}
        \trop \colon & \cV & \rightarrow & \overline\R_{\geq 0}^{m+1} \\
        & \nu & \mapsto & (\nu (x_0), \dots, \nu (x_m) )
\end{array}
\]
\end{definition}

\begin{lemma}\label{lem:tropval}
With the previous notation, we have $\Trop_{\geq 0} (I_S) = \trop ( \cV)$. 
\end{lemma}
\begin{proof}
If $\mu$ is a semivaluation of $\k [[X_0, \dots, X_m]]$ verifying that 
$\mu_{| I_S} = \infty $, then 
one has a semivaluation $\bar{\mu}$ of the quotient ring $\k [[ X_0, \dots, X_m ]] / I_S =  R$ 
by setting 
$\bar{\mu} ( \phi + I_S ) := \mu( \phi)$ for every  $\phi \in \k [[ X_0, \dots, X_m ]]$, and  every semivaluation of $R$ can be obtained in this way.

\end{proof}

The aim of this section is to describe the \textit{finite local 
tropicalization} of $I_S$, which is defined as
$\Trop_{\geq 0} (I_S) \cap \R_{\geq0}^{m+1}$. By Lemma \ref{lem:tropval}, 
we have the equality 
$
\Trop_{\geq 0} (I_S) \cap \R_{\geq0}^{m+1}=\trop(\mathcal{U} )$,
where $\mathcal{U}$ denotes the set of semivaluations $\nu \in \cV$ such 
that $\nu (x_j) < \infty$, for $j =0, \dots, m$. 

\medskip

In what follows we denote by $e_0,\ldots,e_m$ the canonical basis of $\R^{m+1}$. Given $w\in\R^{m+1}$, we write $w=\sum_{0\leq j\leq m}w_j e_j$. 
If $X \subset \R^{m+1}$, we denote by $\R_{>0} X $ the set 
$\{ a x \mid a \in \R_{>0}, x \in X \}$, and 
similarly $\R_{\geq 0} X := \R_{>0} X \cup \{ 0 \}$.

\begin{lemma} \label{lem:tropU}
We have the following equalities:
\begin{enumerate}[label=(\roman*)]
\item\label{item:tropU1} $\trop(\mathcal{U} ) 
=\left(\trop(\mathcal{U} ) \cap \R_{>0}^{m+1}\right) 
\cup \R_{\geq0} \, e_0 \cup \cdots \cup \R_{\geq0} \, e_m$.
\item\label{item:tropU2} 
$\trop(\mathcal{U} ) \cap \R_{>0}^{m+1}=
\R_{>0}
\left(\trop\left(\cV_{L_0}\setminus
\set{\ord^{L_0},I_{L_0}^{L_1},\ldots,I_{L_0}^{L_m}}\right)\right)$.
\end{enumerate}
\end{lemma}

\begin{proof}
(i) 
If $\nu \in \mathcal{U}$ is a
trivial semivaluation then $\trop(\nu)= 0$. 
Assume that $\nu \in \mathcal{U}$ is non-trivial. If $\nu (x_j) \ne 0 $ for all $j \in \{ 0, \dots, m\}$,
then $\trop (\nu) \in \R_{>0}^{m+1}$. If $\nu (x_i) =  0 $ for some $i \in \{ 0, \dots, m\}$, then 
$\nu$ is not centered at $O$. Since $\nu$ is non-trivial it must be proportional to 
$\ord^A$ for some branch $A$.  If $A=L_j$ for some 
$0\leq j\leq m$ then the vector $\trop(\nu)$ belongs to $\R_{>0}\, e_j$, 
otherwise $\trop(\nu) =0$.

(ii) 
Let us check the inclusion $\subset$. If $\nu \in \mathcal{U}$ and $\trop(\nu) \in \R_{>0}^{m+1}$, then we can write 
$\nu = \nu(x_0) \frac{ \nu}{\nu(x_0)} $, where $\frac{ \nu}{\nu(x_0)}$ belongs to 
$\cV_{L_0}$. In addition, if $\nu' \in \{ \ord^{L_0},I_{L_0}^{L_1},\ldots,I_{L_0}^{L_m} \}$ then
$\nu \ne \nu'$ since $\nu(x_j) \ne \nu'(x_j) $ for some $j \in \{0, \dots, m \}$.

Next we check the inclusion $\supset$. If $\nu \in \trop (\cV_{L_0} \setminus
\{\ord^{L_0},I_{L_0}^{L_1},\ldots,I_{L_0}^{L_m}\})$ then
$\nu(x_0) = 1$ and for $j \in \{1, \dots, m \}$ we have that $\nu(x_j) \ne 0 $
since $\nu \ne \ord^{L_0}$, while $\nu(x_j) \ne \infty$ since $\nu \ne I_{L_0}^{L_j}$
(see Lemma \ref{rem:intsemivalLj}).
This implies that $\trop (\nu) $ belongs to $\trop(\mathcal{U} ) \cap \R_{>0}^{m+1}$.
\end{proof}

We will study the set $\trop(\mathcal{U})$ by using properties of the 
semivaluation spaces $\cV$ and $\mathcal V_{L_0}$, in connexion with the 
Eggers-Wall trees which we introduce breafly as follows.

\subsubsection{Eggers-Wall trees}\label{sec:EW}

We consider a plane curve singularity $C = \cup_{j=1}^r C_j$ with $r\geq1$ branches and
a smooth branch $L_0$. We describe some features of the Eggers-Wall tree associated to $C$ with respect to $L_0$, see  \cite{GGP19} for precise definitions and details.

\medskip

Let us introduce first some basic notation about trees (see \cite[section 2]{GGP19}).
If $a, b$ and $c$ are points in a tree $T$, 
we denote by $\langle a, b ,c \rangle$ 
the unique point of intersection of the segments 
$[a, b]$, $[b,c]$ and $[a,c]$. We call 
the point $\langle a, b ,c \rangle$
the \emph{center} of the \emph{tripod}
$\{a, b, c \}$. If the tree $T$ is rooted at $r$, 
we set 
$a\preceq_{r}b$ if $[r,a]\subseteq [r,b]$ for $a, b\in T$. 

\medskip

The Eggers-Wall tree $\Theta_{L_0} (C)$ is a tree endowed 
with an \emph{exponent function} $\ex_{L_0}$ and an \emph{index function} $\de_{L_0}$:
\[
\ex_{L_0}\colon\Theta_{L_0} (C)\to\overline\R_{\geq 0}, \quad 
\de_{L_0}\colon\Theta_{L_0} (C)\to\Z_{>0}.
\]
If $\k$ is of characteristic zero, this tree 
and these functions are usually defined in terms of the Puiseux 
expansions with respect to $x_0$ of the branches of $C$ different from $L_0$. 
If $\k$ is of arbitrary characteristic 
$\Theta_{L_0} (C)$ may be defined as the convex hull of the points labelled by $L_0$ and the branches 
of $C$  in the  \emph{fan tree} associated to a 
toroidal pseudo-resolution of $C$ with respect to $L_0$. 
The functions $\ex_{L_0}$ and $\de_{L_0}$
are determined then by the \emph{slope function} on the 
fan tree (see \cite[Remark 5.37]{GGP19b}).

\medskip

The  \emph{marked points} of $\Theta_{L_0} (C)$ are its ends (i.e.\ the root labelled by $L_0$ 
and the leaves labelled by $C_j$, $1\leq j \leq r$), its ramification points, and the points of discontinuity 
of the index function. 
The index function is constant 
on any segment of the form $(P,Q]$, where $P$ and $Q$ 
are two consecutive marked points. We consider $\Theta_{L_0} (C)$ as a poset with respect to $\preceq_{L_0}$. If  $C_j$ is a branch of $C$ and $C_j \ne L_0$, then the restriction of the exponent function to the segment $[L_0, C_j] \subset \Theta_{L_0} (C)$ is a homeomorphism of posets onto $[0, \infty]$. A point $P \in \Theta_{L_0} (C)$ is \emph{rational} if $\ex_{L_0} (P) \in \Q_{> 0}$.
 
\medskip 
 
The functions $\ex_{L_0}$ and $\de_{L_0}$ 
determine the \emph{contact complexity function},
\[
{\ic}_{L_0}\colon\Theta_{L_0} (C)\to\overline\R_{\geq 0},
\]  defined by 
$\ic_{L_0} (P) = \int_{L_0}^P \frac{d \ex_{L_0}}{\de_{L_0}}$. 
The knowledge of the functions $\de_{L_0}$ and $\ic_{L_0}$ allow us to recover 
the intersection multiplicity  
of pairs of different branches of $C \cup L_0$. 
One has that $\de_{L_0} (L_0) = 1 $ while 
$\de_{L_0} (C_j)  = (L_0 \cdot C_j)$ for  $j \in \{1, \dots, r \}$ (see \cite[Rem. 3.25]{GGP19}). 
In addition, by  \cite[Cor. 3.26]{GGP19}, 
we have 
\begin{equation} \label{eq:int3}
(C_i \cdot C_j) = \de_{L_0} (C_i )  \de_{L_0} (C_j) \ic_{L_0} ( \langle L_0, C_i, C_j \rangle ) 
\mbox{ for }  i, j  \in \{1, \dots, r \}. 
\end{equation}
For any branch $A$ different from $L_0$ the \emph{attaching point}  of $A$ to the tree $\Theta_{L_0} (C)$
is 
\[P_A  := \max \{ \langle L_0, A , C_j \rangle \mid  j =1, \dots, r \} \in \Theta_{L_0} (C),\] where the tripods are viewed in  $\Theta_{L_0} (C \cup A)$ and the maximum is taken with respect to ${\prec_{L_0}}$.

\begin{example} \label{rem:chexp}
Assume that $L_0, L_1, \dots, L_g$ is a sequence of maximal contact curves of $\psi_{A}$, for a branch $A$ of $C$. 
Then, the discontinuity locus of the restriction of the index function  $\de_{L_0}$ to the segment $[L_0, A]$
is equal to the set of points  $\langle L_0, L_j, A \rangle$, for $j =1, \dots,g$, where the tripods 
are viewed  
in the tree $\Theta_{L_0} (C \cup L_1 \cup \dots \cup L_g )$.
With Notation \ref{not:GA}, one has that  $\de_{L_0} (\langle L_0, L_1, A \rangle) =1 = \de_{L_0} (L_1)$ and 
$\de_{L_0} (\langle L_0, L_j, A \rangle) = n_1 \dots n_{j-1} = \de_{L_0} (L_j)$ for $j=2, \dots, g$.
 By \eqref{eq:int3} we get that $\ic_{L_0} (\langle L_0, L_1, A \rangle) = \g_0^{-1} \g_1$ and
$\ic_{L_0} (\langle L_0, L_j, A \rangle) = \g_0^{-1} (n_1 \dots n_{j-1})^{-1} \g_j$ for $j=2, \dots, g$.
In addition, we get
$\ex_{L_0} (\langle L_0, L_j, A \rangle)  = \beta_j/\beta_0 $, for $j=1, \dots, g$,
where $\beta_0 = \g_0$, 
$\beta_1 = \g_1$ and $\beta_j - \beta_{j-1} = \g_j - n_{j-1} \g_{j-1} $
for $j=2, \dots, g$. 
\end{example}

If $\k$ is of characteristic zero,
then the sequence $(\beta_0, \beta_1, \dots, \beta_g)$ is the \emph{characteristic}
of the branch $A$, which can be defined in terms of the Newton Puiseux series of $A$ with respect to $L_0$ (see Definition 3.2 of \cite{Z}). 
 If $\k$ is a field of positive characteristic
 then the notion of characteristic exponents  of a branch was introduced in Campillo's book \cite[Chapter 3]{Campillo}.

\begin{figure}[h!] 
\begin{center}
\begin{tikzpicture}[scale=0.45]

     \node[right] at (-1,-1) {$\Theta_{L_0} (D)$};

  \draw [-, color=black, thick](0,0) -- (0, 8) ; 
   \node[draw,circle, inner sep=1.5pt,color=black, fill=black] at (0,0){};
   \node [right] at (0,0) {$L_0$};
    \node[draw,circle, inner sep=1.5pt,color=black, fill=black] at (0,8){};
   \node [right] at (0,8) {${L_{5}}$};
 
    \node[draw,circle, inner sep=1.5pt,color=black, fill=black] at (0,2){};
  %  \node [left] at (0,2) {$\mathbf {\frac{3}{2}}$};
   \node [left] at (0,2) {$P_1$};
   
     \node[draw,circle, inner sep=1.5pt,color=black, fill=black] at (0,5){};
   % \node [left] at (0,5) {$\mathbf {\frac{7}{4}}$};
     \node [left] at (0,5) {$P_4$};

     \node [above] at (1,5.5) {$2$};
   
    \draw [-, color=black, thick](0,5) -- (2,6) ; 
     \node[draw,circle, inner sep=1.5pt,color=black, fill=black] at (2,6){};
  \node [right] at (2,6) {${L_{4}}$};
   
 \draw [-, color=black, thick](0,2) -- (6,5) ;   
     \node[draw,circle, inner sep=1.5pt,color=black, fill=black] at (6,5){};
  \node [right] at (6,5) {${L_{3}}$};
  
%    \node [above] at (4,6) {$\infty$}; 
 \node[draw,circle, inner sep=1.5pt,color=black, fill=black] at (4,4){};  
 %  \node [below] at (4,4) {$\mathbf {\frac{5}{2}}$};
    \node [below] at (4,4) {$P_3$};

     \draw [-, color=black, thick](4,4) -- (6,2) ; 
       \node[draw,circle, inner sep=1.5pt,color=black, fill=black] at (6,2){};  
          \node [below] at (5,3) {$3$};
    \node [right] at (6,2) {${L_{2}}$};
   
   \node[draw,circle, inner sep=1.5pt,color=black, fill=black] at (2,3){};  
     %\node [below] at (2,3) {$\mathbf {\frac{5}{3}}$};
      \node [below] at (2,3) {$P_2$};
  \draw [-, color=black, thick](2,3) -- (4,1) ; 
       \node[draw,circle, inner sep=1.5pt,color=black, fill=black] at (4,1){};  
          \node [below] at (3,2) {$1$};
          
       \node [right] at (4,1) {${L_1}$};
     \node [below] at (1,3.5) {$1$};
     \node [above] at (3,3.5) {$3$};
       \node [below] at (5,5.5) {$6$};
     \node[right] at (0,1) {$1$};
       \node[right] at (0,4) {$2$};
      \node[right] at (0,7) {$4$};
\end{tikzpicture}
\end{center}
\caption{The Eggers-Wall tree of Example      
         \ref{ex:EW}} 
\label{fig:EW}
\end{figure}
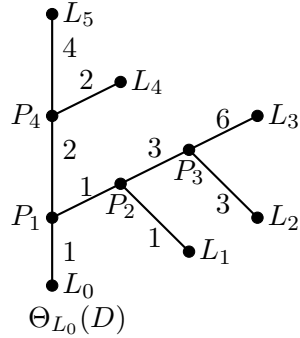

\begin{example} \label{ex:EW} In this example we assume for simplicity 
that the field $\k$ is of characteristic zero.
Let us consider the branches $L_0 = Z(X)$, $L_1 = Z(Y)$, and
$L_j$, for $j=2,\dots, 5$, parametrized respectively by the 
Newton-Puisseux series
$\zeta_2 := X^{5/3}$, $\zeta_3 := X^{5/3} + X^2 + X^{5/2}$, 
$\zeta_4 := X^{3/2}$, and $\zeta_5 :=  X^{3/2} + X^{7/4}$.
We have represented the Eggers-Wall tree $\Theta_{L_0} (D)$ of 
the curve $D:=\cup_{j=1}^5 L_j$ in Figure \ref{fig:EW}. 
The marked points are $L_0, \dots, L_5, P_1, \ldots, P_4$. 
In Figure \ref{fig:EW} is indicated
the constant value of the index function on every interval $(Q, Q']$ 
between consecutive marked points $Q \preceq_{L_0} Q'$. 
We have that $\ex_{L_0} (P_1) = 3/2$, $\ex_{L_0} (P_2) = 5/3$, 
$\ex_{L_0} (P_3) = 5/2$,  $\ex_{L_0} (P_4)  = 7/4$, 
while  $\ex_{L_0} (L_0) =0$ and  $\ex_{L_0} (L_j) = \infty$ 
for $j=1,\dots, 5$. 
By definition, the value of the contact complexity function $\ic_{L_0}$ 
at a point $P$ of the tree 
is just a finite sum. For instance, one has
$
\ic_{L_0} (P_3) = 
\int_{L_0}^{P_2}  \frac{1}{\de_{L_0}} d \ex_{L_0} +  \int_{P_2}^{P_3}  \frac{1}{\de_{L_0}} d \ex_{L_0} =   \frac{5}{3}
+  \frac{1}{3} \left(\frac{5}{2} - \frac{5}{3}\right) = \frac{35}{18}
% \frac{1}{1} (\ex_{L_0} (P_2) -  \ex_{L_0} (L_0) ) + \frac{1}{3} (\ex_{L_0} (P_3) -  \ex_{L_0} (P_2) )
$, $\ic_{L_0} (P_4) =
\int_{L_0}^{P_1}  \frac{1}{\de_{L_0}} d \ex_{L_0} +  \int_{P_1}^{P_4}  \frac{1}{\de_{L_0}} d \ex_{L_0}  =
 \frac{3}{2}  +  \frac{1}{2} \left(\frac{7}{4} - \frac{3}{2}\right) =
\frac{13}{8} $, and similarly
$\ic_{L_0} (P_1) = \frac{3}{2}$,  $\ic_{L_0} (P_2) = \frac{5}{3} $,   $\ic_{L_0} (L_0)=0$, 
and $\ic_{L_0} (L_j)=\infty$ for $j=1, \dots, 5$.
\end{example}

\medskip

If $P$ is a rational point of an Eggers-Wall tree,
we denote by $\de_{L_0}^+ (P)$ 
the lowest common multiple of
$\de_{L_0} (P)$ and 
of the denominator of $\ex_{L_0} (P)$ written as an irreducible fraction  (see \cite[Definition 3.14 and Proposition 3.16]{GGP19}).
 Let $A$ and $B$ be two branches such that 
 $P = \langle L_0, A, B \rangle$ in $\Theta_{L_0} ( A \cup B)$
 and the restriction of the index function to the segment 
 $[L_0, B ]$ is continuous at $P$. Then, the indices
 $\de_{L_0} (P)$ and $\de_{L_0}^+ (P)$ can be seen as 
 limits when $P'$ tends to $P$: 
 \[
 \lim_{P' \in (L_0, B), \, P' \to P} \de_{L_0} (P') = \de_{L_0} (P)
 \mbox{ and }
 \lim_{P' \in (P, A),  \, P' \to P} \de_{L_0} (P') = \de_{L_0}^+ (P)
 \]
 (see \cite[Lemma 3.15]{GGP19}).
 In particular, if $\de_{L_0} (P) < \de_{L_0}^+ (P)$
 then $P$ is a point of discontinuity of the index function in restriction to $[L_0, A]$.
 This implies that 
 $\ex_{L_0} (P)$ is a characteristic exponent of $A$ with respect to $L_0$.

\begin{definition} \label{def:rootlevel}
Let $P \in \Theta_{L_0} (C)$ and set
$Q_P$ the minimum, with respect to the order $\prec_{L_0}$, of the closure of the connected component of the  set $\{ P' \in  \Theta_{L_0} (C) \mid 
\de_{L_0} (P') =  \de_{L_0} (P) \}$
which contains the point $P$. 
If $P$ is a rational point we set
\begin{equation} \label{eq: Qrenor}
 \ex_{Q_P} (P): = \de_{L_0} (P) \left( \ex_{L_0} (P) - \ex_{L_0} (Q_P) \right) = \frac{m_P}{n_P},
\end{equation}
where $(n_P, m_P)$ is a couple of coprime integers. 
\end{definition}

Notice that $Q_P= L_0$ if $\de_{L_0} (P) = 1$ and then  $\ex_{Q_P} (P) = \ex_{L_0} (P)$.
Otherwise,
$Q_P$ is a point of discontinuity of the index function, hence it is a rational point,  and 
$\de_{L_0} (Q_P) < \de_{L_0}^+ (Q_P) = \de_{L_0} (P)$. 
\begin{remark} \label{rem: QP2}
Let $P_1, P_2 \in \Theta_{L_0} (C) $ with $P_1 \preceq_{L_0} P_2$. 
If the index function is constant on $(P_1, P_2]$ then
$Q_{P_2} \preceq_{L_0} P_1$
with equality if and only if $P_1 = L_0$ or if $\de_{L_0} (P_1) < \de_{L_0}(P_2)$.
\end{remark}

If $P$ is a rational point it follows from the definitions that 
\begin{equation} \label{eq: indexplus}
\de_{L_0}^+ (P) = \de_{L_0} (P) \cdot n_P.
\end{equation}

\subsubsection{The finite local tropicalization} 

In this section we prove that the finite local tropicalization 
$\trop(\mathcal{U})$  of $I_S$ is the support of a fan and 
we describe it in terms of  the Eggers-Wall tree.
We use the  
 functions $\ex_{L_0}$, $\de_{L_0}$ and 
$\ic_{L_0}$, which can be expressed in terms
of natural "coordinate functions" on the valuative tree $\cV_{L_0}$:
the \textit{log-discrepancy},  the \textit{multiplicity} and 
 the \textit{self-interaction} (we use here the terminology of  \cite{GGP19}).
 These functions are crucial in the work of 
 Favre and Jonsson \cite{FJ}, see also Jonsson's survey  \cite[section 7]{Jon15}. 
 
\medskip 

An embedding of an Eggers-Wall tree of a plane curve 
in the semivaluation space 
is described in \cite[Section 8]{GGP19} if $k$ is a field of characteristic zero, and
in \cite[Remark 5.37]{GGP19b} if $\k$ has positive characteristic. In the following proposition 
we describe this embedding for the curve $D = \cup_{j=0}^m L_j$
(recall the notation at the beginning of Section \ref{sec:TEpairs}).

\begin{proposition} \cite[Prop. 8.10 and Th. 8.11]{GGP19}
There is an embedding of rooted trees: 
\begin{equation} \label{eq:emb}
V_{L_0} : \Theta_{L_0} (D) \hookrightarrow \cV_{L_0},  \quad P \mapsto V_{L_0}^P
\end{equation}
such that for any $P \in  \Theta_{L_0} (D) $ and any plane branch $A$ 
on $S$ we have 
\begin{equation} \label{VPA}
V_{L_0}^P(A) =   \left\{ 
                \begin{array}{ccl}
                          (L_0 \cdot A) \cdot {\ic}_{L_0} (\langle L_0, P, A \rangle)  &  
                          \mbox{ if } & A \neq L_0, \\
                       1 &  \mbox{ if } & A = L_0,  
                  \end{array}  \right. 
\end{equation}
where $\langle L_0, P, A \rangle$ is the center of the tripod
defined by $L_0$, $P$, and $A$ on the Eggers-Wall tree $\Theta_{L_0} (D \cup A)$. In addition, $V_{L_0}^{L_j} = I_{L_0}^{L_j}$ 
for $1\leq j\leq m$ and $V_{L_0}^{L_0}=\ord^{L_0}$.
\end{proposition}

Notice that if $P$ is an interior 
point of $\Theta_{L_0}(D)$ then $V_{L_0}^P(A) \in \R_{>0}$ for any branch $A$, since the function ${\ic}_{L_0}$ restricted to the interior of the tree
$\Theta_{L_0}(D)$ has values in $\R_{>0}$.

\begin{definition} \label{def:WP}
If $P$ is an interior point of the tree $\Theta_{L_0} (D)$ , or if $P = L_0$, then we set:
\begin{equation} \label{defwp}
w^P := \trop (V_{L_0}^P) = (V_{L_0}^P (L_0), \dots,  V_{L_0}^P (L_m)) \in \R^{m+1}_{>0}.
\end{equation}
For $0 \leq j \leq m$, we define 
\begin{equation} \label{defwp-2}
w^P:= e_j,
\end{equation}
where $e_0, \dots, e_{m}$
 is the canonical basis of $\R^{m+1}$. 
 \end{definition}

\begin{remark} \label{r-inf} $\,$
Notice that $w^{L_0}=\trop(V_{L_0}^{L_0}) =\trop(\ord^{L_0})$. 
However, for $1\leq j\leq m$, 
we have $w^{L_j}=\trop(\ord^{L_j})$, while 
$\trop(V_{L_0}^{L_j}) = \trop (I_{L_0}^{L_j}) \notin\R^{m+1}_{\geq 0}$. 
\end{remark}

\begin{lemma} \label{block1}
Let $P_1, P_2$ be two consecutive marked points of $\Theta_{L_0} (D)$. 
Then, the cone spanned by $\{ w^P \mid P \in [P_1, P_2] \}$ is equal to 
$\R_{\geq0} w^{P_1} + \R_{\geq0} w^{P_2}$. 
\end{lemma}
\begin{proof}
We may assume that $P_1 \preceq_{L_0} P_2$. 
By definition of the contact complexity function 
we get the following equality for $P \in [P_1, P_2]$: 
\begin{equation} \label{eq:icP}
 {\ic}_{L_0} (P) = \int_{L_0}^{P_1} \frac{d \ex_{L_0}}{\de_{L_0}} + 
 \int_{P_1}^P \frac{d \ex_{L_0}}{\de_{L_0}} = 
 {\ic}_{L_0} (P_1) + \frac{1}{\de_{L_0} (P_2)} \left(\ex_{L_0} (P) - 
 \ex_{L_0} (P_1)  \right), 
\end{equation}
where the second equality holds since $\de_{L_0}$ has constant value 
${\de_{L_0} (P_2)}$ on the segment $(P_1, P_2]$.

We discuss first the value $V_{L_0}^P(L_j)$ according to the relative position of $L_j$ with 
respect to $P_1$ and $P_2$, for $j\in\set{0,\ldots,m}$. 
By hypothesis, no ramification 
point of the tree $\Theta_{L_0}(D)$ belongs to the segment $(P_1, P_2)$. 
Then, we have two cases for $L_j$: 
\begin{enumerate}
\item \label{first} $\langle L_0, L_j , P \rangle = P$ if $P_2 \preceq_{L_0} L_j$;
\item \label{second} $\langle L_0, L_j , P \rangle = \langle L_0, L_j , P_1 \rangle $,  otherwise. 
\end{enumerate}
In case \eqref{first}, it follows from the definition of $V_{L_0}$ and 
the equality \eqref{eq:icP} that 
\begin{equation} \label{aux2}
V^{P}_{L_0} (L_j) =
(L_0 \cdot L_j)  \left( {\ic}_{L_0} (P_1) + \frac{1}{\de_{L_0} ( P_2)} 
 \left(\ex_{L_0} (P) - \ex_{L_0} (P_1)  \right) \right).
\end{equation} 
In case \eqref{second}, from the definition of $V_{L_0}$ we see that
\begin{equation} \label{aux1}
V^P_{L_0} (L_j) =
V^{P_1}_{L_0} (L_j) .
\end{equation}

We prove the lemma now by distinguishing two cases:

$\bullet$ The point $P_2$ is not a leaf of $\Theta_{L_0} (D)$. 
If $P$ runs through $[P_1, P_2]$
then the number 
\[
s(P)  := \frac{\ex_{L_0} (P) - \ex_{L_0} (P_1)}{\ex_{L_0} (P_2) - \ex_{L_0} (P_1)},
\]
runs through the interval $[0,1]$. 
From  \eqref{defwp}, \eqref{aux1} 
and \eqref{aux2} we 
check coordinate-wise that the following equality holds 
$$
w^P = s(P) w^{P_2} + (1 - s(P)) w^{P_1},
$$
which implies the assertion in this case.

$\bullet$ The point $P_2$ is a leaf of $\Theta_{L_0} (D)$. Then 
$P_2 = L_i$ for some $1 \leq i \leq m$. 
If $j \in \{0, \dots, m \}$ and $j \ne i$
then $\langle L_0, L_j , P \rangle = \langle L_0, L_j , P_1 \rangle $ for $P \in [P_1, P_2]$, i.e., we are in the case \eqref{second} discussed above, while we are in case \eqref{first} for $L_i$.
If $P$ runs through $[P_1, P_2)$ then the number
\[
\tilde{s}(P)  := \frac{(L_0 \cdot L_i) }{\de_{L_0} (P_1)} 
 \left(\ex_{L_0} (P) - \ex_{L_0} (P_1) \right),
\]
runs through the interval $[0, \infty)$. 
By \eqref{defwp}, \eqref{defwp-2},  \eqref{aux1} and \eqref{aux2} we get that 
$$
w^P = w^{P_1} + \tilde{s}(P) w^{P_2}.
$$
This ends the proof in this case. 
\end{proof}

\begin{lemma} \label{block2}
The map $\Theta_{L_0} (D) \to \R^{m+1}_{\geq 0}$, 
$P \mapsto w^P$ is injective and its restriction to the set of non-leaf 
points of the tree $\Theta_{L_0} (D)$ is continuous.
\end{lemma}
\begin{proof}
Take $j\in\set{1,\ldots,m}$. The restriction of $\ic_{L_0}$ to 
$[L_0, L_j]$ is an increasing homeomorphism onto $[0,\infty]$. For $P\in[L_0,L_j)$, 
we have $w_j^P=V_{L_0}^P(L_j)=(L_0 \cdot L_j)\ic_{L_0}(P)$. 
Notice that $w_0^P = 1$ for every non-leaf point $P$.
Hence 
the assertion on the continuity of this map follows. 
In addition, in order to prove the injectivity it is enough to show that 
$w^P \ne w^Q$ for any pair of non comparable points $P, Q \in \Theta_{L_0} (D) $ with respect to the order
$\preceq_{L_0}$. 
The special case where one of these points is an end 
follows from Definition \ref{def:WP}, using 
that $w^U = \trop(V_{L_0}^U) \in \R^{m+1}_{>0}$ for any interior point $U$. 
Suppose now that $P, Q$ are interior points. Let $1 \leq  j \leq m$ be such that 
$P \in [L_0, L_j]$. Then, if 
$F := \langle L_0, Q, P \rangle$ we have that $F \prec_{L_0} P$,  $F \prec_{L_0} Q$ and $F= \langle L_0, Q, L_j \rangle$. Therefore, 
$w_j^Q = (L_0 \cdot L_j)\ic_{L_0}(F) = w_j^F < w_j^P$, 
since the map $\ic_{L_0}$ is strictly increasing on $[L_0, L_j]$.
\end{proof}

\begin{notation} \label{not:aug}
The augmented set of marked points
$\mathcal{A}_{L_0}^C (D)$ 
consists of those points $Q \in \Theta_{L_0} (D)$ which are marked points when seen in 
 $\Theta_{L_0} (D \cup C)$. 
See Remark \ref{rem:aug} below for a geometrical interpretation of the set $\mathcal{A}_{L_0}^C (D)$.
\end{notation}

As a consequence of Lemmas \ref{block1} and \ref{block2} we get that:

\begin{proposition} \label{prop:complexT}
There exists a fan $\cT$ with respect to the lattice  $\Z^{m+1}$ consisting of cones of dimension at most $2$, 
whose support is equal to 
\[
|\cT| := \R_{\geq 0} \{ w^P \mid P \in \Theta_{L_0} (D) \}.
\]
Its $2$-dimensional cones are
$\R_{\geq0} w^{P_1} + \R_{\geq0} w^{P_2}$, for $P_1, P_2 \in \mathcal{A}_{L_0}^C (D)$  consecutive points on 
the tree $\Theta_{L_0} (D)$.
\end{proposition}

Now, we prove that: 
\begin{theorem} \label{loc-trop} 
With the previous notation, we have 
$ \trop (\mathcal{U}) = |\cT|.$  
\end{theorem} 

\begin{proof} Since $w^{L_j}=e_j$ for $j=0,\ldots,m$, 
it is enough to show that 
\begin{equation} \label{eq:tropi}
\trop\left(\cV_{L_0}\setminus
\set{\ord^{L_0},I_{L_0}^{L_1},\ldots,I_{L_0}^{L_m}}\right)=
\set{ w^P \mid P  
\text{ is an interior point of } \Theta_{L_0} (D) }
\end{equation}
(see Lemma \ref{lem:tropU}). 
For any interior point $P \in \Theta_{L_0} (D)$, we have 
$w^P=\trop(V_{L_0}^P)$, with 
$ V_{L_0}^P\in 
\cV_{L_0} \setminus \set{\ord^{L_0}, I_{L_0}^{L_1}, \dots, I_{L_0}^{L_m}}$. Therefore, 
the inclusion $\supset$ in the equality \eqref{eq:tropi} follows. 

Next we show the other inclusion.
Take $\nu \in \cV_{L_0} \setminus \{ \ord^{L_0} , I_{L_0}^{L_1}, \dots, I_{L_0}^{L_m}\}$. 
Since the valuative tree  $\cV_{L_0}$ is an $\R$-tree, and 
the set $V_{L_0} (\Theta_{L_0} (D))$ is a closed subtree of it,  
there exists a unique \emph{attaching point} 
$\nu'$ of $\nu$ to $V_{L_0} (\Theta_{L_0} (D))$ 
(see \cite[sections 2, 6 and 8]{GGP19}). 
The attaching point $\nu'$ is the unique point of $\cV_{L_0}$
such that the intersection of the interval $[\nu,\nu']$ with the subtree 
is reduced to $\set{\nu'}$. 
Notice that the semivaluations
$\ord^{L_0} , I_{L_0}^{L_1}, \dots, I_{L_0}^{L_m}$ are the ends of $V_{L_0} (\Theta_{L_0} (D))$, and these 
semivaluations are also ends of $\cV_{L_0}$. 
We get that $\nu'$ must be of the form 
$\nu' = V_{L_0}^Q$ for 
some interior point $Q \in \Theta_{L_0} (D)$. Since 
$\nu$ belongs to $\cV_{L_0}$ one has that $\nu(L_0) = 1 = w^Q_0$.
Let us prove that $\trop(\nu)=w^Q$, 
that is, $\nu (L_j)  = w^Q_j$
for $1 \leq j \leq m$. 

We rely below on some results of \cite{FJ} and also of \cite{GR}, which we use following the presentation of  \cite{GGP19}.
There exists a function 
$\langle \cdot, \cdot \rangle : \cV^* \times \cV^* \to [0, \infty]$, called the \emph{bracket}, such that 
for any branch $A$ on $S$ and any $\mu \in \cV^*$ we have 
\begin{equation} \label{eq:muA}
\langle I^A, \mu \rangle = \mu (A) 
\end{equation}
(see \cite[Prop. 7.5 and 7.11]{GGP19}).
In particular, if $B$ is another plane branch
we get that 
\begin{equation}\label{eq:AB}
\langle I^A, I^B \rangle = I^B (A) = (A\cdot B).
\end{equation}
The   \emph{relative self-interaction} function $\si_{L_0} : \cV_{L_0} \to [0,\infty]$ is defined by $\si_{L_0} (\mu) = \langle \mu, \mu \rangle$.

Let us fix $j \in \{ 1, \dots, m \}$. 
By the \emph{generalized tripod formula} applied to the semivaluations $I^{L_j}$ and $\nu$ with respect to the observer $L_0$, 
the following equality holds
(see \cite[Prop. 7.18]{GGP19}):
\begin{equation} \label{eq2}
\langle I^{L_j}, \nu \rangle= \langle I^{L_0}, \nu \rangle   \langle I^{L_0},  I^{L_j} 
\rangle \, \si_{L_0} ( \langle \ord_{L_0}, I_{L_0}^{L_j}, \nu  \rangle  ).
\end{equation}
By definition of the attaching point, we have that 
\begin{equation} \label{block3} 
\langle \ord_{L_0}, I_{L_0}^{L_j}, \nu \rangle = \langle \ord_{L_0}, I_{L_0}^{L_j}, V_{L_0}^Q  \rangle.
\end{equation}
Taking into account \eqref{eq:muA}, \eqref{eq:AB} and \eqref{block3}, we can 
reformulate \eqref{eq2} as follows:
\begin{equation} \label{eq:inter}
\nu(L_j)=\nu(L_0) (L_0 \cdot L_j) \si_{L_0}  (\langle \ord_{L_0}, I_{L_0}^{L_j}, V_{L_0}^Q  \rangle ) = (L_0 \cdot L_j) \si_{L_0} ( \langle \ord_{L_0}, I_{L_0}^{L_j}, V_{L_0}^Q  \rangle ).
\end{equation}
On the one hand, we have that $V_{L_0}^{L_0} = \ord^{L_0}$ and $V_{L_0}^{L_j} = I_{L_0}^{L_j}$. On the other hand
by \cite[Th. 8.18]{GGP19} one has
$\si_{L_0}  \circ V_{L_0} = \ic_{L_0}$.
Combining these facts we reformulate the equality \eqref{eq:inter} as
$
\nu (L_j) = (L_0 \cdot L_j) \,  \ic_{L_0} (\langle L_0, L_j, Q \rangle) = w^Q_j$.
\end{proof}

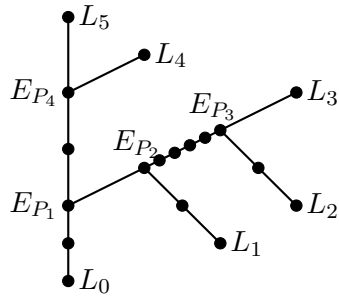
\begin{figure}[h!] 
\begin{center}
\begin{tikzpicture}[scale=0.5]

   %   \node[right] at (-1,-1) {$\Theta_L (C)$};

  \draw [-, color=black, thick](0,0) -- (0, 7) ; 
   \node[draw,circle, inner sep=1.5pt,color=black, fill=black] at (0,0){};
   \node [right] at (0,0) {$L_0$};
 \node[draw,circle, inner sep=1.5pt,color=black, fill=black] at (0,7){};
   \node [right] at (0,7) {${L_{5}}$};
 
    \node[draw,circle, inner sep=1.5pt,color=black, fill=black] at (0,2){};
    \node[draw,circle, inner sep=1.5pt,color=black, fill=black] at (0,3.5){};
     \node[draw,circle, inner sep=1.5pt,color=black, fill=black] at (3,2){};
      \node[draw,circle, inner sep=1.5pt,color=black, fill=black] at (5,3){};
  \node[draw,circle, inner sep=1.5pt,color=black, fill=black] at (2.4,3.2){};
    \node[draw,circle, inner sep=1.5pt,color=black, fill=black] at (2.8,3.4){};
     \node[draw,circle, inner sep=1.5pt,color=black, fill=black] at (3.2,3.6){};
  \node[draw,circle, inner sep=1.5pt,color=black, fill=black] at (3.6,3.8){};
  %  \node [left] at (0,2) {$\mathbf {\frac{3}{2}}$};
   \node [left] at (0,2) {$E_{P_1}$};
   
     \node[draw,circle, inner sep=1.5pt,color=black, fill=black] at (0,5){};
     \node[draw,circle, inner sep=1.5pt,color=black, fill=black] at (0,1){};
   % \node [left] at (0,5) {$\mathbf {\frac{7}{4}}$};
     \node [left] at (0,5) {$E_{P_4}$};

  %   \node [above] at (1,5.5) {$2$};
   
    \draw [-, color=black, thick](0,5) -- (2,6) ; 
     \node[draw,circle, inner sep=1.5pt,color=black, fill=black] at (2,6){};
  \node [right] at (2,6) {${L_{4}}$};
   
 \draw [-, color=black, thick](0,2) -- (6,5) ;   
     \node[draw,circle, inner sep=1.5pt,color=black, fill=black] at (6,5){};
  \node [right] at (6,5) {${L_{3}}$};
  
%    \node [above] at (4,6) {$\infty$}; 
 \node[draw,circle, inner sep=1.5pt,color=black, fill=black] at (4,4){};  
 %  \node [below] at (4,4) {$\mathbf {\frac{5}{2}}$};
    \node [above] at (3.8,4) {$E_{P_3}$};

     \draw [-, color=black, thick](4,4) -- (6,2) ; 
       \node[draw,circle, inner sep=1.5pt,color=black, fill=black] at (6,2){};  
      %    \node [below] at (5,3) {$3$};
    \node [right] at (6,2) {${L_{2}}$};
   
   \node[draw,circle, inner sep=1.5pt,color=black, fill=black] at (2,3){};  
     %\node [below] at (2,3) {$\mathbf {\frac{5}{3}}$};
      \node [above] at (1.8,3) {$E_{P_2}$};
  \draw [-, color=black, thick](2,3) -- (4,1) ; 
       \node[draw,circle, inner sep=1.5pt,color=black, fill=black] at (4,1){};  
  %         \node [below] at (3,2) {$1$};
          
       \node [right] at (4,1) {${L_1}$};
    %  \node [below] at (1,3.5) {$1$};
   %   \node [above] at (3,3.5) {$3$};
   %     \node [below] at (5,5.5) {$6$};
  %    \node[right] at (0,1) {$1$};
  %      \node[right] at (0,4) {$2$};
  %     \node[right] at (0,7) {$4$};
\end{tikzpicture}
\end{center}
\caption{The dual graph of the total transform of the curve in Example \ref{ex:EW} 
under its minimal resolution.} 
\label{fig:totaltrans}
\end{figure}

\begin{example}
Let us come back to the case of Example \ref{ex:EW}. 
We have represented in Figure \ref{fig:totaltrans}
the dual graph of the total transform of $D\cup L_0$ 
by its minimal embedded resolution. 
Take $x_j\in m_R$ such that $L_j=Z(x_j)$, for $j=0, \dots, 5$.
By Corollary \ref{lem:genseqratpt}, the sequence $x_0,\dots, x_5$ is 
a generating sequence 
of $(\nu_{E_{P_1}},\ldots,\nu_{E_{P_4}})$. 
Then, we use formula \eqref{VPA} to find the vector 
$w^{P_i} = (w^{P_i}_0, \dots, w^{P_i}_5)$, for $i=1, \dots,4$. 
We obtain $w^{P_i}_0 =1$ and 
$w^{P_i}_j = (L_0 \cdot L_j) \ic_{L_0} (\langle L_0, L_j, P_i \rangle)$, 
for $j=1, \dots, 5$.
In particular, for $i=4$, we get that 
$\langle L_0, L_j, P_4 \rangle = P_1$ for $j=1, 2,3$ while 
$\langle L_0, L_j, P_4 \rangle = P_4$ for $j=4,5$. 
Then, one obtains $w^{P_4} = (1, {3}/{2}, {9}/{2}, 9, {13}/{4}, {13}/{2})$, 
and similarly
$w^{P_1} = (1, {3}/{2}, {9}/{2}, 9, 3, 6)$, 
$w^{P_2} = (1, {5}/{3}, 5, 10, 3, 6)$,
and 
$w^{P_3} = (1,{5}/{3}, {25}/{6},  {25}/{3}, 3, 6)$.
\end{example}

We will make use of the following lemma in section \ref{sec:lastone}.

\begin{lemma}\label{lem:charrat}
Let $P$ be an interior point of $\Theta_{L_0}(D)$. The point $P$ is 
rational if and only if $w^P\in\Q_{>0}^{m+1}$.
\end{lemma}

\begin{proof}
 The point $Q$ is a rational point if and only if $\ex_{L_0}(Q)\in\Q_{>0}$, 
which is equivalent to $\ic_{L_0}(Q)\in\Q_{>0}$ since $\ic_{L_0}$ is defined by the integral of a lower semicontinuous function, which 
is not continuous at a finite number of rational points, and which has a finite number of rational values.
Taking into account this 
observation and \eqref{VPA}, we see that $w^P$ is rational if 
$P$ is rational. Let us now show the converse. 
Assume that $w^P\in\Q_{>0}^{m+1}$. Take $P_1, P_2$ two consecutive 
marked points of $\Theta_{L_0} (D)$ such that $P\in[P_1,P_2]$ and $j\in\set{1,\ldots,m}$ such that $P_2\preceq_{L_0}L_j$. By 
\eqref{aux2}, we have that $w^P_j\in\Q_{>0}$ if and only if 
$\ex_{L_0}(P)\in\Q_{>0}$.
\end{proof}

We need the following particular property of the Eggers-Wall tree 
$\Theta_{L_0} (D)$, which does not hold for Eggers-Wall trees of 
arbitrary plane curves singularities.

\begin{lemma} \label{lem:levels}
The leafs of the closures of the levels of the index function $\de_{L_0}$ on the tree $\Theta_{L_0} (D)$ are 
leafs of $\Theta_{L_0} (D)$. 
\end{lemma}
\begin{proof}
The closure of a connected component $\Xi$ of the index function on the tree $\Theta_{L_0} (D)$ is a tree
rooted at the point of minimal exponent of $\Xi$. If $P$ is an end of $\Xi$ which is not an end of $\Theta_{L_0} (D)$ then 
$P$, seen on $\Theta_{L_0} (D)$, is a point of discontinuity of $\de_{L_0}$. 
Take a branch $A$ of $D$ such that $P \prec_{L_0} A$. 
By Lemma \ref{lem:mgs}, there exists $1 \leq g \leq m$ such 
that, up to relabelling $L_1, \dots, L_m$, we have that 
$L_0, \dots, L_{g}$ is a sequence of maximal contact curves of $\psi_A$.  
Then, $P$ is also a point of 
discontinuity of the restriction of $\de_{L_0}$ to the segment $[L_0, A]$. 
By Example \ref{rem:chexp} applied to the branch $A$ of $D$ 
and the sequence of maximal contact curves 
$L_0, \dots, L_{g}$ of $\psi_A$, there exists 
$j \in \{ 1, \dots,  g \}$ such that 
$P = \langle L_0, L_j, A \rangle$  and then
$\de_{L_0} (P) = \de_{L_0} (L_j)$. 
We get that  $P \prec_{L_0} L_j$ with $L_j \in \Xi$, which is a contradiction.  
\end{proof}
%---------------------------------------------------------

\subsection{Eggers-Wall trees and the minimal embedded resolution}

In this section we describe the divisors of the minimal embedded resolution of $C$, 
and more generally of the 
\emph{representing divisors} of rational points in the Eggers-Wall tree $\Theta_{L_0} (D)$.

\medskip 

Let $P$ be a rational point of the Eggers-Wall tree $\Theta_{L_0} (D)$. 
As it is explained in \cite[Def. 8.14]{GGP19}, 
there exists a unique exceptional prime divisor $E_P$ 
(up to birational transformation) such that 
its associated divisorial valuation $\nu_{E_P}$ is proportional to the valuation $V_{L_0}^P$, that is, 
\begin{equation}\label{eq: EP}
 \nu_{E_P}  = \nu_{E_P} (L_0) \cdot  V_{L_0}^P . 
\end{equation}
We say that $E_P$ is the \emph{representing divisor} of $P$
(see \cite[Prop. 8.16]{GGP19} for some of its properties). 

\medskip

\begin{lemma} \label{lem:curvetta}
Let $A$ be a branch, $L_0$ a smooth branch, and  $Q \in \Theta_{L_0} (A)$  a rational point.
Assume that 
$ 
\lim_{P \in (Q, A), \, P \to Q} \de_{L_0} (P) = \de_{L_0}^+ (Q)$. 
If $\de_{L_0}^+ (Q) = \de_{L_0} (A)$,
then $A$ is a curvetta of $E_P$ at the minimal resolution of $\nu_{E_Q}$. 
\end{lemma}
\begin{proof} 
Denote by  $\psi \colon (S(\psi ),E(\psi))\to(\S,O)$  the minimal resolution of $\nu_{E_Q}$. We have the following formula 
for the intersection number at the surface $S(\psi)$: 
\begin{equation} \label{f:curv}
(A^\psi \cdot E_Q^\psi)_{S(\psi)}  = (\de_{L_0}^+ (Q))^{-1}  \cdot  \de_{L_0} (A). 
\end{equation}
Formula \eqref{f:curv} is a consequence of the 
renormalization formulas (see  Propositions 1.6.20 (1)  and 1.6.22 (1) of \cite{GGP19b}, applied 
to $Q$ and the points of discontinuity of the index function $\de_{L_0}$ on 
the segment $[L_0, Q)$). If   
$\de_{L_0}^+ (Q)  =  \de_{L_0} (A)$, then
\eqref{f:curv} implies that $A^\psi$ is smooth and intersects transversally $E_{Q}^\psi$. 
Then, applying \cite[Lemma 1.6.18]{GGP19b} and the renormalization formulas as before, 
we get that $A^\psi$ does not intersect any other component of $E(\psi)$. Hence 
 $A$ belongs to  $\mathcal{C}_{E_{Q}} (\psi)$. 
\end{proof}

\begin{remark} \label{rem: EP}
Let us give another description of 
the divisorial valuation $\nu_{E_P}$ associated to a rational point $P$ of $\Theta_{L_0} (D)$
(see \cite[Section 4.7]{Bu19}). 
By Lemma \ref{lem:levels} there is a branch $L'$ of $D$ such that $P \in [L_0, L']$ and  $\de_{L_0} (P) = \de_{L_0} (L')$.
Let us denote by $Q$ the point $Q_P$ introduced in Definition \ref{def:rootlevel}.
We distinguish two cases: 

- If $Q = L_0$  then $\de_{L_0} (P) = 1$ 
and $(L_0, L')$ is a cross at $(S,O)$. 
Take $(x_0, y_0)$ a local coordinate system such that $Z(x_0) = L_0$ 
and 
$Z(y_0) = L'$. 
Then, $\nu_{E_P}$ is the monomial valuation of $R$ 
determined by $\nu_{E_P} (x_0) = n_P$ and $\nu_{E_P} (y_0) = m_P$.

- If $Q \ne L_0$, then $\de_{L_0}^+ (Q) = \de_{L_0} (P) > 1$. We denote by $\psi$ the minimal resolution of $\nu_{E_{Q}}$. 
By Lemma \ref{lem:curvetta}, we have that $L'$ belongs to $\mathcal{C}_{\psi} (E_{Q})$. Let $(u, v)$ be a local coordinate system 
on $S(\psi)$ defining the cross $(E_{Q}, (L')^{\psi})$, 
that is, $E_{Q} =Z(u)$ and $(L')^{\psi} =Z(v)$ on this model.
Then, the 
divisorial valuation $\nu_{E_P}$ is \emph{monomial} in terms of the cross
$(E_{Q}, (L')^{\psi})$.  
This is a particular case of the \emph{quasi-monomial} valuations considered in \cite{FJ}. 
For any $0 \ne h \in R$, if 
one has the expansion $h \circ \psi = \sum a_{s,t} u^s v^t$ in 
$\k [[u,v]]$, then
$
\nu_{E_P} (h) = \min \{ s n_P  + t m_P  \mid a_{s, t} \ne 0 \}$. 

In both cases, 
the pair $(n_P, m_P) \in \mathbb{N}^2$ is the one determined by \eqref{eq: Qrenor}.
\end{remark}

\medskip

In the following definition we use the notion of minimal regularization of a 
fan with respect to a rank two lattice (see \cite[Section 1.3.1]{GGP19b}).

\begin{definition} \label{not:dis} Let $P_1, P_2 \in \mathcal{A}_{L_0}^C (D)$ be consecutive points in $\Theta_{L_0} (D)$ with $P_1 \preceq_{L_0} P_2$. Denote $(n_{L_j}, m_{L_j}) := (0,1)$, for 
$j \in \{ 1, \dots, m \}$. 
We set 
\[
\sigma_{P_1, P_2} := 
\left\{ 
\begin{array}{lcl}
\R_{\geq 0} (n_{P_1}, m_{P_1}) + \R_{\geq 0} (n_{P_2}, m_{P_2}) & \mbox{ if } &  P_1 \ne Q_{P_2},
\\
\R_{\geq 0} (1, 0) + \R_{\geq 0} (n_{P_2}, m_{P_2}) & \mbox{ if } &  P_1 = Q_{P_2}.
\end{array}
\right.
\]

A point $P$ of $\Theta_{L_0} (D)$ is \emph{distinguished by $C$} if 
$P \in \mathcal{A}_{L_0}^C (D)$ or 
$P \in (P_1, P_2)$ where $P_1$ and $P_2$ are consecutive points of $\mathcal{A}_{L_0}^C (D)$
and $(n_P, m_P)$ defines a ray of the minimal regularization of the fan of the cone $\sigma_{P_1, P_2}$ with respect to the lattice $\Z^2$ (see Notation \ref{not:aug}).
\end{definition}

The minimal embedded resolution $\psi_C$
of $C$ can be seen as a toroidal embedded resolution 
associated to a smooth branch
in such a way that the sequence of 
auxiliary curves defining local coordinates in the process is $L_0, \dots, L_m$   (see  \cite[Theorem 3.12]{LO}).
 As a consequence of \cite[Prop. 1.4.35 and Th. 1.6.27]{GGP19b} 
we deduce
the following result:

\begin{proposition} \label{prop:exc-dis2}
The dual graph of the total transform of $D$ by the minimal embedded resolution 
of $C$ is isomorphic to the tree $\Theta_{L_0} (D)$ with marked vertices running through the set of points distinguished by $C$, by an isomorphism which preserves the labels by
the branches of $D$
and which sends a rational point $P$
distinguished by $C$ to the vertex 
of $G(\psi_C, D)$ labelled by $E_P$.
\end{proposition}

\begin{remark} \label{rem:aug}
We have also a similar notion of \textit{toroidal pseudo-resolutions} of $C$ with respect to $L_0$. By 
\cite{GGP19b} 
one has a toroidal pseudo-resolution $\phi$ of $C$ such that 
the dual graph of the total transform of $D$ by $\phi$ is isomorphic to 
the tree $\Theta_{L_0} (D)$ with vertices running through the set of 
augmented marked points $\mathcal{A}_{L_0}^C (D)$.
\end{remark}

\begin{corollary}\label{lem:genseqratpt} (see  \cite[Theorem 4.160]{Bu19})
Let $P \in \Theta_{L_0} (D)$ be a rational point.  Then, the sequence $x_0,\ldots,x_m$ is a generating sequence of the 
divisorial valuation $\nu_{E_P}$.
\end{corollary}

\begin{proof}
By Lemma \ref{lem:levels} there is a branch $L' \in \{ L_1, \dots, L_m \}$ such that 
$\de_{L_0} (P) = \de_{L_0} (L')$. Let us take $Q = Q_P$. 
By definition we have that $Q$ and $L' $ belong to 
$ \mathcal{A}_{L_0}^D (D)$.  By Remark 
\ref{rem: EP} the valuation $\nu_{E_P}$
is monomial with respect to the cross 
defined by $E_Q$ and $L'$ at some model $S(\psi)$ (where $E_Q = L_0$ and  $S(\psi) = S$
if $Q = L_0$).
Take a minimal model $\psi'$ dominating $\psi_C$ where $E_P$
appears. Then, the dual graph $G(\psi', D)$ is obtained 
from $G(\psi_C, D)$ by subdividing the unique path joining 
the vertex labelled by $E_Q$ with the vertex labelled by $L'$.
This implies that $L_0, \dots, L_m$ is a sequence of maximal contact curves of $\psi'$. Therefore,
$x_0, \dots, x_m$
is a minimal generating sequence of the tuple 
$(\nu_{E_P}, \nu_{E_{C_1}}, \dots, \nu_{E_{C_r}})$ by Proposition \ref{gen-cur}. By Remark \ref{rem:globalgivesgs}, $x_0, \dots, x_m$
is a generating sequence of $\nu_{E_P}$.
\end{proof}

In the following lemma we describe the value of 
$\nu_{E_P}$ on $L_0, \dots, L_m$. 

\begin{lemma} \label{lem:linear}
Let $P_1, P_2 \in \mathcal{A}_{L_0}^C (D)$ be consecutive points in $\Theta_{L_0} (D)$ with $P_1 \preceq_{L_0} P_2$.
For $P \in [P_1, P_2]$, $P \ne Q_{P_2}$ and for 
 $0 \leq j \leq m$, there exists non negative integers $a_j, b_j$
such that 
\[
\nu_{E_P} (L_j)= a_j n_P+b_j m_P.
\]
In addition, if $P_1 = Q_{P_2}$ then we have $ 
\nu_{E_{P_1}} (L_j) = a_j
$, for $j=0, \dots, m$. 
\end{lemma}

\begin{proof}
Let us set $Q = Q_{P_2}$.
By Remark \ref{rem: QP2} we have that 
 $Q \preceq_{L_0} P_1$
with equality if $P_1 = L_0$ or if $\de_{L_0} (P_1) < \de_{L_0}(P_2)$.
Let  $P \in [P_1, P_2]$  be a rational point, $P \ne Q$. 
By \cite[Prop. 8.16 (3) and Def. 3.14]{GGP19} we have that 
\begin{equation} \label{eq: L0}
\nu_{E_P} (L_0) = 
\de_{L_0}^+ (P) \stackrel{\eqref{eq: indexplus}}{=}
n_P  \, \de_{L_0} ( P) = n_P  \, \de_{L_0} ( P_2). 
\end{equation}
We set $a_0 := \de_{L_0} ( P_2)$ and $b_0 := 0$.

Take $j \in \{1, \dots, m\}$.
We distinguish two cases as in the proof of Lemma \ref{block1}.

- 
If $P_2 \prec_{L_0} L_j$ then by \eqref{aux2} we get $
V^{P}_{L_0} (L_j) =
(L_0 \cdot L_j)  \left( {\ic}_{L_0} (Q) + \frac{1}{\de_{L_0} ( P_2)} 
 \left(\ex_{L_0} (P) - \ex_{L_0} (Q)  \right) \right)$.
By \eqref{eq: Qrenor} we have that
$
\ex_{L_0} (P) - \ex_{L_0} (Q) =    \frac{1}{ \de_{L_0} ( P_2)} \cdot \frac{m_P}{n_P}$.
Since $(L_0 \cdot L_j) = \de_{L_o} (L_j)$, 
we obtain: 
\[
V^{P}_{L_0} (L_j) = \de_{L_o} (L_j) \left( {\ic}_{L_0} (Q) + \frac{1}{\de_{L_0}^2 ( P_2)}  \frac{m_P}{n_P}
\right). 
\]
By \eqref{eq: EP} and \eqref{eq: L0} it follows that 
$
\nu_{E_P} (L_j) = n_P  \, \de_{L_0} ( P_2) \, \de_{L_0} (L_j) \ic_{L_0} (Q)
+ m_P \frac{\de_{L_0} (L_j)}{\de_{L_0} (P_2)}$.
Notice that 
$b_j := \frac{\de_{L_0} (L_j)}{\de_{L_0} (P_2)}$ is an integer since
$P_2 \preceq_{L_0} L_j$.
Let us check that the coefficient $a_j:= \de_{L_0} ( P_2) \, \de_{L_0} (L_j) \ic_{L_0} (Q)$ is also an integer. 
If $Q = L_0$, then $a_j=0$. If $Q \ne L_0$,  let us take 
a branch $A$ 
such that $\langle L_0, L_j , A \rangle = Q$ and 
$\lim_{P' \in (Q, A), \, P' \to Q} \de_{L_0} (P') = \de_{L_0}^+ (Q) =  \de_{L_0} (A)$.
Since $\de_{L_0}^+ (Q) = \de_{L_0} (P_2)$ by definition,
we get from \eqref{eq:int3} that $a_j = (L_j \cdot A)$.

\medskip 

- Assume now that condition $P_2 \prec_{L_0} L_j$ does not hold. 
By Lemma \ref{lem:levels}
we can take 
$L' \in \{ L_0, \dots, L_m  \}$
such that $P_2 \in [L_0, L']$ and $\de_{L_0} (P_2) = \de_{L_0} (L') $. 
This implies that 
$\langle L_0, L_j , P \rangle =  \langle  L_0, L_j , L' \rangle$.
By \eqref{VPA}
we get $
V^{P}_{L_0} (L_j) =
\de_{L_0} ( L_j)    {\ic}_{L_0} ( \langle L_0, L_j , L' \rangle  ) $.
By \eqref{eq: EP} and \eqref{eq: L0} 
we obtain
\[
\nu_{E_P} (L_j) 
=  n_P  \, \de_{L_0} ( L')  \de_{L_0} ( L_j)    
{\ic}_{L_0} ( \langle L_0, L_j , L' \rangle ).
\]
 Then, by 
\eqref{eq:int3}
the coefficient 
$a_j:=  \de_{L_0} ( L' )  \de_{L_0} ( L_j)    
{\ic}_{L_0} ( \langle L_0, L_j , P \rangle )$ is equal to $(L_j, L')$, while $b_j :=0$.

\medskip 

Finally, let us check the assertion in the case $P_1 = Q$.
By \cite[Prop. 8.16 (3) and Def. 3.14]{GGP19} we have that 
$ 
\nu_{E_Q} (L_0) = \de_{L_0}^+ (Q)
$.
We have also that $\de_{L_0}^+(Q) = \de_{L_0} (P_2)$ by hypothesis.
By \eqref{eq: EP} and the discussion in the previous cases we obtain 
$\nu_{E_Q} (L_j) = a_j$.
\end{proof}

We will apply Lemma \ref{lem:linear} and Proposition \ref{prop:linear2}
in Section \ref{sec:min_res} below. 

\begin{proposition} \label{prop:linear2}
Let $P_1, P_2 \in \mathcal{A}_{L_0}^C (D)$ be consecutive points in $\Theta_{L_0} (D)$.
Take the integers $a_j, b_j \in {\mathbb N}$ given by Lemma \ref{lem:linear}, for $0 \leq j \leq m$.
Then, the map
\[
\phi_{P_1, P_2} \colon \sigma_{P_1, P_2} \cap \Z^2 \longrightarrow (\R_{\geq 0} w^{P_1} + \R_{\geq 0} w^{P_2} ) \cap \Z^{m+1}, \quad 
 (p, q) \mapsto (a_0 p + b_0 q, \dots, a_m p + b_m q ), 
\]
is an isomorphism. 
\end{proposition}
\begin{proof} 
Let us show first 
that the group homomorphism 
\begin{equation} \label{eq: grouphom}
\Z^2  \longrightarrow \Z^{m+1}, \quad 
 (p, q) \mapsto (a_0 p + b_0 q, \dots, a_m p + b_m q ), 
\end{equation}
is injective and its image is a direct factor of $\Z^{m+1}$.
Its associated matrix $G$, with respect to the canonical basis,
has rows $(a_j, b_j)$ for $j \in \{ 0, \dots, m \}$.
By Lemma \ref{lem:levels} we can take 
$L' \in \{ L_1, \dots, L_m  \}$
such that $P_2 \in [L_0, L']$ and $\de_{L_0} (P_2) = \de_{L_0} (L') $. 
By Lemma \ref{lem:mgs},  there exists $1  \leq g \leq m-1$ such that, up to relabelling $L_1, \dots, L_m$, one has $L' =L_{g+1}$ and 
$L_0, \ldots, L_{g}$ is a sequence of maximal contact curves of  $\psi_{L'}$.
By the proof of Lemma \ref{lem:linear} we have that 
\[ 
a_0 = \de_{L_0} (P_2) = \de_{L_0} (L') = (L_0 \cdot L')
\]
and
\[
a_j = (L_j \cdot L') \mbox{ and } b_j = 0 \mbox{ for } j \in \{ 1, \dots, g \}.
\]
We get by the first case of the proof of Lemma \ref{lem:linear} that 
$b_{g+1} = \frac{\de_{L_0} (L')}{\de_{L_0} (P_2) } = 1$. 
This implies that
\[
\left|
\begin{array}{cc}
a_j & b_j
\\
a_{g+1} & b_{g+1}
\end{array}
\right|= a_j, \mbox{ for } j \in \{ 0, \dots, g \}.
\]
The numbers $a_0, \dots, a_g$ are coprime since they form a system of generators 
of the semigroup $\Gamma_{L'}$ (see  Proposition \ref{gen-cur}).
This implies that the matrix $G$ has rank two, and 
its minors of rank two are coprime, hence 
the homomorphism \eqref{eq: grouphom} is injective 
and its image is a direct factor of $\Z^{m+1}$.

  By \eqref{defwp} and Lemma  \ref{lem:linear}, the vector
  $w^{P_2}$ is a positive multiple of    $\phi_{P_1, P_2} (n_{P_2}, m_{P_2})$.
  If $Q_{P_2} \ne P_1$ (resp. $Q_{P_2} =P_1$)  
  the same happens when we compare the vector $w^{P_1}$ with  $\phi_{P_1, P_2} (n_{P_1}, m_{P_1})$  (resp.  $\phi_{P_1, P_2} (1, 0)$).
Taking into account Definition \ref{not:dis}, 
  this implies that the map  $\phi_{P_1, P_2}$ is an isomorphism of semigroups.
\end{proof}

%---------------------------------------------------

\begin{remark} \label{rem:ray}

Take $j \in \{0, \dots, m \}$. By Proposition \ref{prop:exc-dis2} the exceptional divisor 
$E_{C_j}$ is the representing divisor of the attaching point $Q_j$
of $C_j$ to the tree $\Theta_{L_0} (D)$, i.e., $ E_{C_j}  = E_{Q_j}$.
This implies that  $(L_i \cdot C_j) = \nu_{E_{C_j}} (L_i)$, for $i =0, \dots, m$,
since the sequence of maximal contact curves $L_0, \dots, L_m$ is
generic for $C$ 
(see \eqref{f:generic} and Remark \ref{rem:val-int}).
Thus, we get
$\ord (\eta_j ) = ( \nu_{E_{C_j}} (x_0), \dots , \nu_{E_{C_j}} (x_m))$.
Then, $\ord (\eta_j )$
defines a ray of the fan $\cT$ since $Q_j \in \mathcal{A}_{L_0}^C (D)$
(see Proposition \ref{prop:complexT}).
\end{remark}

%------------------------------------------------------------------------

\subsection{Initial ideals associated to the local tropicalization and toric resolutions}\label{sec:lastone}

In this section we consider the embedding $S= \k^2  \subset \k^{m+1}$ 
defined by the map  \eqref{emb-plane}. 
We describe the initial ideals $in_w (I_S)$ when $w \in \Trop_{\geq0}(I_S)$. 
Then, we prove that there exists a regular subdivision $\Sigma$ of the positive quadrant 
$\Sigma_{0, m+1}$
such that the restriction of $\pi_{\Sigma}$ to the strict transform of $S$ 
is an embedded resolution of the given plane curve singularity $C\subset S$. 

\medskip 

Given $\alpha=(\alpha_0,\ldots,\alpha_{m})\in\mathbf Z_{\geq0}^{m+1}$, 
we denote $X^\alpha=X_0^{\alpha_0}\cdots X_{m}^{\alpha_{m}}$. 
Then we can write a power series $g\in\k[[X_0,\ldots,X_{m}]]$, $g\ne 0$, as 
$g=\sum_{\alpha\in\mathbf Z_{\geq0}^{m+1}}{c_\alpha X^\alpha}$. 
Let $w\in\R_{\geq0}^{m+1}$. 
The $w$-\emph{weight} of $g$ is 
$$\nu_w(g):=\min\set{w\cdot\alpha\;/\;c_\alpha\neq0},$$ 
where $\cdot$ stands for the usual scalar product of $\R^{m+1}$. 
The $w$-\emph{initial part} of $g$, denoted $in_w(g)$, 
is the sum of all terms in $g$ of $w$-weight $\nu_w(g)$. Set 
$\nu_w(0)=\infty$ and $in_w(0)=0.$ 
Actually, the map 
$\nu_w$ defines a valuation that is called the monomial valuation of weight 
$w.$ Such a valuation is divisorial if and only if $w\in \Q^{m+1}.$
 Similarly,  if $\prec$ is a monomial order on  $\k[X_0,\dots,X_{m}]$,
 the $\prec$-initial part of $g$ is 
$c_{\alpha_0} X^{\alpha_0}$, where $\alpha_0 = 
\min_{\prec} \{ \alpha \mid c_{\alpha} \ne 0 \}$.

\medskip

Let $I$ be an ideal of $\k[[X_0,\ldots,X_{m}]]$. 
The $w$-\emph{initial ideal} of $I$ is the ideal generated by the set 
$\set{in_w(g)\;/\;g\in I}$. We denote this ideal by $in_w(I)$. 
Note that the $w$-initial parts of a system of generators 
$g_1,\ldots,g_\ell$ of $I$ do not generate in general $in_w(I)$. 
When this occurs, $g_1,\ldots,g_\ell$ is a $w$-\emph{standard basis} 
of $I$.

\medskip

Given $w=(w_0,\ldots,w_m)\in\R_{\geq0}^{m+1}$, we call $\text{supp}(w)$ 
the set of indices $i$ such that $w_i \neq 0$. We define $C_w(I)$ as 
the closure  of the
set:
$$\set{w'\in\R_{\geq0}^{m+1}\mid \text{supp}(w)=\text{supp}(w') 
\text{ and }in_w(I)=in_{w'}(I)},$$
with respect to the usual Euclidean topology on $\R^{m+1}$.
Then, the family $\set{C_w(I)\mid w\in\R_{\geq0}^{m+1}}$ is finite and forms 
a rational polyhedral fan $\GF(I)$ subdividing the fan $\Sigma_{0,m+1}$ 
(see \cite{BT07,PS13,Tou05}). This fan is called the \emph{Gr\"{o}bner fan} 
of $I$ and is a local version of the notion of Gr\"{o}bner fan studied 
in \cite{MR88}.

\medskip

The following result relates the local tropicalization and the 
Gr\"{o}bner fan of $I$. For the proof, we refer the reader to 
\cite[Theorem 5.8]{Tou05} and \cite[Theorem 11.2]{PS13}.

\begin{proposition} \label{prop:nomonomials}
Let $I$ be an ideal of $\k[[X_0,\ldots,X_m]]$. 
\begin{enumerate}
\item Given $w\in\R_{\geq0}^{m+1}$, 
$w\in \Trop_{\geq0}(I)$ if and only if $in_w(I)$ contains no monomial.
\item $\Trop_{\geq0}(I)\cap\R_{\geq0}^{m+1}$ is a union of cones of 
$\GF(I)$.
\end{enumerate}
\end{proposition}

The following definition 
is based upon that of non-degenerate functions 
with respect to their Newton polyhedra of 
Khovanskii and Kouchnirenko \cite{Kou,Kho}.
\begin{definition}[\cite{AGS,Te1}]\label{Newton}
An ideal $I\subset\k[[X_0,\ldots,X_m]]$ is \emph{Newton non-degenerate} if 
for every $w\in\R_{\geq0}^{m+1}$, the variety $Z(in_w(I))$ does not have 
singularities in the torus $(\k^*)^{m+1}$. 
\end{definition}
The notion of Newton non-degenerate given in Definition \ref{Newton} may 
be seen as a particular case of the notion of sch\"{o}n compactification of a subvariety of a torus introduced in 
\cite{Te1}.
\begin{remark}\label{rem:rationals}
It is a straightforward consequence of Proposition 
\ref{prop:nomonomials} and the Nullstellensatz that 
$\Trop_{\geq0}(I)\cap\R_{\geq0}^{m+1}$ can also be described as the set of all 
$w\in\R_{\geq0}^{m+1}$ such that 
$Z(in_w (I))\cap(\k^*)^{m+1}\neq\emptyset$. Hence in order to 
decide whether the ideal $I$ is Newton non-degenerate it is enough to 
verify that the condition given in Definition \ref{Newton} is satisfied 
for those $w$ in $\Trop_{\geq0}(I) \cap\R_{\geq0}^{m+1}$. In addition, 
it is sufficient to check the condition for every vector in 
$\Trop_{\geq0}(I) \cap\Q_{\geq0}^{m+1}$. 
Indeed, if $w\in\Trop_{\geq0}(I) \cap\R_{\geq0}^{m+1}$ and 
$w\notin \Q_{\geq0}^{m+1}$ then, by Proposition \ref{prop:nomonomials}, 
there exists a cone $\sigma\in\GF(I)$ such that 
$w\in\stackrel{\circ}{\sigma} \subset \Trop_{\geq0}(I)$. Since $\sigma$ 
is a rational cone, we can take $w'\in\stackrel{\circ}{\sigma}\cap\Q_{\geq0}^{m+1}$ 
and we get $in_w(I)=in_{w'}(I)$.
\end{remark}

\begin{theorem}[See \cite{AGS,Te1}] \label{res:surf} 
Let $I$ be a Newton non-degenerate ideal of $\k[[X_0,\ldots,X_m]]$. 
Let $\Sigma$ be a regular subdivision of the Gr\"{o}bner fan of $I$ and let 
$\pi^\Sigma_{\Sigma_{0,m+1}}\colon X_\Sigma\to\k^{m+1}$ be the associated 
toric modification. 
Then, the strict transform of $Z(I)$ by $\pi^\Sigma_{\Sigma_{0,m+1}}$ 
is non-singular and transversal to the orbit stratification of 
the exceptional locus of 
$\pi^\Sigma_{\Sigma_{0,m+1}}$.
\end{theorem}

In what follows we keep the notation of section \ref{sec:localtrop}. 
Let $P$ be a rational point of $\Theta_{L_0}(D)$. Then, $w^P$ belongs 
to $\Trop_{\geq0}(I_S)\cap\Q_{>0}^{m+1}$. Let us relabel the functions $x_0,\ldots,x_m$ in such a way that $x_0, \dots, x_{g}$ is a minimal 
generating sequence of the divisorial valuation $\nu_{E_P}$ (see Corollary \ref{lem:genseqratpt}). With the notation of Section \ref{sec:dival}, 
$g=g(E_P)$. Next we give a system of generators $H_2,\ldots,H_m$ 
of $I_S$ that is well adapted to the computation of $w^P$-initial parts:
\medskip 

\begin{itemize}
\item For $j\in\set{1,\ldots,g-1}$, we consider the 
$(x_0,\ldots,x_j)$-expansion of $x_{j+1}$ given in Proposition 
\ref{prop:adicexpxj} and set 
\begin{equation}\label{ISgenerators1}
H_{j+1}:=-X_{j+1}+X_{j}^{n_j}-\theta_j \cdot  X_{0}^{b_{0}^j}
X_{1}^{b_{1}^j}\cdots X_{j-1}^{b_{j-1}^j}
+p_j(X_0,\ldots,X_j).
\end{equation}
\item For $j\in\set{g,\ldots,m-1}$, we set 
\begin{equation}\label{ISgenerators2}
H_{j+1}:=-X_{j+1} + \sum_{I = (i_0, \dots, i_{g})} 
c_{j,I }X_0^{i_0} X_1^{i_1} \cdots X_{g}^{i_{g}},
\end{equation}
where $x_{j+1}=\sum_{I = (i_0, \dots, i_{g})} 
c_{j,I }x_0^{i_0} x_1^{i_1} \cdots x_{g}^{i_{g}}$ is 
the $(x_0, \dots, x_g)$-adic expansion of 
$x_{j+1}$ (see Remark 
\ref{rem:expxj}).
\end{itemize}

\medskip

By definition, all the $H_{j+1}$ are in the kernel $I_S$ of the 
homomorphism \eqref{emb-plane}. They generate $I_S$ because 
any element of $I_S$ in congruent modulo $(H_2,\ldots,H_m)$ to 
some element of $\k[[X_0,X_1]]$ and the restriction of the 
homomorphism \eqref{emb-plane} to this ring is injective.

\begin{lemma}\label{lem:inpartH} 
Let $P$ be a rational point of $\Theta_{L_0}(D)$. With the previous 
notation, we have:
\begin{enumerate}
\item\label{enum:if1} For $j\in\set{1,\ldots,g-1}$,  
$
in_{w^P}(H_{j+1})=X_{j}^{n_j}-\theta_j \cdot  X_{0}^{b_{0}^j}
X_{1}^{b_{1}^j}\cdots X_{j-1}^{b_{j-1}^j}$.

\item\label{enum:if2} For $j\in\set{g,\ldots,m-1}$, $
in_{w^P}(H_{j+1})=-X_{j+1} + in_{w^P}\left(\sum_{I = (i_0, \dots, i_{g})} 
c_{j,I }X_0^{i_0} X_1^{i_1} \cdots X_{g}^{i_{g}}\right)$.
\end{enumerate}
\end{lemma}

\begin{proof}
With the notation of Section \ref{sec:expgensec}, one has that 
$w^P=\left(1,\g_1/\g_0,\ldots,\g_g/\g_0,w^P_{g+1},\ldots,w^P_m\right),$
where $w^P_i=\nu_{E_P}(x_j)/\g_0$, for $g+1\leq i\leq m$. 
Assertion \eqref{enum:if1} follows from Proposition \ref{prop:adicexpxj}. 
Assertion \eqref{enum:if2} follows from Remark \ref{rem:expxj} 
applied to $h=x_{j+1}$, for $g\leq j\leq m-1$.
\end{proof}

\begin{remark} If $P$ is a rational point of $\Theta_{L_0}(D)$, then the vector $\trop(\nu_{E_P})=\nu_{E_P}(x_0)w^P$ is primitive (see the proof of Lemma \ref{lem:inpartH}).
\end{remark}

\begin{proposition} \label{prop:Hilbert}
Let $P$ be a rational point of $\Theta_{L_0}(D)$. With the previous 
notation, 
\begin{equation} \label{eq: Hilbert}
in_{w^P}(I_S)=(in_{w^P}(H_2),\ldots, in_{w^P}(H_m)).
\end{equation}
\end{proposition}

\begin{proof}
Let $\prec$ be the monomial order on $\k[X_0,\ldots,X_m]$ defined 
by $X^\alpha\prec X^\beta$ if $w^P\cdot\alpha<w^P\cdot\beta$ or 
$w^P\cdot\alpha=w^P\cdot\beta$ and $X^\alpha\prec_{lex}X^\beta$, 
where $\prec_{lex}$ is the lexicographic order with respect to 
$X_m<\cdots<X_1<X_0$. By Lemma \ref{lem:inpartH}, the $\prec$-initial 
parts of $H_2,\ldots,H_m$ are pairwise coprime. Therefore, 
$H_2,\ldots,H_m$ is a $\prec$-standard basis of $I_S$. 
By \cite[Theorem 2.1]{KTV}, $H_2,\ldots,H_m$ is also a $w^P$-standard 
basis of $I_S$. 
\end{proof}

\begin{proposition} \label{prop:CI}
Let $P$ be a rational point of $\Theta_{L_0}(D)$. 
Then, 
the  singular locus 
of $Z(in_{w^P}(I_S))$ does not meet the torus $(\k^*)^{m+1}$.
\end{proposition}
\begin{proof}
The variety $\mathcal C \subset \k^{g}$ defined by 
the ideal $(in_{w^P}(H_2),\ldots,in_{w^P}(H_g))$, seen as an ideal of $\k[X_0, \dots, X_{g-1}]$, is a monomial curve with semigroup generated by $\g_0/n_g,\ldots,\g_{g-1}/n_g$ (see \cite{T3}). 
By the form of $in_{w^P}(H_{j+1})$, for $j=g,\ldots,m-1$, we get that the restriction of 
the projection $(X_0, \dots, X_m) \mapsto (X_0, \dots,X_g)$ to the variety $Z(in_{w^P}(I_S))$
is an isomorphism whose image is the irreducible surface $\mathcal C \times \k \subset\k^{g}\times\k$
defined by the ideal $(in_{w^P}(H_2),\ldots,in_{w^P}(H_g))$ of  $\k[X_0, \dots, X_{g}]$. 
 As the singular locus of 
$\mathcal C \times \k $
does not meet the torus $(\k^*)^{g+1}$ the assertion follows. 
\end{proof}

Next we prove the following:

\begin{proposition}\label{prop:ISisNnd}
The ideal $I_S\subset\k[[X_0,\ldots,X_m]]$ is Newton non-degenerate.
\end{proposition}
\begin{proof} 
By Remark \ref{rem:rationals}, it is enough to prove that the singular 
locus of $Z(in_w(I_S))$ and the torus $(\k^*)^{m+1}$ have empty intersection 
for every $w\in\Trop_{\geq0}(I_S)\cap\Q_{\geq 0}^{m+1}$. 
Take $w\in\Trop_{\geq0}(I_S)\cap\Q_{> 0}^{m+1}$. By Theorem \ref{loc-trop}, 
it must be of the form $w=w^P$ for some interior point $P$ of $\Theta_{L_0}(D)$. Moreover, $P$ is a rational point (see Lemma \ref{lem:charrat}). According to Proposition \ref{prop:CI}, we have that
$\mathrm{Sing}(Z(in_w(I_S)))\cap(\k^*)^{m+1}=\emptyset$.

Now take $w\in\Trop_{\geq0}(I_S)\cap\Q_{\geq0}^{m+1}$ such that $w_i=0$ 
for some $i \in\set{0,\ldots,m}$. By Lemma 
\ref{lem:tropU}.\eqref{item:tropU1}, up to replacing $w$ by  
a positive multiple, we can assume that it
is one of the vectors of 
the canonical basis of $\R^{m+1}$. By definition, 
$w=w^{L}$ for some $L \in \{ L_0, \dots, L_m \}$. 
Up to relabelling $L_0,\ldots,L_m$ we can assume that 
$L = L_{g+1}$
for some $g \in \{ 1, \dots, m-1 \}$, and $L_0, \dots, L_{g}$ is a sequence of maximal 
contact curves of $\psi_{L}$ (see Lemma \ref{lem:mgs}), and then $w = e_{g+1}$.
Take 
$H_2, \dots, H_m $ defined by \eqref{ISgenerators1} and \eqref{ISgenerators2}. 
By Proposition \ref{expansions} 
applied to  $E = E_{L}$, we get that
\[
H_{g+1}:=-X_{g+1}+X_{g}^{n_g}-\theta_g \cdot  X_{0}^{b_{0}^g}
X_{1}^{b_{1}^g}\cdots X_{g-1}^{b_{g-1}^g}
+p_g(X_0,\ldots,X_g).
\]

Let us check that \eqref{eq: Hilbert} holds also in this case.  If $j \ne g +1$ the variable $X_{g+1}$ does not appear on 
$H_j$,  hence $ in_w (H_j) =  H_j$, while
$in_w (H_{g+1}) = H_{g+1} + X_{g+1}$.
Denote by $\g_0, \dots, \g_{g}$ the minimal sequence of
generators of $\Gamma_{L}$ and 
set $w' = (\g_0, \dots, \g_{g}, 0, \dots, 0) \in \Z^{m+1}$.
Let $\prec'$ be the monomial order on $\k[X_0,\ldots,X_m]$ defined 
by $X^\alpha\prec' X^\beta$ if $w'\cdot\alpha<w'\cdot\beta$ or 
$w'\cdot\alpha=w'\cdot\beta$ and $X^\alpha\prec_{lex}X^\beta$, 
where $\prec_{lex}$ is the lexicographic order with respect to 
$X_m<\cdots<X_1<X_0$.
Consider then the monomial order $\prec$ 
defined by $X^\alpha\prec X^\beta$ if $w\cdot\alpha<w\cdot\beta$ or 
$w\cdot\alpha=w\cdot\beta$ and $X^\alpha\prec' X^\beta$.
We obtain that: 
\[
in_{\prec} (H_2) = X_1^{n_1}, \, \dots, \,
in_{\prec} (H_{g+1}) = X_g^{n_g}, \, in_{\prec} (H_{g+2}) = -X_{g+2}, \, \dots, \, 
in_{\prec} (H_{m}) = -X_{m}.
\]
As the initial forms $in_{\prec} (H_j)$, for $j \in \{2, \dots,m\}$
are coprime, 
we end the proof of \eqref{eq: Hilbert} in this case by arguing as
 in the proof of Proposition \ref{prop:Hilbert}. 

Notice that the algebroid subvariety of $\k^{g+2}$
defined by the ideal $(in_w (H_2), \dots, in_w (H_{g+1}))$  is the product of the embedding of $L \subset \k^{g+1}$
defined by $(x_0, \dots, x_g)$ times an affine line. 
Then, we check that 
$\mathrm{Sing}(Z(in_w(I)))\cap(\k^*)^{m+1}=\emptyset$
by arguing as in the proof of Proposition \ref{prop:CI}.
\end{proof}

\begin{theorem} \label{th-res}
Let $C=\cup_{j=1}^r C_j$ be a reduced plane curve singularity. Take a sequence  $L_0, \dots, L_m$
of maximal contact curves of the minimal embedded resolution of $C$, which is generic for $C$.
Let $\Sigma$ be a regular subdivision of $\GF(I_S)$ inducing a 
regular subdivision $\cT_{reg}$ of 
the fan $\cT$. Then the restriction of 
$\pi^{\Sigma}_{\Sigma_{0,m+1}}$ to the strict transform of $Z(I_S)$ 
induces an embedded resolution of $C \subset S$.
\end{theorem}
\begin{proof}
By Remark \ref{rem:ray}, the fan $\cT$  contains the rays
spanned by $\ord (\eta_j)$ 
for all the branches $C_j$ of $C$. 
Let us denote by $\pi: S(\pi)  \to S$ the restriction of $\pi^\Sigma_{\Sigma_{0,m+1}}$ to 
the strict transform $S(\pi)$ of $Z(I_S)$.
By Proposition \ref{prop:ISisNnd} and Theorem \ref{res:surf}, 
$S(\pi)$ is non-singular and transversal to the orbit 
stratification of the exceptional locus. 
In addition, by  Theorem \ref{th:resgenplane} the modification 
$\pi^\Sigma_{\Sigma_{0,m+1}}$ is a toric 
embedded resolution of 
$C \subset \k^{m+1}$. 
Since $C \subset S$ it follows that the strict transform 
$C^{\pi}$ of $C$ is contained in $S(\pi)$ and 
the transversality properties above imply 
that $C^{\pi}$ is transversal to the exceptional divisor 
of $\pi$. 
\end{proof}

\subsection{Minimal embedded resolution with one toric morphism} \label{sec:min_res}

In this section we prove that 
when the regular subdivision  $\cT_{reg}$
in Theorem \ref{th-res} is the minimal one, then the restriction of 
$\pi^{\Sigma}_{\Sigma_{0,m+1}}$ to the strict transform of $Z(I_S)$ 
induces the minimal embedded resolution of $C \subset S$.

\medskip

Let us denote by $\cT_{reg}^{min}$ the minimal regularization 
 of the two dimensional fan $\cT$.
Its projectivization is the set of images of its non-zero cones 
in the projective space $\mathbb{P}^m (\R)$.
The projectivization of $\cT_{reg}^{min}$ is a graph 
since the fan $\cT_{reg}^{min}$ consist of cones 
of dimension at most two contained in $\R^{m+1}_{\geq 0}$. 
The following result  is a consequence of Propositions \ref{prop:exc-dis2} and \ref{prop:linear2},  taking into account Definition \ref{not:dis}
and the definition of $\cT$ (see Proposition \ref{prop:complexT}).

\begin{proposition} \label{prop:exc-dis} 
Let $\psi$ be the minimal embedded resolution of $C$. 
The projectivization of the fan $\cT_{reg}^{min}$
is isomorphic to the dual graph of the total transform of $D$ under $\psi$, by an isomorphism 
which respects the labels by the components of $D$. 
\end{proposition}

%-------------------------------------------------

Keep the notation of Theorem \ref{th-res}. Let us study 
the intersection of $S(\pi)$ with 
the orbits defined by cones of $\mathcal{T}_{reg}$. This 
will allow us to describe when $\pi:S(\pi)\to S$ is the 
minimal embedded resolution of $C\subset S$.

If $e_j$ is a vector of the canonical basis of $\Z^{m+1}$, then
the orbit $O ( \R_{\geq 0} e_j)_{\Sigma_{0,m+1}}$ is not contained in 
the discriminant locus of $\pi^\Sigma_{\Sigma_{0,m+1}}$ and 
the  intersection of $O ( \R_{\geq 0} e_j)_{\Sigma}$  with $S (\pi)$ is equal to 
the strict transform $L_j^{\pi}$ of $L_j$. 
The orbits defined by other rays of $\cT_{reg}$ are contained in 
the critical locus of $\pi^\Sigma_{\Sigma_{0,m+1}}$.
Any such ray is of the form 
$\rho =  \R_{\geq 0} w^{P} $
for a unique rational point $P$ of $\Theta_{L_0} (D)$ (see Lemma \ref{block2}).

\begin{theorem} \label{prop:div}
With the hypothesis of Theorem \ref{th-res}, assume that   $\cT_{reg} = \cT_{reg}^{min}$. Then: 
\begin{enumerate}
\item 
If  $P$ is a rational point of  $\Theta_{L_0} (D)$, and 
$\rho =  \R_{\geq 0} w^{P} \in \cT_{reg}^{min} $, then the closure of 
$S(\pi) \cap O(\rho)$
is equal to $E_P$. 
\item  The restriction of 
$\pi^{\Sigma}_{\Sigma_{0,m+1}}$ to the strict transform of $Z(I_S)$ 
induces the minimal  embedded resolution of $C \subset S$.
\end{enumerate}
\end{theorem}
\begin{proof} 
We prove first that $S(\pi) \cap O(\rho)$ is irreducible. 
In order to study this intersection, we can see the orbit $O(\rho)$ 
in an affine open toric subvariety $X_\sigma$ of $X_\Sigma$, 
where $\sigma\in\Sigma$ is spanned by a basis $v_0,\ldots,v_m$ 
of $N=\Z^{m+1}$ and $v_0:=\trop(\nu_{E_P})=\nu_{E_p} (x_0)w^{P}$.
Then, the affine toric variety $X_\sigma$ 
is isomorphic to $\k^{m+1}$ with coordinates $(U_0,\ldots,U_m)$, 
where $U_i=\chi^{\check v_i}$ for $i=0,\ldots,m$. 
In addition, the restriction $\pi_\sigma$ of the toric modification $\pi^\Sigma_{\Sigma_{0,m+1}}$ 
to $X_\sigma$ is given by 
\[X_i\circ\pi_\sigma=U_0^{v_{0,i}} U_1^{v_{1,i}} \cdots U_m^{v_{m,i}},
\text{ for }0\leq i\leq m,\]
where $v_j=\sum_{j=0}^m v_{j,i}e_i$,
and the orbit $O(\rho)$ is defined by $U_0 = 0$ and 
$U_1\cdots U_m\ne 0$.
Given $g\in\k[[X_0,\ldots,X_m]]$, we write 
$g\circ\pi_\sigma=U_0^{\nu_{v_0}(g)}\cdot g_\rho$, where 
 $g_\rho \in\k[[U_0,\ldots,U_m]]$ and $U_0$ does not 
divide $g_\rho$.
If $h=in_{v_0}(g)$, then $h_\rho \in\k[U_1,\ldots,U_m]$ and $h_\rho =g_\rho(0,U_1,\ldots,U_m)$. 
This fact and Proposition \ref{prop:Hilbert} imply 
that 
$S (\pi) \cap O (\rho)$ is defined by the ideal of $\k[U_1^{\pm 1},\ldots,U_m^{\pm 1}]$ generated by the set 
$\set{(in_{v_0}(H_2))_\rho, \dots, (in_{v_0}(H_m) )_\rho}$. 
Since the polynomials $(in_{v_0}(H_j))_\rho$ belong to $\k[U_1, \dots, U_m]$ for $j =2, \dots, m$, we get that the variety $S (\pi) \cap O (\rho)$ is irreducible if and only if 
the variety $V \subset (k^*)^{m+1}$ defined by these polynomials, 
seen in $\k [U_0^{\pm 1}, U_1^{\pm 1}, \dots, U_m^{\pm 1}]$,
is irreducible.
The variety $V$ is also defined by the ideal
$(in_{v_0}(H_2) \circ \pi_\sigma ,\ldots,in_{v_0}(H_m) \circ \pi_\sigma)$. The monomial map 
$\pi_\sigma$ induces an isomorphism of tori, therefore 
the variety $V$
is irreducible by Proposition \ref{prop:CI}. This shows that $S (\pi) \cap O (\rho)$ is irreducible.

We can take the regular cone $\sigma\in\Sigma$ in such a way that $v_1 \in |\mathcal{T}| $ and $v_2, \dots, v_m \notin |\mathcal{T}|$. 
Denote by $\rho_i $ the ray $\R_{\geq 0} v_i$ for $i=1,\ldots,m$ and set $\rho_0=\rho$. 
Take $i \in\set{ 2, \dots, m}$. The condition $v_i \notin |\mathcal{T}|$ implies that 
there exists a monomial in $\k [X_0, \dots, X_m]$ which belongs to the $v_i$-initial ideal of $I_S$
(see  Proposition \ref{prop:nomonomials}). 
Hence the intersection $S(\pi) \cap O(\rho_i)$ is empty.

Denote by $E_i$ the closure of  $S(\pi) \cap O(\rho_i)$ on $X_\sigma$ for $i=0,1$. 
As the map $\pi: S(\pi) \to S$ is a model of the smooth surface $(S,O)$, 
its dual graph $G(\pi)$ is a tree. This implies that 
the divisors $E_0$ and $E_1$ intersect at at most one point, 
which belongs to the orbit $O(\rho_0 + \rho_1)$. We show now that 
this intersection consists of exactly one point $o_1$.

Assume by contradiction that $E_0 \cap E_1 =\emptyset$. Then 
as the dual graph $G(\pi)$ is connected, 
there exists a sequence of 
rays $\rho'_0, \dots, \rho'_s$ of $\cT_{reg}^{min}$, with $s>1$, 
$\rho_0 = \rho_0'$, $\rho_1 = \rho'_s$, such that $\rho_{i-1}' + \rho_{i}' \in \cT_{reg}^{min}$, and $\rho_{i-1}' + \rho_{i}'\ne \rho_0 + \rho_1$,
for $i = 1, \dots, s$.
Then, the projectivization of the nonzero elements of the set 
 $(\rho_0 + \rho_1) \cup \bigcup_{i=1}^s (\rho_{i-1}' + \rho_{i}')$
contains a non trivial cycle. 
We get a contradiction by Proposition  \ref{prop:exc-dis}, since the projectivization of fan  $\cT_{reg}^{min}$ is a tree.

Denote by $u_j$, for $j=0, \dots, m$, the restriction of the function $U_j$ to the surface $S(\pi)$. 
We obtain that $(u_0, u_1)$ are local coordinates  at $o_1$.
It follows that the order of vanishing 
of $x_j \circ\pi_\sigma $ along the divisor $E_0$ is
equal to $v_{0,j} $, for $0\leq j\leq m$. Therefore, 
$\nu_{E_P}(x_j)=\nu_{E_0}(x_j)$, for $0\leq j\leq m$.

By the transversality of $S(\pi)$ to the toric stratification of the exceptional divisor of $\pi$ we
have that any irreducible component of $E(\pi)$ must be the closure of 
$S (\pi) \cap O (\rho')$, for $\rho'$ running through the 
rays of $\cT_{reg}$ which intersect $\R_{>0}^{m+1}$. 
By Proposition 
\ref{prop:exc-dis}
we have a bijection 
 between 
the set of irreducible components of the exceptional divisor of the minimal embedded resolution of $C$
and the rays of  $\cT_{reg}$ which intersect $\R_{>0}^{m+1}$, 
which is given by 
$E_{P'} \mapsto \R_{\geq 0} w^{P'}$, 
where $P'$ runs through the set of 
rational points of $\Theta_{L_0} (D)$ which are distinguished by $C$
(see Definition \ref{not:dis}). 
It follows that $\pi$ is the minimal embedded resolution of $C$.
This implies also 
that  $E_0 = E_P$. 
\end{proof}

%--------------------------------------------------------------------------------

\providecommand{\bysame}{\leavevmode\hbox to3em{\hrulefill}\thinspace}
\providecommand{\MR}{\relax\ifhmode\unskip\space\fi MR }

\providecommand{\href}[2]{#2}

\end{document}